\newcommand{\MailFile}[1]{}
\def\PaperDate{2011/01/30}
\documentclass[a4paper,12pt]{article}
  \usepackage{buxmath}
  \usepackage[DinA4]{buxlayout}
  \usepackage[quotes]{bux2ref}
  \usepackage{graphicx}
  \usepackage[all]{xy}
  \usepackage{datum}
  \usepackage[german,american]{babel}
  \usepackage{float}

  \renewcommand{\qedsymbol}{\ensuremath{\Box}}
  \newcommand{\ignore}[1]{}

  \usepackage{email}
  \MailFile{all_files.tex}
  \DoMail{makeall.make}
  \DoMail{hyperbolic_3_3_4.mp}
  \DoMail{w_one.mp}
  \DoMail{w_two.mp}
  \DoMail{almost_rich.mp}
  \DoMail{rich.mp}
  \DoMail{min_01.pov}
  \DoMail{min_02.pov}
  \DoMail{min_03.pov}
  \DoMail{move.mp}
  \DoNotMail{ams*, *bref, graphicx.sty, *cfg, *def, *fd}
  \DoNotMail{keyval.sty, trig.sty, xy*, babel*, *ldf}
  \DoNotMail{verbatim.sty}

  \DeclareRobustCommand{\change}[2][]{%
    \begingroup
      \def\empty{}%
      \def\arg{#2}%
      \ifx\empty\arg
        \def\next{\unskip}%
      \else
        \def\next{#2}%
      \fi
    \expandafter\endgroup\next
  }
  \let\xxx\iffalse
\xxx
  \let\oldlabel\label
  \def\label{%
    \begingroup
      \catcode`\_=11
      \dolabel
  }
  \def\dolabel#1{%
    \endgroup\oldlabel{#1}%
    \fbox{\texttt{#1}}%
  }
  \usepackage{color}
  \usepackage[normalem]{ulem}
  \usepackage{changebar}
  \DeclareRobustCommand{\delete}[1]{\begin{color}{blue}\sout{#1}\end{color}}
  \DeclareRobustCommand{\insert}[1]{\begin{color}{red}#1\end{color}}
  \DeclareRobustCommand{\swap}[2]{\delete{#1} \insert{#2}}
  \global\let\change\undefined
  \DeclareRobustCommand{\change}[2][]{%
    \cbstart
    \begingroup
      \def\next{}%
      \def\deleted{#1}%
      \def\newmat{#2}%
      \def\empty{}%
      \ifx\deleted\empty
        \def\next{\insert{#2}}%
      \else
        \ifx\newmat\empty
          \def\next{\delete{#1}}%
        \else
          \def\next{\swap{#1}{#2}}%
        \fi
      \fi
    \expandafter\endgroup\next
    \cbend
  {}}
\fi

  \InputIfFileExists{basic_notation.tex}{}{}

  \newcommand{\Period}{\text{.}}
  \newcommand{\Comma}{\text{,}}

  \newvariable{\ThePlaceSet}{S}
  \newvariable{\ThePlace}{p}
  \newvariable{\TheLattice}{\Gamma}
  \newvariable{\TheImage}{\Gamma_0}
  \newvariable{\TheImageElement}{\gamma}
  \newvariable{\TheLatticeElement}{\gamma}
  \newvariable{\TheGroupScheme}{\mathcal{G}}
  \newvariable{\TheGlobalField}{K}
  \newvariable{\TheLocalField}{\TheGlobalField}
  \newvariable{\TheRing}{R}
  \newvariable{\TheOkaRing}{\Oka_{\ThePlaceSet}}
  \newvariable{\TheEuBuild}{X}
  \newvariable{\TheEuSubset}{V}
  \newvariable{\TheEuPoint}{x}
  \newvariable{\AltEuPoint}{y}
  \newvariable{\TheEuChamber}{C}
  \newvariable{\TheStar}{B}

  \newvariable{\TheGroup}{G}
  \newvariable{\FinType}{m}
  \newvariable{\FinitenessLength}{\phi}

  \newvariable{\Sl}{\operatorname{SL}}
  \newvariable{\FinField}{\FFF}
  \newvariable{\ThePrimePower}{q}
  \newvariable{\TheIndeterminate}{t}

  \newvariable{\TheDim}{d}
  \newvariable{\LocalConnectivity}{m}
  \newvariable{\HomotopyGroup}{\pi}

  \newcommand{\CardOf}[1]{\left|#1\right|}
  \newvariable{\TheSize}{n}

  \newcommand{\CAT}[1]{{\small CAT}$(#1)$}
  \newvariable{\PosReals}{\RRR_{\geq\Zero}}
  \newvariable{\TheRay}{\rho}
  \newvariable{\TheTime}{t}
  \newvariable{\TheRatBuild}{\Delta}
  \newvariable{\Chambers}{\mathcal{C}}
  \newvariable{\TheRatChamber}{c}
  \newvariable{\AltRatChamber}{c'}
  \newvariable{\Vertices}{\mathcal{V}}
  \newvariable{\TheRatVertex}{v}
  \newvariable{\AltRatVertex}{w}
  \newvariable{\TheRatSimplex}{\tau}
  \newvariable{\AltRatSimplex}{\sigma}
  \newvariable{\TheInfBuild}{\TheRatBuild}
  \newvariable{\Boundary}{\partial}
  \newvariable{\SphSimplex}{\tau}
  \newvariable{\TheReals}{\RRR}
  \newvariable{\TheBusemann}{\beta}
  \newvariable{\TheBound}{t}
  \newvariable{\AltBound}{t'}
  \newvariable{\CloseSimplex}{\sigma}
  \newvariable{\TheDiam}{d}
  \newvariable{\TheRedFiltr}{Y}
  \newvariable{\Simplices}{\mathcal{S}}
  \newvariable{\TheBaseValue}{b}
  \newvariable{\TheLowerBound}{r}
  \newvariable{\TheUpperBound}{R}
  \newcommand{\compose}{\circ}
  \newvariable{\DummyVertex}{*}

  \newvariable{\TheMorse}{h}
  \newvariable{\RoughMorse}{\hat{\TheMorse}}
  \newvariable{\TheEuApp}{\Sigma}
  \newvariable{\AltEuApp}{\TheEuApp'}
  \newvariable{\TheEuVertex}{\TheEuPoint}
  \newvariable{\HorEuCell}{\tau}
  \newvariable{\HorEuFace}{\sigma}
  \newcommand{\faceof}{\leq}
  \newvariable{\ConvexChunk}{Y}
  \newvariable{\TheEnd}{e}
  \newvariable{\RoughEnd}{\hat{e}}
  \newcommand{\ConvexHullOf}[1]{\overline{#1}}
  \newcommand{\DistanceOf}[3][]{\operatorname{dist}_{#1}(#2,#3)}
  \newvariable{\TheIsometry}{\iota}

  \newcommand{\SubLevelSet}[2]{{#1}^{-1}([\Zero,#2])}
  \newvariable{\TheLevel}{s}
  \newvariable{\AltLevel}{\TheLevel'}
  \newvariable{\TheAngle}{\vartheta}
  \newcommand{\LineSegment}[2]{[#1,#2]}
  \newcommand{\RayFromTo}[2]{[#1,#2)}
  \newcommand{\RayFromThrough}[2]{\overrightarrow{#1,#2}}
  \newvariable{\GradOf}{\nabla}
  \newvariable{\Link}{\operatorname{lk}}
  \newcommand{\crossprod}{\times}

  \newvariable{\TheLastIndex}{n}
  \newvariable{\TheIndex}{i}

  \newcommand{\FType}[1][]{\ensuremath{\text{F}_{#1}}}

  \newvariable{\DummyArg}{-}

  \newvariable{\TheEuSpace}{\EEE}
  \newcommand{\InnerProd}[2]{\langle #1, #2 \rangle}
  \newvariable{\SpZero}{\mathbf{0}}
  \newvariable{\ThePolytope}{P}
  \newvariable{\TheFace}{F}
  \newvariable{\ThePolyVertex}{\mathbf{v}}
  \newvariable{\AltPolyVertex}{\mathbf{v'}}
  \newvariable{\TheNormalVect}{\mathbf{n}}
  \newvariable{\ThePolytopePoint}{\mathbf{z}}
  \newvariable{\AltPolytopePoint}{\mathbf{z'}}
  \newvariable{\TheSpPoint}{\mathbf{x}}
  \newvariable{\AltSpPoint}{\mathbf{y}}
  \newvariable{\ThrSpPoint}{\mathbf{z}}
  \newvariable{\TheSpVector}{\mathbf{x}}
  \newvariable{\AltSpVector}{\mathbf{x'}}
  \newvariable{\TheSpVertex}{\mathbf{v}}
  \newvariable{\AltSpVertex}{\mathbf{v'}}
  \newvariable{\TheSpSimplex}{\mathbf{\tau}}
  \newvariable{\ClosestPoint}{\operatorname{pr}}
  \newvariable{\NormalCone}{N}
  \newvariable{\TheZonoDef}{D}
  \newvariable{\TheZonotope}{Z}
  \newvariable{\TheZonoPoint}{\mathbf{z}}
  \newvariable{\TheSpEdge}{e}
  \newvariable{\TheZonoVect}{\mathbf{d}}
  \newvariable{\TheCoefficient}{a}
  \newvariable{\TheDistFct}{f}
  \newvariable{\MinSet}{\operatorname{min-set}}
  \newvariable{\MaxSet}{\operatorname{max-set}}
  \newvariable{\TheSpVertexSet}{\mathbf{V}}
  \newcommand{\NormOf}[1]{\|#1\|}
  \newvariable{\TheConvexMap}{f}
  \newvariable{\TheConvexSet}{C}

  \newvariable{\TheSpace}{M}
  \newvariable{\NewPart}{A}
  \newvariable{\BadSet}{B}
  \newvariable{\ClosedBadSet}{A}
  \newvariable{\TheRadialCenter}{x}
  \newvariable{\TheGeometricComplex}{L}
  \newvariable{\TheGeometricSubcomplex}{K}
  \newvariable{\AltGeometricComplex}{L'}
  \newvariable{\ThrGeometricComplex}{L''}
  \newvariable{\TheDeformation}{\rho}
  \newvariable{\GeometricCell}{\tau}
  \newvariable{\TheSphBuild}{\Delta}
  \newvariable{\TheSphChamber}{C}
  \newvariable{\TheSphApp}{\Sigma}
  \newvariable{\TheSphComplement}{Y}
  \global\let\Dim\undefined
  \newvariable{\Dim}{\operatorname{dim}}
  \newvariable{\TheSphere}{S}
  \newcommand{\ClosureOf}[1]{\overline{#1}}
  \newvariable{\NorthPole}{n}
  \newcommand{\ver}{\operatorname{ver}}
  \newcommand{\hor}{\operatorname{hor}}
  \newcommand{\join}{*}

  \newvariable{\TheModelApp}{\hat\Sigma}
  \newvariable{\TheFullWeyl}{W}
  \newvariable{\BigChunk}{Y}
  \newvariable{\TheDistBound}{{\RoughMorse}_{\TheZonoDef}}
  \newvariable{\TheZonoDiam}{R}
  \newcommand{\AngleOf}[3][]{\angle_{#1}(#2,#3)}
  \newcommand{\PiHalf}{\frac{\pi}{2}}
  \newvariable{\TheSegment}{\varepsilon}
  \newvariable{\TheChamberDiam}{\operatorname{chamber-diam}}
  \newvariable{\TheEuBrick}{B}
  \newvariable{\TheHugeBound}{{\TheUpperBound^*}}
  \newvariable{\TheHausdorffDist}{L}
  \newvariable{\TheShiftConst}{C}
  \newvariable{\MetricNbhd}{\operatorname{Nbhd}}

  \newcommand{\goesup}{\nearrow}
  \newcommand{\goesdown}{\searrow}
  \newcommand{\strictfaceof}{<}
  \newvariable{\Busemann}{\beta}
  \newvariable{\Depth}{\operatorname{dp}}
  \newcommand{\RoofOf}[2][]{\hat{#2}^{#1}}
  \newvariable{\TheEuCell}{\tau}
  \newvariable{\AltEuCell}{\sigma}
  \newvariable{\CombinedHeight}{f}
  \newvariable{\Cells}{\mathcal{C}}

  \newvariable{\TheEuSubdiv}{{\topcirc{\TheEuBuild}}}
  \newcommand{\CenterOf}[1]{\topcirc{#1}}

  \newvariable{\Star}{\operatorname{st}}
  \newvariable{\TheDescBoundary}{K}
  \newvariable{\TheAscVertex}{v}
  \newvariable{\DescLink}{\operatorname{lk}^{\downarrow}}
  \newvariable{\TheEuCoface}{\xi}
  \newvariable{\AltEuCoface}{\zeta}
  \newcommand{\strictcofaceof}{>}
  \newvariable{\TheHemisphereCx}{H}
  \newvariable{\hh}{\mathrm{h}}
  \newvariable{\vv}{\mathrm{v}}
  \newvariable{\UpSet}{L^{\uparrow}}
  \newvariable{\TheAppVertex}{v}
  \newcommand{\commoncell}{\vee}
  \newvariable{\TheHull}{A}
  \newvariable{\TheEpsilon}{\varepsilon}
  \newvariable{\TheHullProj}{\tilde{\TheHull}}
  \newvariable{\TheBoundaryChamber}{C}
  \newvariable{\TheRetraction}{\rho}
  \newvariable{\TheIsomorphism}{\TheRetraction|_{\AltEuApp}}
  \newvariable{\TheLinkChamber}{\tilde{\TheBoundaryChamber}}
  \newvariable{\TheLinkApp}{\tilde{\TheEuApp}}
  \newvariable{\AltLinkApp}{\tilde{\AltEuApp}}
  \newvariable{\TheLinkRetraction}{\tilde{\TheRetraction}}

  \newvariable{\TheFiltration}{F}
  \newvariable{\Image}{\operatorname{im}}

  \newvariable{\TheConstField}{k}
  \newvariable{\TheResField}{\TheConstField}
  \newvariable{\TheAdeles}{\AAA}
  \newvariable{\TheMaxIdeal}{\mathfrak{m}}
  \newvariable{\TheOka}{\mathcal{O}}
  \newvariable{\TheLocalDeg}{d}
  \newcommand{\ExtDegreeOf}[2]{[#1:#2]}
  \newvariable{\TheFieldElement}{f}
  \newcommand{\AbsValueOf}[2][]{\|#2\|_{#1}}
  \newvariable{\TheAdel}{\mathbf{f}}
  \newcommand{\LogOf}[1]{\log_{\CardOf{\TheConstField}}(#1)}
  \newvariable{\ExtGlobalField}{\TheGlobalField'}
  \newvariable{\ExtAdeles}{\TheAdeles'}
  \newvariable{\ExtPlace}{\ThePlace'}
  \newvariable{\ExtNorm}{N}
  \newvariable{\ExtFieldElement}{\TheFieldElement'}
  \newvariable{\TheCharacter}{\chi}
  \newcommand{\CharKernel}[1]{#1^{\circ}}

  \newvariable{\TheSpecVertexSet}{\mathfrak{X}}
  \newvariable{\TheBaseVertex}{x}
  \newvariable{\TheParabolic}{\mathcal{P}}
  \newvariable{\AltParabolic}{\mathcal{P}'}
  \newvariable{\TheLevy}{\mathcal{L}}
  \newvariable{\TheTorus}{\mathcal{T}}
  \newvariable{\UniRadical}{\mathcal{R}_{\mathrm{u}}}
  \newvariable{\TheRadical}{\mathcal{R}}
  \newvariable{\TheSpecVertex}{\mathfrak{x}}
  \global\let\Stab\undefined
  \newvariable{\Stab}{\operatorname{Stab}}
  \newvariable{\TheCompactSubgroup}{\mathcal{C}}
  \newvariable{\TheHarderForm}{\omega}
  \newvariable{\TheHarderMeasure}{\operatorname{d}{\TheHarderForm}_{\TheAdeles}}
  \newvariable{\Volume}{\operatorname{vol}}
  \newvariable{\HarderPi}{\pi}
  \newvariable{\TheGroupElement}{\gamma}
  \newcommand{\Conj}[2]{{}^{#1}#2}
  \newvariable{\AlmostBusemann}{\tilde\beta}
  \newvariable{\TheHarderChar}{\chi}
  \newvariable{\Mult}{\mathfrak{Mult}}
  \newvariable{\TheAdeleGroupElement}{\mathbf{\gamma}}
  \newvariable{\TheSpecVertexApp}{\mathfrak{S}}
  \newvariable{\TheScaleFactor}{c}
  \newvariable{\Characters}{\operatorname{X}}
  \newvariable{\TheSplitTorus}{\TheTorus'}
  \newvariable{\TheRationals}{\QQQ}
  \newvariable{\TheRoot}{\alpha}
  \newvariable{\TheLastRoot}{r}
  \newvariable{\TheRootIndex}{i}
  \newcommand{\subgroup}{\leq}
  \global\let\Hom\undefined
  \newvariable{\Hom}{\operatorname{Hom}}
  \newvariable{\TheType}{j}
  \newvariable{\AltType}{k}
  \newvariable{\TheCotype}{i}
  \newvariable{\AltCotype}{j}
  \newvariable{\TheBaseChange}{c}
  \newvariable{\InvBaseChange}{n}
  \newvariable{\HarderNu}{\nu}
  \newvariable{\LogNu}{\mu}
  \newvariable{\LogPi}{\beta}
  \newvariable{\HarderConstOne}{C_1}
  \newvariable{\HarderConstTwo}{C_2}
  \newvariable{\TheLogNu}{\mu}
  \newcommand{\TypeOf}[1]{\operatorname{type}(#1)}
  \newvariable{\LogConst}{R}
  \newvariable{\TheRootBound}{c}
  \newvariable{\AdelesTorusElement}{\mathbf{t}}
  \newvariable{\TheLeftCompact}{L}
  \newvariable{\TheRightCompact}{R}
  \newvariable{\TheMahlerSet}{M}
  \newvariable{\TheMahlerElement}{\TheAdeleGroupElement}
  \newvariable{\TheLowerCompactnessBound}{c_{-}}
  \newvariable{\TheUpperCompactnessBound}{c_{+}}
  \newvariable{\TheStandardSubgroup}{\mathcal{K}}
  \newvariable{\TheBorel}{\mathcal{B}}
  \newvariable{\AltIndex}{j}

  \newvariable{\TheAdmissibleSet}{A}
  \newvariable{\TheBase}{\alpha}
  \newvariable{\DualBase}{\chi}
  \newvariable{\Const}{c}
  \newvariable{\TheLogNuSet}{\mathfrak{M}}
  \newvariable{\TheLogPiSet}{\mathfrak{B}}

  \newvariable{\acts}{\cdot}
  \newvariable{\TheBasePolyVertex}{*}
  \newcommand{\fieldextends}{/}
  \newcommand{\sep}{\operatorname{s}}
  \newvariable{\ExtSplitTorus}{\mathcal{T}^{\sep}}
  \newvariable{\ExtVectorSpace}{V^{\sep}}
  \newvariable{\TheVectorSpace}{V}
  \newvariable{\Gal}{\operatorname{Gal}}
  \newcommand{\subspace}{\leq}
  \global\let\Aut\undefined
  \newvariable{\Aut}{\operatorname{Aut}}
  \newvariable{\TheHoroCenter}{e}
  \newvariable{\TheSpeed}{s}
  \newvariable{\StdEuApp}{\TheEuApp}
  \newvariable{\AltLowerBound}{\TheLowerBound'}
  \newvariable{\ThrLowerBound}{\TheLowerBound''}
  \newvariable{\AltUpperBound}{\TheUpperBound'}
  \newvariable{\TheCone}{H}
  \newvariable{\TheOrigin}{0}
  \newcommand{\isometric}{\cong}
  \newvariable{\DualBusemann}{\mu}
  \newvariable{\DualChunk}{Z}
  \newvariable{\NormalChunk}{N}

  \newvariable{\TheLastPoint}{s}
  \newvariable{\TheSector}{S}
  \newvariable{\TheFundamentalDomain}{D}
  \newvariable{\TheTips}{Q}
  \newvariable{\TheIndexBound}{K}
  \newvariable{\TheBall}{B}
  \newvariable{\TheTree}{T}
  \newvariable{\TheRepSet}{\mathcal{R}}
  \newvariable{\TheChRepSet}{\mathcal{X}}
  \newvariable{\TheDistance}{d}

\begin{document}
  \title{Higher Finiteness Properties of Reductive Arithmetic Groups in
         Positive Characteristic:\\the Rank Theorem}
  \author{Kai-Uwe Bux \and
          Ralf Gramlich \and
          Stefan Witzel}
  \date{\datum\PaperDate}
  \maketitle
  \begin{abstract}
    \noindent
    We show that the finiteness length of an $\ThePlaceSet$-arithmetic
    subgroup $\TheLattice$ in a non-commutative isotropic absolutely
    almost simple group $\TheGroupScheme$ over
    a global function field is one less than the sum of the local ranks of
    $\TheGroupScheme$ taken over the places in $\ThePlaceSet$. This
    determines the finiteness properties for arithmetic subgroups
    in isotropic reductive groups,
    confirming the conjectured finiteness properties for this class
    of groups.

    Our main tool is Behr-Harder reduction theory which we recast
    in terms of the metric structure of euclidean buildings.
  \end{abstract}

  \noindent
  Let $\TheGlobalField$ be a global function field and let
  $\TheGroupScheme$ be a linear algebraic group defined
  over $\TheGlobalField$.   We fix a finite set $\ThePlaceSet$ of
  places over $\TheGlobalField$ and let $\TheOkaRing$ denote the
  subring of $\ThePlaceSet$-integers in $\TheGlobalField$. We want
  to talk about the group $\TheLattice:=\TheGroupSchemeOf{\TheOkaRing}$,
  but as an algebraic variety, the $\TheRing$-points of $\TheGroupScheme$
  are only well-defined for $\TheGlobalField$-algebras $\TheRing$.
  However, we regard $\TheGroupScheme$ as a concrete matrix group
  defined by polynomial equations in the matrix coefficients. That is,
  we choose a particular realization of the variety $\TheGroupScheme$
  as an algebraic set in some affine space. Given this realization,
  we define $\TheLattice$ as its set of $\TheOkaRing$-points.
  The subgroup $\TheLattice$ obtained this way is called an
  \notion{$\ThePlaceSet$-arithmetic subgroup} of $\TheGroupScheme$.
  Of course, the arithmetic group $\TheLattice$ depends on the chosen
  realization of $\TheGroupScheme$, but any two choices lead to
  $\ThePlaceSet$-arithmetic subgroups of $\TheGroupScheme$ that share a
  common subgroup of finite index in both; see, e.g., \cite[\S 1]{Serre:1979}.
  Hence, the commensurability  class of $\TheLattice$ depends only
  on the group scheme $\TheGroupScheme$ and the $\ThePlaceSet$-arithmetic
  ring $\TheOkaRing$.

  We are interested in finiteness properties of the group
  $\TheLattice$. Recall that a group $\TheGroup$ is of
  \notion{type \FType[\FinType]} if it admits
  a classifying space
  with finite $\FinType$-skeleton. The \notion{finiteness
  length} $\FinitenessLengthOf{\TheGroup}$
  of $\TheGroup$ is the largest $\FinType$ such that
  $\TheGroup$ is of type \FType[\FinType]. We say that
  $\FinitenessLengthOf{\TheGroup}= \infty$ if
  $\TheGroup$ is of type \FType[\FinType] for all $\FinType$.
  The finiteness length is a commensurability invariant:
  in fact, it is invariant under
  quasi-isometries \cite{Alonso:1994}.
  In particular, the finiteness length
  $\FinitenessLengthOf{\TheLattice}$ of the
  $\ThePlaceSet$-arithmetic group $\TheLattice$
  depends only on $\TheGroupScheme$ and $\TheOkaRing$ but not
  on the particular chosen realization of $\TheGroupScheme$ as
  a matrix group.

  Let $\TheLocalField[\ThePlace]$ be the completion of the field
  $\TheGlobalField$ at the place $\ThePlace$. The \notion{local
  rank} of $\TheGroupScheme$ at the place $\ThePlace$ is the
  rank of $\TheGroupScheme$ over the field
  $\TheLocalField[\ThePlace]$. If $\TheGroupScheme$ is isotropic
  and absolutely almost simple, the group
  $\TheGroupSchemeOf{\TheLocalField[\ThePlace]}$ acts on its
  associated Bruhat-Tits building $\TheEuBuild[\ThePlace]$
  and the dimension $\DimOf{\TheEuBuild[\ThePlace]}$ is the
  local rank of $\TheGroupScheme$ at the place
  $\ThePlace$. We prove the following theorem, which answers
  \cite[Question~13.20]{Abramenko.Brown:2008}.
  \begin{NewTh}<statement>{Rank Theorem}
    Let $\TheGroupScheme$ be a connected non-commutative absolutely
    almost simple $\TheGlobalField$-isotropic $\TheGlobalField$-group.
    Then the finiteness length
    \(
      \FinitenessLengthOf{
        \TheLattice
      }
    \)
    of the $\ThePlaceSet$-arithmetic group
    $\TheLattice=\TheGroupSchemeOf{\TheOkaRing}$ is
    $\TheDim-\One$ where $\TheDim:=
    \Sum[\ThePlace\in\ThePlaceSet]{\DimOf{\TheEuBuild[\ThePlace]}}$
    is the sum of the local ranks of $\TheGroupScheme$.
  \end{NewTh}
  It was shown in \cite{Bux.Wortman:2007} that
  \(
    \FinitenessLengthOf{
      \TheLattice
    }
    \leq
    \TheDim-\One.
  \)
  Hence, only the positive statement is new:
  $\TheLattice$ is of type \FType[\TheDim-\One].

  If $\TheGroupScheme$ is reductive and anisotropic,
  J.-P.\,Serre showed that $\ThePlaceSet$-arithmetic
  subgroups are of type \FType[\infty]. More precisely,
  $\TheGroupSchemeOf{\TheOkaRing}$ has a torsion free subgroup
  of finite index that admits a finite Eilenberg-Mac\,Lane
  complex \cite[Cas~(b), p.~126--127]{Serre:1971}.

  The Rank Theorem contrasts with the number field case where
  $\ThePlaceSet$-arithmetic subgroups of reductive groups are
  of type \FType[\infty] \cite[\S\,11]{Borel.Serre:1976}. We refer
  to the introduction of \cite{Bux.Wortman:2007} for some
  conjectures about a more quantitative account which should
  reveal indeed deep similarities in the geometric underpinnings of
  both cases.

  Interest in finiteness properties of $\TheLattice$ started
  in 1959 when H.\,Nagao \cite{Nagao:1959} showed that
  $\SlOf[\Two]{\FinFieldAd[\ThePrimePower]{\TheIndeterminate}}$ is
  not finitely generated. In this case, there is a single
  place and the corresponding euclidean building is a tree,
  thus $\TheDim=\One$.

  In 1969, H.\,Behr \cite{Behr:1969} proved that $\TheLattice$ is
  finitely generated if and only if $\TheDim>1$. He had to exclude
  a few cases. However, as he pointed out, those restrictions can be
  removed by appealing to Harder's version \cite{Harder:1969} of
  reduction theory.
  Using Harder's reduction theory again, Behr \cite{Behr:1998} showed in
  1998 that $\TheLattice$ is finitely presented if and only if
  $\TheDim>2$.

  Concerning higher finiteness properties (i.e., beyond finite
  presentability), U.\,Stuhler \cite{Stuhler:1980} showed that
  $\SlOf[\Two]{\TheOkaRing}$ has finiteness length
  $\CardOf{\ThePlaceSet}-1=\TheDim-\One$. H.\,Abels
  \cite{Abels:1991} and P.\,Abramenko \cite{Abramenko:1987} showed
  that $\SlOf[\TheSize]{\FinFieldAd[\ThePrimePower]{\TheIndeterminate}}$
  have finiteness length $\TheSize-\Two=\TheDim-\One$ provided
  that $\ThePrimePower$ is large enough depending on $\TheSize$.
  Abramenko \cite{Abramenko:1996} extended this result to classical
  groups by recasting it in the context of groups acting
  on twin buildings. The need to exclude small $\ThePrimePower$
  arises from the method of proof: it involves the analysis of
  certain subcomplexes in spherical buildings which only have
  sufficient topological connectivity if the underlying buildings
  are sufficiently thick. The articles
  \cite{Devillers.Gramlich.Muehlherr:2009}
  and
  \cite{Gramlich.Witzel:2009}
  also suffer from this shortcoming as they use a filtration
  modeled upon the filtration by combinatorial codistance
  introduced by Abels and Abramenko. Hence, similar relative
  links arise.

  In 2005, Bernd~Schulz analyzed in his PhD thesis
  \cite{Schulz:2010}  a class of
  subcomplexes of spherical buildings which have the right
  topological connectivity without restrictions. Immediately following
  his discovery, several results were obtained in the positive
  direction (i.e., establishing that $\TheLattice$ is of
  type \FType[\TheDim-\One]).
  In \cite{Bux.Wortman:2008}, K.\,Wortman
  and the first author have proved the Rank~Theorem for $\TheGlobalField$-groups
  of $\TheGlobalField$-rank one using reduction theory as a source
  for a filtration that leads to relative links analyzed by Schulz.
  In a previous version of this
  paper \cite{Bux.Gramlich.Witzel:2009},
  the authors eliminated the restriction on $\ThePrimePower$
  in Abramenko's results. The third author \cite{Witzel:2010}
  extended the analysis to
  \(
    \TheGroupSchemeOf{\FinFieldAd[\ThePrimePower]{
      \TheIndeterminate,\TheIndeterminate[][-\One]
    }}
  \).
  The basic idea was to replace combinatorial codistance by a
  metric codistance that leads to better behaved relative links.
  This paper extends the reduction theory approach
  from \cite{Bux.Wortman:2008} and removes the restriction on
  the global rank.

  The Rank~Theorem allows one to deduce finiteness properties of
  arbitrary reductive groups as described
  in \cite[2.6(c), page~91]{Behr:1998}:
  First pass to the connected
  component of the identity. For the arithmetic subgroup this means
  passing to a subgroup of finite index which does not change the
  finiteness length. Reducing to a semi-simple group scheme by
  splitting off a central torus does also not affect the finiteness
  length of $\TheLattice$. Using a central
  isogeny \cite[Satz~1]{Behr:1968}, we can
  assume that the group scheme is simply connected still without
  any change of the finiteness length. A simply connected
  semi-simple group is
  the direct product of its almost simple factors (which remain
  simply connected). The finiteness length of a direct product
  is the minimum of the finiteness lengths of its factors. Finally,
  by restriction of scalars, one may assume that each factor of the
  direct product
  is absolutely almost simple (and still simply connected); see, e.g.,
  \cite[Hilfssatz~7.4 and~7.5]{Kneser:1965}. Kneser treats the number
  field case, but his arguments can be extended without difficulty.
  Finally, the Rank~Theorem applies to the absolut almost simple factors.

  \paragraph{Acknowledgments.}
    The authors wish to thank the
    Hausdorff Research Institute for Mathematics (Bonn) for their
    hospitality within the trimester programs ''Rigidity''
    and ``Algebra and number theory''. We also gratefully
    acknowledge the hospitality of the
    Mathematisches Forschungsinstitut Oberwolfach
    during several workshops and their RiP program.
    The {\small DFG} provided appreciated financial and structural
    support via {\small SFB 701} at Bielefeld University
    and via {\small GR~2077/5} and {\small GR~2077/7} at Darmstadt.

    We are also indebted to
    Bernd~Schulz for explaining his thesis.
    The participants of the Bielefeld seminar on reduction theory
    in summer~2010
    (in particular Herbert~Abels, Werner~Hoffmann, Gregory~Margulis, and
    Andrei~Rapinchuk) provided valuable expertise.
    We thank
    Herbert~Abels,
    Peter~Abramenko,
    Richard~Weiss,
    Bernhard~M\"uhlherr,
    Michael~Joswig,
    Sven~Herrmann,
    Bernd~Schulz, and
    Hendrik~Van~Maldeghem
    for helpful discussions.
    Andrei~Rapinchuk has helped us with Lemma~\ref{harder:2.2.3}.
    We also thank Bertrand~R\'e{}my for sharing with us his
    knowledge on Kac-Moody groups, algebraic groups, and their
    interaction.

  \section{Behr-Harder reduction theory}\label{sec:behr}
    In this section, we collect the results of reduction theory.
    Let $\TheGroupScheme$ be a connected, reductive, non-commutative
    $\TheGlobalField$-isotropic $\TheGlobalField$-group.
    In particular, it has proper $\TheGlobalField$-parabolic subgroups.
    In Section~\ref{sec:us}, we shall state and prove results
    in the slightly more general setting of reductive but not
    necessarily isotropic $\TheGroupScheme$.

    The euclidean building
    \(
      \TheEuBuild := \Product[\ThePlace\in\ThePlaceSet]{
        \TheEuBuild[\ThePlace]
      }
    \)
    associated to $\TheGroupSchemeOf{\TheOkaRing}$
    is a \CAT{\Zero}-space
    \cite[Theorem~11.16]{Abramenko.Brown:2008}.
    Its \notion{visual boundary} consists
    of parallelism classes of geodesic rays
    \(
      \TheRay\mapcolon\PosReals \rightarrow \TheEuBuild
    \).
    To each such ray, one associates a \notion{Busemann function}
    \begin{eqnarray*}
      \TheBusemann\mapcolon\TheEuBuild &\longrightarrow&\TheReals\\
      \TheEuPoint & \mapsto &
      \Limit[\TheTime\rightarrow\infty]{(
        \DistanceOf{\TheRayOf{\Zero}}{\TheRayOf{\TheTime}}
        -
        \DistanceOf{\TheEuPoint}{\TheRayOf{\TheTime}}
      )}\Period
    \end{eqnarray*}
    Two rays are parallel if and only if their corresponding Busemann
    functions differ by an additive constant. In particular, a Busemann
    function determines a unique point in the visual boundary, its
    \notion{center}. See, e.g.,
    \cite[pages~267ff, in particular~8.20]{Bridson.Haefliger:1999}.

    Let $\TheRatBuild:=\TheRatBuild[\TheGlobalField]$ be the spherical building
    for the group $\TheGroupSchemeOf{\TheGlobalField}$. It is an irreducible
    building whose chambers correspond to the minimal $\TheGlobalField$-parabolic
    subgroups and whose vertices correspond to the
    maximal $\TheGlobalField$-parabolic subgroups of $\TheGroupScheme$.
    Let $\ChambersOf{\TheRatBuild}$ denote its set of chambers and
    $\VerticesOf{\TheRatBuild}$ its set of vertices.


    The visual boundary $\BoundaryOf{\TheEuBuild}$ is the spherical join
    of the boundaries $\BoundaryOf{\TheEuBuild[\ThePlace]}$ where the join is
    taken over all $\ThePlace\in\ThePlaceSet$.
    In Proposition~\ref{us:metric}, we describe how one can construct
    a $\TheLattice$-invariant isometric embedding
    \(
      \TheRatBuild\monorightarrow\BoundaryOf{\TheEuBuild}
    \).
    Using this embedding, we identify vertices of $\TheRatBuild$
    with points in the visual boundary $\BoundaryOf{\TheEuBuild}$.
    It turns out that Busemann functions centered at vertices
    of $\TheRatBuild\subset\BoundaryOf{\TheEuBuild}$ are
    not constant on any of the factors $\TheEuBuild[\ThePlace]$
    of $\TheEuBuild$.

    \begin{figure}[t]
      \null\hfill
      \setlength{\unitlength}{1cm}
      \begin{picture}(10,8)
        \put(-5,-0.2){\includegraphics{reduction_01.ps}}
        \put(4.4,1){$\ConvexChunkOf[\TheEuApp,\TheRatSimplex]{\TheBound}$}
      \end{picture}
      \hfill\null
      \caption{the closest point projection $\ClosestPoint[\TheEuApp,\TheRatSimplex][\TheBound]$\label{behr:01.fig}}
    \end{figure}

    \begin{figure}[t]
      \null\hfill
      \setlength{\unitlength}{1cm}
      \begin{picture}(10,8)
        \put(-5,-0.2){\includegraphics{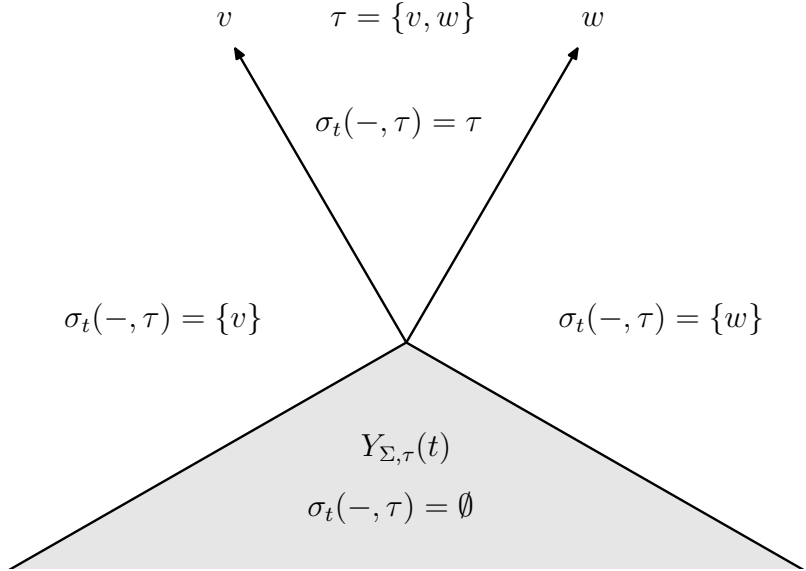}}
        \put(2.5,7){$\TheRatVertex$}
        \put(7.3,7){$\AltRatVertex$}
        \put(4.0,7){$\TheRatSimplex=\SetOf{\TheRatVertex,\AltRatVertex}$}
        \put(4.4,1.3){$\ConvexChunkOf[\TheEuApp,\TheRatSimplex]{\TheBound}$}
        \put(0.5,3){$\CloseSimplexOf[\TheBound]{\DummyArg,\TheRatSimplex}=
                  \SetOf{\TheRatVertex}$}
        \put(7,3){$\CloseSimplexOf[\TheBound]{\DummyArg,\TheRatSimplex}=
                  \SetOf{\AltRatVertex}$}
        \put(3.8,5.6){$\CloseSimplexOf[\TheBound]{\DummyArg,\TheRatSimplex}=
                  \TheRatSimplex$}
        \put(3.7,0.5){$\CloseSimplexOf[\TheBound]{\DummyArg,\TheRatSimplex}=
                  \emptyset$}
      \end{picture}
      \hfill\null
      \caption{decomposition of an apartment\label{behr:03.fig}}
      {\small
      Shown is an apartment
      whose visual boundary contains the simplex $\TheRatSimplex$ and the
      decomposition of the apartment into regions according to the values of
      $\CloseSimplexOf[\TheBound]{\DummyArg,\TheRatSimplex}$.
      Note in particular how moving up the tip will affect the picture.}
    \end{figure}

    Now, fix a family
    \(
      \FamOf[{
        \TheRatVertex \in \VerticesOf{\TheRatBuild}
      }]{
        \TheBusemann[\TheRatVertex] \mapcolon
        \TheEuBuild\rightarrow \TheReals
      }
    \)
    of Busemann functions so that each $\TheBusemann[\TheRatVertex]$
    is centered at $\TheRatVertex$. For any simplex $\TheRatSimplex$
    of $\TheRatBuild$, put
    \(
      \TheBusemannOf[\TheRatSimplex]{\TheEuPoint}
      :=
      \MaxOf[\TheRatVertex\in\TheRatSimplex]{
        \TheBusemannOf[\TheRatVertex]{\TheEuPoint}
      }.
    \)
    For an apartment $\TheEuApp$ of $\TheEuBuild$ and a simplex
    $\TheRatSimplex$ of $\TheRatBuild$ contained in the visual
    boundary $\BoundaryOf{\TheEuApp}$, we consider the convex
    cones
    \[
      \ConvexChunkOf[\TheEuApp,\TheRatSimplex]{\TheBound}
      :=
      \SetOf[\TheEuPoint\in\TheEuApp]{
        \TheBusemannOf[\TheRatSimplex]{\TheEuPoint}\leq\TheBound
      }
    \]
    as dependent on a real parameter $\TheBound$. Let
    \[
      \ClosestPoint[\TheEuApp,\TheRatSimplex][\TheBound]
      \mapcolon
      \TheEuApp
      \longrightarrow
      \ConvexChunkOf[\TheEuApp,\TheRatSimplex]{\TheBound}
    \]
    denote the closest point projection.
    \begin{observation}\label{behr:eliminate.apartment}
      Any two apartments $\TheEuApp$ and $\AltEuApp$ containing
      $\TheEuPoint\in\TheEuBuild$ and with
      $\TheRatSimplex$ in their visual boundary can be identified
      via an isometry that commutes with the Busemann functions
      $\TheBusemann[\TheRatVertex]$ for
      all $\TheRatVertex\in\TheRatSimplex$. More precisely,
      there is an isomorphism of Coxeter complexes
      \[
        \TheIsometry \mapcolon
        \TheEuApp \longrightarrow \AltEuApp
      \]
      such that the following diagram commutes:
      \[
        \xymatrix{
          {\TheEuApp}
          \ar[rr]^{\TheIsometry}
          \ar[dr]_{\TheBusemann[\TheRatVertex]}
          &
          &
          {\AltEuApp}
          \ar[dl]^{\TheBusemann[\TheRatVertex]}
          \\
          &
          {\TheReals}
          &
        }
      \]
      In particular, $\TheIsometry$ identifies
      $\ConvexChunkOf[\TheEuApp,\TheRatSimplex]{\TheBound}$
      and
      $\ConvexChunkOf[\AltEuApp,\TheRatSimplex]{\TheBound}$.
      Moreover,
      \[
        \ClosestPoint[\AltEuApp,\TheRatSimplex][\TheBound]
        \compose
        \TheIsometry
        =
        \TheIsometry
        \compose
        \ClosestPoint[\TheEuApp,\TheRatSimplex][\TheBound]
      \]
      and the values
      \(
        \TheBaseValueOf[\TheRatSimplex,\TheRatVertex][\TheBound]{
          \TheEuPoint
        }
        :=
        \TheBusemannOf[\TheRatVertex]{
          \ClosestPointOf[\TheEuApp,\TheRatSimplex][\TheBound]{
            \TheEuPoint
          }
        }
      \)
      are independent of the apartment $\TheEuApp$.\qed
    \end{observation}
    We put
    \[
      \CloseSimplexOf[\TheBound]{\TheEuPoint,\TheRatSimplex}
      :=
      \SetOf[\TheRatVertex\in\TheRatSimplex]{
        \TheBaseValueOf[\TheRatSimplex,\TheRatVertex][\TheBound]{
          \TheEuPoint
        }
        =
        \TheBound
      }\Period
    \]
    Thinking within a given apartment $\TheEuApp$ containing
    $\TheEuPoint$ and $\TheRatSimplex$, the set
    \(
      \CloseSimplexOf[\TheBound]{\TheEuPoint,\TheRatSimplex}
    \)
    collects precisely those vertices $\TheRatVertex\in\TheRatSimplex$
    whose associated inequalities
    \(
      \TheBusemannOf[\TheRatVertex]{\DummyArg}
      \leq
      \TheBound
    \)
    are sharp at the point
    \(
      \ClosestPointOf[\TheEuApp,\TheRatSimplex][\TheBound]{
        \TheEuPoint
      }
    \).
    Hence, we may delete the other inequalities:
    \begin{observation}\label{behr:restricting.simplex}
      For any subsimplex $\AltRatSimplex\subseteq\TheRatSimplex$
      containing $\CloseSimplexOf[\TheBound]{\TheEuPoint,\TheRatSimplex}$,
      the closest point to $\TheEuPoint$ in
      \(
        \ConvexChunkOf[\TheEuApp,\TheRatSimplex]{\TheBound}
      \)
      is also the closest point to $\TheEuPoint$ in
      \(
        \ConvexChunkOf[\TheEuApp,\AltRatSimplex]{\TheBound}
      \),
      i.e.,
      \(
        \ClosestPointOf[\TheEuApp,\TheRatSimplex][\TheBound]{\TheEuPoint}
        =
        \ClosestPointOf[\TheEuApp,\AltRatSimplex][\TheBound]{\TheEuPoint}
      \).
      In particular, it follows that
      \(
        \CloseSimplexOf[\TheBound]{\TheEuPoint,\TheRatSimplex}
        =
        \CloseSimplexOf[\TheBound]{\TheEuPoint,\AltRatSimplex}
      \).\qed
    \end{observation}

    We say that a chamber $\TheRatChamber\in\ChambersOf{\TheRatBuild}$
    \notion{$\TheBound$-reduces}
    $\TheEuPoint\in\TheEuBuild$
    if
    \(
      \CloseSimplexOf[\TheBound]{\TheEuPoint,\TheRatChamber}
      =
      \TheRatChamber
    \).
    \begin{figure}[ht]
      \null\hfill
      \setlength{\unitlength}{1cm}
      \begin{picture}(10,8)
        \put(0,-0.2){\includegraphics{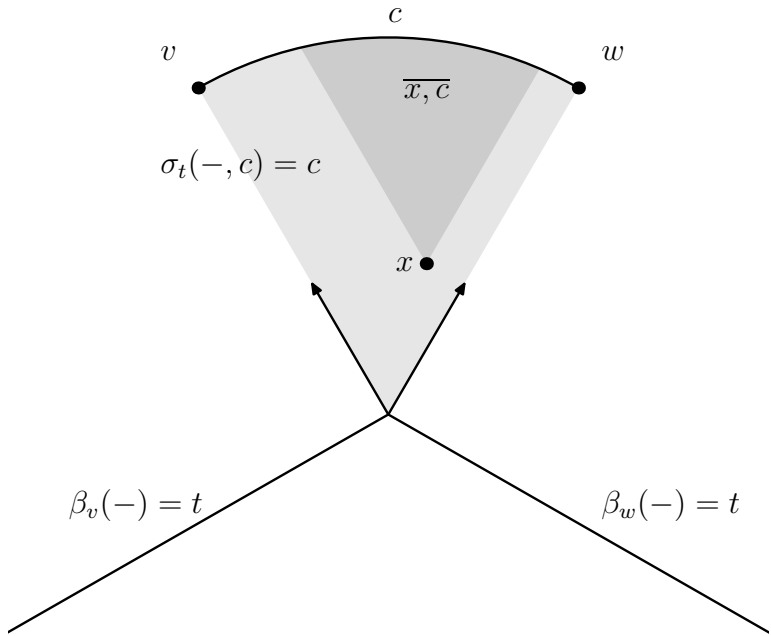}}
        \put(2,7.5){$\TheRatVertex$}
        \put(0.8,1.5){$\TheBusemannOf[\TheRatVertex]{\DummyArg}=\TheBound$}
        \put(7.8,7.5){$\AltRatVertex$}
        \put(7.8,1.5){$\TheBusemannOf[\AltRatVertex]{\DummyArg}=\TheBound$}
        \put(5.1,4.7){$\TheEuPoint$}
        \put(2,6){$\CloseSimplexOf[\TheBound]{\DummyArg,\TheRatChamber}=\TheRatChamber$}
        \put(5.2,7){$\ConvexHullOf{\TheEuPoint,\TheRatChamber}$}
        \put(5,8){$\TheRatChamber$}
      \end{picture}
      \hfill\null
      \caption{reducing convex hulls\label{behr:02.fig}}
      {\small
      The lightly shaded area is the region of points in the
      depicted apartment reduced by $\TheRatChamber$. The darker region
      is the cone $\ConvexHullOf{\TheEuPoint,\TheRatChamber}$.}%
    \end{figure}

    \begin{observation}\label{behr:close.nonempty}
      For $\TheBound\leq\AltBound$,
      \[
        \CloseSimplexOf[\AltBound]{\TheEuPoint,\TheRatSimplex}
        \subseteq
        \CloseSimplexOf[\TheBound]{\TheEuPoint,\TheRatSimplex}
        \subseteq
        \TheRatSimplex
        \Period
        \qed
      \]
    \end{observation}
    \begin{cor}\label{behr:independent.max}
      Assume $\AltBound:=\TheBusemannOf[\TheRatSimplex]{\TheEuPoint}
      \geq\TheBound$. Then each vertex $\TheRatVertex$ with
      $\TheBusemannOf[\TheRatVertex]{\TheEuPoint}=\AltBound$
      belongs to $\CloseSimplexOf[\TheBound]{\TheEuPoint,\TheRatSimplex}$.
      In particular,
      \(
        \CloseSimplexOf[\TheBound]{\TheEuPoint,\TheRatSimplex}
        \neq\emptyset
      \)
      and
      \(
        \TheBusemannOf[\TheRatSimplex]{\TheEuPoint}
        =
        \TheBusemannOf[{
          \CloseSimplexOf[\TheBound]{\TheEuPoint,\TheRatSimplex}
        }]{
          \TheEuPoint
        }
      \).
    \end{cor}
    \begin{proof}
      We have
      \(
        \SetOf[\TheRatVertex\in\TheRatSimplex]{
          \TheBusemannOf[\TheRatVertex]{\TheEuPoint}=\AltBound
        }
        =
        \CloseSimplexOf[\AltBound]{\TheEuPoint,\TheRatSimplex}
        \subseteq
        \CloseSimplexOf[\TheBound]{\TheEuPoint,\TheRatSimplex}
      \).
    \end{proof}
    \begin{observation}\label{behr:strahlensatz}
      Assume that $\TheRatChamber$ $\TheBound$-reduces
      $\TheEuPoint$. As illustrated in Figure~\ref{behr:02.fig},
      the chamber
      $\TheRatChamber$ $\TheBound$-reduces
      every point in the
      sector $\ConvexHullOf{\TheEuPoint,\TheRatChamber}$.\qed
    \end{observation}

      A \notion{reduction datum} consists of a family
      \(
        \FamOf[{
          \TheRatVertex \in \VerticesOf{\TheRatBuild}
        }]{
          \TheBusemann[\TheRatVertex] \mapcolon \TheEuBuild\rightarrow \TheReals
        }
      \)
      of Busemann functions on the euclidean building $\TheEuBuild$
      and two constants $\TheLowerBound<\TheUpperBound$ so that
      the following holds:
      \begin{quote}
          For any chamber $\TheRatChamber$ that $\TheLowerBound$-reduces
          $\TheEuPoint$, the simplex
          $\CloseSimplexOf[\TheUpperBound]{\TheEuPoint,\TheRatChamber}$
          is contained in \emph{any} chamber $\AltRatChamber$ that
          $\TheLowerBound$-reduces $\TheEuPoint$.
      \end{quote}
    \begin{observation}\label{behr:enlarge}
      If
      \(
        \TupelOf{
          \FamOf{\TheBusemann[\DummyVertex]},
          \TheLowerBound,
          \TheUpperBound[\Zero]
        }
      \)
      is a reduction datum and $\TheUpperBound[\One]>\TheUpperBound[\Zero]$,
      then so is
      \(
        \TupelOf{
          \FamOf{\TheBusemann[\DummyVertex]},
          \TheLowerBound,
          \TheUpperBound[\One]
        }
      \).\qed
    \end{observation}

    \begin{figure}[ht]
      \null\hfill
      \setlength{\unitlength}{1cm}
      \begin{picture}(9,7)
        \put(0,0){\includegraphics{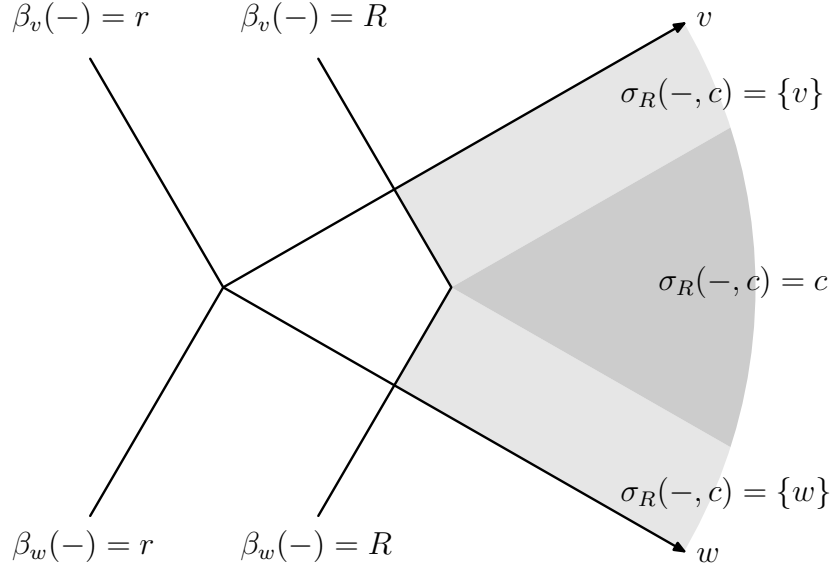}}
        \put(8,7){$\TheRatVertex$}
        \put(8,-0.1){$\AltRatVertex$}
        \put(-1,7){$\TheBusemannOf[\TheRatVertex]{\DummyArg}=\TheLowerBound$}
        \put(-1,0){$\TheBusemannOf[\AltRatVertex]{\DummyArg}=\TheLowerBound$}
        \put(2,7){$\TheBusemannOf[\TheRatVertex]{\DummyArg}=\TheUpperBound$}
        \put(2,0){$\TheBusemannOf[\AltRatVertex]{\DummyArg}=\TheUpperBound$}
        \put(7.5,3.5){$\CloseSimplexOf[\TheUpperBound]{\DummyArg,\TheRatChamber}=\TheRatChamber$}
        \put(7,6){$\CloseSimplexOf[\TheUpperBound]{\DummyArg,\TheRatChamber}=\SetOf{\TheRatVertex}$}
        \put(7,0.7){$\CloseSimplexOf[\TheUpperBound]{\DummyArg,\TheRatChamber}=\SetOf{\AltRatVertex}$}
      \end{picture}
      \hfill\null
      \caption{the upper reduction bound\label{behr:04.fig}}
      {\small
      Shown is the partition of the reduced region as induced by the
      upper reduction bound $\TheUpperBound$. The chamber is
      $\TheRatChamber=\SetOf{\TheRatVertex,\AltRatVertex}$.}
    \end{figure}

    We now fix a reduction datum and say $\TheRatChamber$ \notion{reduces}
    the point $\TheEuPoint$ if it $\TheLowerBound$-reduces the point.
    \begin{observation}\label{behr:independence}
      If the chambers $\TheRatChamber$ and $\AltRatChamber$ both
      reduce the point $\TheEuPoint\in\TheEuBuild$,
      then
      \[
        \CloseSimplexOf[\TheBound]{\TheEuPoint,\TheRatChamber}
        =
        \CloseSimplexOf[\TheBound]{\TheEuPoint,\AltRatChamber}
        \Period
      \]
      for any $\TheBound\geq\TheUpperBound$. Also,
      if $\TheBusemannOf[\TheRatChamber]{\TheEuPoint}\geq\TheUpperBound$,
      then
      \(
        \TheBusemannOf[\TheRatChamber]{\TheEuPoint}
        =
        \TheBusemannOf[\AltRatChamber]{\TheEuPoint}
      \)
      and
      \(
        \CloseSimplexOf[\TheBound]{\TheEuPoint,\TheRatChamber}
        \neq\emptyset
      \).
    \end{observation}
    \begin{proof}
      First,
      \(
        \CloseSimplexOf[\TheBound]{\TheEuPoint,\TheRatChamber}
        \subseteq
        \TheRatSimplex
        :=
        \TheRatChamber\intersect\AltRatChamber
      \).
      Therefore Observation~\ref{behr:restricting.simplex}
      yields
      \(
        \CloseSimplexOf[\TheBound]{\TheEuPoint,\TheRatChamber}
        =
        \CloseSimplexOf[\TheBound]{\TheEuPoint,\TheRatSimplex}
        =
        \CloseSimplexOf[\TheBound]{\TheEuPoint,\AltRatChamber}
      \).

      Furthermore,
      if $\TheBusemannOf[\TheRatChamber]{\TheEuPoint}\geq\TheUpperBound$,
      Corollary~\ref{behr:independent.max}
      implies
      \(
        \TheBusemannOf[\TheRatChamber]{\TheEuPoint}
        =
        \TheBusemannOf[{
          \CloseSimplexOf[\TheBound]{\TheEuPoint,\TheRatChamber}
        }]{
          \TheEuPoint
        }
        =
        \TheBusemannOf[{
          \CloseSimplexOf[\TheBound]{\TheEuPoint,\AltRatChamber}
        }]{
          \TheEuPoint
        }
        =
        \TheBusemannOf[\AltRatChamber]{\TheEuPoint}
      \).
    \end{proof}
    \begin{cor}\label{behr:low}
      For any two chambers $\TheRatChamber,\AltRatChamber
      \in\ChambersOf{\TheRatBuild}$
      both reducing the point $\TheEuPoint\in\TheEuBuild$ and any bound
      $\TheBound\geq\TheUpperBound$, we have
      \[
        \TheBusemannOf[\TheRatChamber]{\TheEuPoint}\leq\TheBound
        \qquad\text{if and only if}
        \qquad
        \TheBusemannOf[\AltRatChamber]{\TheEuPoint}\leq\TheBound
        \Period
        \qed
      \]
    \end{cor}
    We define
    \(
      \CloseSimplexOf[\TheBound]{\TheEuPoint}
      :=
      \Intersection[\TheRatChamber\text{\ reduces\ }\TheEuPoint]{
        \CloseSimplexOf[\TheBound]{\TheEuPoint,\TheRatChamber}
      }
    \).
    We call the (possibly empty) simplex
    \[
      \CloseSimplexOf{\TheEuPoint}
      :=
      \CloseSimplexOf[\TheUpperBound]{\TheEuPoint}
    \]
    \notion{close to} $\TheEuPoint$. It is non-empty
    if $\TheBusemannOf[\TheRatChamber]{\TheEuPoint}\geq\TheUpperBound$ for
    some (and hence any) chamber $\TheRatChamber$ that reduces
    $\TheEuPoint$. In this case, it equals
    $\CloseSimplexOf[\TheUpperBound]{\TheEuPoint,\TheRatChamber}$
    for any reducing chamber.

      A reduction datum is
      \notion{$\TheLattice$-invariant} if
      for each $\TheLatticeElement\in\TheLattice$, each vertex
      $\TheRatVertex\in\VerticesOf{\TheRatBuild}$, and each point
      $\TheEuPoint\in\TheEuBuild$, we have
      \(
        \TheBusemannOf[\TheLatticeElement\TheRatVertex]{
          \TheLatticeElement\TheEuPoint
        }
        =
        \TheBusemannOf[\TheRatVertex]{
          \TheEuPoint
        }
      \).
      A $\TheLattice$-invariant
      reduction datum is called \notion{$\TheLattice$-cocompact} if
      for each $\TheBound\geq\TheUpperBound$, the set
      \begin{eqnarray*}
        \TheRedFiltr[\TheBound]
        &:=&
        \SetOf[\TheEuPoint\in\TheEuBuild]{
          \TheBusemannOf[\TheRatChamber]{\TheEuPoint}
          \leq\TheBound
          \text{\ for all\ }
          \TheRatChamber\in\ChambersOf{\TheRatBuild}
          \text{\ reducing\ }
          \TheEuPoint
        }
        \\
        &=&
        \SetOf[\TheEuPoint\in\TheEuBuild]{
          \TheBusemannOf[\TheRatChamber]{\TheEuPoint}
          \leq\TheBound
          \text{\ for some\ }
          \TheRatChamber\in\ChambersOf{\TheRatBuild}
          \text{\ reducing\ }
          \TheEuPoint
        }
      \end{eqnarray*}
      has compact quotient modulo the action of $\TheLattice$.

      We say that a subset $\TheStar\subseteq\TheEuBuild$ can
      be \notion{uniformly reduced} if there is a chamber
      $\TheRatChamber\in\ChambersOf{\TheRatBuild}$ reducing simultaneously
      all points in $\TheStar$.
      Let $\TheDiam$ be a non-negative real number.
      We say that a reduction datum is \notion{$\TheDiam$-uniform}
      if every subset $\TheStar\subset\TheEuBuild$ of diameter
      at most $\TheDiam$ can be uniformly reduced.

    With these notions, we can rephrase the main theorems of
    Behr-Harder reduction theory as follows:
    \begin{theorem}\label{behr:main}
      For every diameter $\TheDiam$, there is a $\TheDiam$-uniform,
      $\TheLattice$-invariant reduction datum. It is
      $\TheLattice$-cocompact provided $\TheGroupScheme$
      is absolutely almost simple.
    \end{theorem}
    We give a proof in Section~\ref{sec:us}.
    This rendering of the statement is loosely inspired by \cite{Behr:1998}
    where a slightly different version of the sets
    \(
      \TheRedFiltr[\TheBound]
    \)
    was used as a filtration of $\TheEuBuild$: just let $\TheBound$
    tend to $\infty$.
    Ultimately, we will choose $\TheDiam$ large,
    although for our immediate needs, we shall only
    require that closed chambers in $\TheEuBuild$ can be uniformly
    reduced.

  \section{A blueprint for the main argument}
    \label{sec:blueprint}
    The proof of finiteness properties for the group $\TheLattice$ centers
    around a $\TheLattice$-invariant Morse function
    \(
      \TheMorse \mapcolon \TheEuBuild \rightarrow \TheReals
    \)
    with $\TheLattice$-cocompact sublevel sets and highly connected
    descending links. Here, we use the term ``Morse function'' loosely:
    the key feature is that directions can be ascending or descending
    and that there is a gradient, which is the direction of steepest
    ascent.

    In this section, we shall construct an approximation
    \(
      \RoughMorse \mapcolon \TheEuBuild \rightarrow \TheReals
    \)
    that almost suffices: $\RoughMorse$ is $\TheLattice$-invariant,
    $\TheLattice$-cocompact and \emph{generically} has highly connected
    descending links. In the following sections, we will perturb
    $\RoughMorse$ so as to make descending links highly connected everywhere.

    For a point $\TheEuPoint\in\TheEuBuild$, an apartment
    $\TheEuApp$ containing $\TheEuPoint$, and a chamber
    $\TheRatChamber\in\ChambersOf{\TheRatBuild}$ in
    the visual boundary $\BoundaryOf{\TheEuApp}$ and reducing
    $\TheEuPoint$, we define
    \(
      \TheEuPoint[\TheEuApp,\TheRatChamber]
      :=
      \ClosestPointOf[\TheEuApp,\TheRatChamber][\TheUpperBound]{
        \TheEuPoint
      }
    \)
    to be the point in
    \(
      \ConvexChunkOf[\TheEuApp,\TheRatChamber]{\TheUpperBound}
    \)
    closest to $\TheEuPoint$. We let
    \(
      \RoughMorseOf[\TheEuApp,\TheRatChamber]{\TheEuPoint}
    \)
    be the euclidean distance from
    \(
      \TheEuPoint[\TheEuApp,\TheRatChamber]
    \)
    to
    \(
      \TheEuPoint
    \).
    \begin{prop}
      With $\TheEuPoint$, $\TheEuApp$, and $\TheRatChamber$ as above,
      \(
        \TheEuPoint\in\TheRedFiltr[\TheUpperBound]
      \)
      if and only if
      \(
        \TheEuPoint\in
        \ConvexChunkOf[\TheEuApp,\TheRatChamber]{\TheUpperBound}
      \).
      The point $\TheEuPoint[\TheEuApp,\TheRatChamber]$ is also
      the closest point of
      \(
        \ConvexChunkOf[\TheEuApp,\CloseSimplexOf{\TheEuPoint}]{\TheUpperBound}
      \)
      to $\TheEuPoint$.

      If $\TheEuPoint\not\in
      \ConvexChunkOf[\TheEuApp,\TheRatChamber]{\TheUpperBound}$,
      the straight ray
      \(
        \RayFromThrough{\TheEuPoint[\TheEuApp,\TheRatChamber]}{\TheEuPoint}
      \)
      from $\TheEuPoint[\TheEuApp,\TheRatChamber]$ through
      $\TheEuPoint$
      within the apartment
      $\TheEuApp$ meets $\BoundaryOf{\TheEuApp}$ in
      $\CloseSimplexOf{\TheEuPoint}$.
    \end{prop}
    In the latter case, we denote the visual end-point of
      \(
        \RayFromThrough{\TheEuPoint[\TheEuApp,\TheRatChamber]}{\TheEuPoint}
      \)
    by
    \(
      \RoughEndOf[\TheEuApp,\TheRatChamber]{\TheEuPoint}
      \in
      \CloseSimplexOf{\TheEuPoint}
    \).
    \begin{proof}
      One could be tempted to argue
      \(
        \ConvexChunkOf[\TheEuApp,\TheRatChamber]{\TheUpperBound}
        =
        \TheRedFiltr[\TheUpperBound]\intersect\TheEuApp
      \),
      but that is generally not true. However, a point $\TheEuPoint$
      reduced by
      $\TheRatChamber$ belongs to either side or neither side.
      By Corollary~\ref{behr:low},
      \begin{eqnarray*}
        &&
        \TheEuPoint\in
        \ConvexChunkOf[\TheEuApp,\TheRatChamber]{\TheUpperBound}
        \\
        &\Longleftrightarrow&
        \TheBusemannOf[\TheRatChamber]{\TheEuPoint}\leq\TheUpperBound
        \\
        &\Longleftrightarrow&
        \TheBusemannOf[\AltRatChamber]{\TheEuPoint}\leq\TheUpperBound
        \text{\ for all\ }\AltRatChamber\text{\ reducing\ }
        \TheEuPoint
        \\
        &\Longleftrightarrow&
        \TheEuPoint\in\TheRedFiltr[\TheUpperBound]
        \Period
      \end{eqnarray*}

      By Observation~\ref{behr:restricting.simplex},
      \(
        \TheEuPoint[\TheEuApp,\TheRatChamber]
        =
        \ClosestPointOf[{
          \TheEuApp,
          \CloseSimplexOf[\TheUpperBound]{\TheEuPoint,\TheRatChamber}
        }]{
          \TheEuPoint
        },
      \)
      i.e., $\TheEuPoint[\TheEuApp,\TheRatChamber]$ is the point
      of
      \(
        \ConvexChunkOf[{
          \TheEuApp,
          \CloseSimplexOf{\TheEuPoint}
        }]{\TheUpperBound}
      \)
      closest to $\TheEuPoint$. If
      \(
        \TheEuPoint\not\in
        \ConvexChunkOf[\TheEuApp,\TheRatChamber]{\TheUpperBound}
      \)
      then $\TheEuPoint\neq\TheEuPoint[\TheEuApp,\TheRatChamber]$
      and, in the apartment $\TheEuApp$, we can draw the ray
      \(
        [\TheEuPoint[\TheEuApp,\TheRatChamber],\TheEuPoint)
      \).
      This ray lies in the normal cone to
      \(
        \TheEuPoint[\TheEuApp,\TheRatChamber]
      \)
      of the convex body
      \(
        \ConvexChunkOf[{
          \TheEuApp,
          \CloseSimplexOf{\TheEuPoint}
        }]{\TheUpperBound}
      \).
      This normal cone is spanned by the gradients of the
      Busemann functions $\TheBusemann[\TheRatVertex]$
      for $\TheRatVertex\in
      \CloseSimplexOf{\TheEuPoint}$.
      Hence
      \(
        [\TheEuPoint[\TheEuApp,\TheRatChamber],\TheEuPoint)
      \)
      meets $\BoundaryOf{\TheEuApp}$ in
      \(
        \CloseSimplexOf{\TheEuPoint}
      \).
    \end{proof}
    Now, we can see that
    \(
      \RoughMorseOf[\TheEuApp,\TheRatChamber]{\TheEuPoint}
    \)
    and
    \(
      \RoughEndOf[\TheEuApp,\TheRatChamber]{\TheEuPoint}
    \)
    are independent of the choices of $\TheEuApp$ and $\TheRatChamber$.
    \begin{cor}
      If $\AltEuApp$ and $\AltRatChamber$ is another pair with
      $\AltRatChamber\subset\BoundaryOf{\AltEuApp}$,
      $\TheEuPoint\in\AltEuApp$, and $\AltRatChamber$ reducing
      $\TheEuPoint$, then
      \(
        \RoughMorseOf[\TheEuApp,\TheRatChamber]{\TheEuPoint}
        =
        \RoughMorseOf[\AltEuApp,\AltRatChamber]{\TheEuPoint}
      \)
      and
      \(
        \RoughEndOf[\TheEuApp,\TheRatChamber]{\TheEuPoint}
        =
        \RoughEndOf[\AltEuApp,\AltRatChamber]{\TheEuPoint}
      \).
    \end{cor}
    \begin{proof}
      First, we have the following chain of equivalences:
      \begin{eqnarray*}
        &&
        \RoughMorseOf[\TheEuApp,\TheRatChamber]{\TheEuPoint}=\Zero
        \\
        &\Longleftrightarrow&
        \TheEuPoint\in\ConvexChunkOf[\TheEuApp,\TheRatChamber]{\TheUpperBound}
        \\
        &\Longleftrightarrow&
        \TheEuPoint\in\TheRedFiltr[\TheUpperBound]
        \\
        &\Longleftrightarrow&
        \TheEuPoint\in\ConvexChunk[\AltEuApp,\AltRatChamber]
        \\
        &\Longleftrightarrow&
        \RoughMorseOf[\AltEuApp,\AltRatChamber]{\TheEuPoint}=\Zero
      \end{eqnarray*}
      Hence, whether or not
      \(
        \RoughMorseOf[\TheEuApp,\TheRatChamber]{\TheEuPoint}
        =\Zero
      \)
      does depend neither on $\TheEuApp$ nor $\TheRatChamber$.

      Both chambers, $\TheRatChamber$ and $\AltRatChamber$, reduce
      $\TheEuPoint$. Hence,
      \(
        \CloseSimplexOf{\TheEuPoint}\subseteq
        \TheRatChamber\intersect\AltRatChamber
      \)
      and the convex hull
      $\ConvexHullOf{\TheEuPoint,\CloseSimplexOf{\TheEuPoint}}$
      of $\TheEuPoint$ and
      $\CloseSimplexOf{\TheEuPoint}$ lies in
      $\TheEuApp\intersect\AltEuApp$. Hence, there is an isometry
      \(
        \TheIsometry \mapcolon
        \TheEuApp \rightarrow \AltEuApp
      \)
      fixing
      $\ConvexHullOf{\TheEuPoint,\CloseSimplexOf{\TheEuPoint}}$
      pointwise. Hence for each
      $\TheRatVertex\in\CloseSimplexOf{\TheEuPoint}
      \subseteq\TheRatChamber\intersect\AltRatChamber$, the
      following diagram commutes:
      \[
        \xymatrix{
          {\TheEuApp}
          \ar[rr]^{\TheIsometry}
          \ar[dr]_{\TheBusemann[\TheRatVertex]}
          &
          &
          {\AltEuApp}
          \ar[dl]^{\TheBusemann[\TheRatVertex]}
          \\
          &
          {\TheReals}
          &
        }
      \]
      Hence, $\TheIsometry$ identifies
      $\ConvexChunk[\TheEuApp,\CloseSimplexOf{\TheEuPoint}]$
      with
      $\ConvexChunk[\AltEuApp,\CloseSimplexOf{\TheEuPoint}]$ and
      $\TheEuPoint[\TheEuApp,\TheRatChamber]$ with
      $\TheEuPoint[\AltEuApp,\AltRatChamber]$. It follows that
      \(
        \RoughEndOf[\TheEuApp,\TheRatChamber]{\TheEuPoint}
        =
        \RoughEndOf[\AltEuApp,\AltRatChamber]{\TheEuPoint}
      \)
      and
      \(
        \RoughMorseOf[\TheEuApp,\TheRatChamber]{\TheEuPoint}
        =
        \RoughMorseOf[\AltEuApp,\AltRatChamber]{\TheEuPoint}
      \).
    \end{proof}
    Thus, we define
    \(
      \RoughMorseOf{\TheEuPoint}
      :=
      \RoughMorseOf[\TheEuApp,\TheRatChamber]{\TheEuPoint}
    \)
    and
    \(
      \RoughEndOf{\TheEuPoint}
      :=
      \RoughEndOf[\TheEuApp,\TheRatChamber]{\TheEuPoint}
    \)
    where $\TheRatChamber$ is a chamber reducing $\TheEuPoint$ and
    $\TheEuApp$ is an apartment containing $\TheEuPoint$ and whose
    boundary contains $\TheRatChamber$.

    The properties of $\RoughMorse$ that we are about to establish roughly
    say that $\RoughMorse$ qualifies as a Morse function for the analysis
    of $\TheLattice$.

    \begin{observation}
      The function $\RoughMorse$ is defined entirely in terms of a
      $\TheLattice$-invariant reduction datum, hence it is  itself
      $\TheLattice$-invariant. I.e.,
      \(
        \RoughMorseOf{\TheLatticeElement\TheEuPoint}
        =
        \RoughMorseOf{\TheEuPoint}
      \)
      for any $\TheLatticeElement\in\TheLattice$ and any point
      $\TheEuPoint\in\TheEuBuild$.\qed
    \end{observation}

    \begin{prop}\label{blueprint:cocompact}
      Each sublevel set $\SubLevelSet{\RoughMorse}{\TheLevel}$ has
      compact quotient modulo the action of $\TheLattice$.
    \end{prop}
    \begin{proof}
      We use the $\TheLattice$-cocompactness of the reduction datum.
      More precisely, we show that for any $\TheLevel\geq \Zero$ the
      sublevel set $\SubLevelSet{\RoughMorse}{\TheLevel}$ is
      contained in the $\TheLattice$-cocompact set
      $\TheRedFiltr[\TheLevel+\TheUpperBound]$. It is easier to
      prove the contrapositive. So assume that we have a point
      $\TheEuPoint\not\in\TheRedFiltr[\TheLevel+\TheUpperBound]$.
      We have to argue that $\RoughMorseOf{\TheEuPoint}>\TheLevel$.

      As
      \(
        \TheEuPoint\not\in
        \TheRedFiltr[\TheLevel+\TheUpperBound]
      \)
      there is a chamber $\TheRatChamber$ reducing $\TheEuPoint$ with
      \(
        \TheBusemannOf[\TheRatChamber]{\TheEuPoint}>
        \TheLevel+\TheUpperBound
      \),
      i.e., there is
      a vertex $\TheRatVertex\in\TheRatChamber$ with
      $\TheBusemannOf[\TheRatVertex]{\TheEuPoint}>\TheLevel+\TheUpperBound$.
      To estimate $\RoughMorseOf{\TheEuPoint}$ we additionally choose an
      apartment $\TheEuApp$ containing $\TheEuPoint$ that has $\TheRatChamber$
      in its visual boundary. Then
      \begin{eqnarray*}
        \RoughMorseOf{\TheEuPoint}
        &=&
        \RoughMorseOf[\TheEuApp,\TheRatChamber]{\TheEuPoint}
        \\
        &=&
        \DistanceOf{\TheEuPoint[\TheEuApp,\TheRatChamber]}{\TheEuPoint}
        \\
        &\geq&
        \TheBusemannOf[\TheRatVertex]{\TheEuPoint}
        -
        \TheBusemannOf[\TheRatVertex]{\TheEuPoint[\TheEuApp,\TheRatChamber]}
        \\
        &>&
        (\TheLevel+\TheUpperBound)
        -
        \TheUpperBound
        =
        \TheLevel
        \Period
      \end{eqnarray*}
      This completes the proof. The first inequality follows in general from
      the definition of Busemann functions and the triangle inequality. However,
      as $\TheEuApp$ is a euclidean space, we have a more precise statement:
      \(
        \TheBusemannOf[\TheRatVertex]{\TheEuPoint}
        -
        \TheBusemannOf[\TheRatVertex]{\TheEuPoint[\TheEuApp,\TheRatChamber]}
        =
        \cos(\TheAngle)
        \DistanceOf{\TheEuPoint[\TheEuApp,\TheRatChamber]}{\TheEuPoint}
      \)
      where $\TheAngle$ is the angle between the ray
      $[\TheEuPoint[\TheEuApp,\TheRatChamber],\TheEuPoint)$ and the
      gradient of $\TheBusemann[\TheRatVertex]$.
    \end{proof}

    \begin{prop}\label{blueprint:continuous}
      The function $\RoughMorse\mapcolon\TheEuBuild\rightarrow\TheReals$
      is continuous.
    \end{prop}
    \begin{proof}
      Here, we use the hypothesis that the chosen reduction datum is
      sufficiently uniform. More precisely, we shall use that every chamber
      in $\TheEuBuild$ can be uniformly reduced. So let $\TheEuChamber$
      be a chamber in $\TheEuBuild$,
      let $\TheRatChamber\in\ChambersOf{\TheRatBuild}$ be a chamber
      uniformly reducing $\TheEuChamber$, and let $\TheEuApp$ be an
      apartment containing $\TheEuChamber$ whose visual boundary
      contains $\TheRatChamber$. Then, the functions
      \(
        \RoughMorse
      \)
      and
      \(
        \RoughMorse[\TheEuApp,\TheRatChamber]
      \)
      agree on $\TheEuChamber$. In particular, $\RoughMorse$ restricts to
      a continuous function on $\TheEuChamber$. As the euclidean building
      $\TheEuBuild$ carries the weak topology with respect to its chambers,
      the function $\RoughMorse$ is continuous.
    \end{proof}

    Morse functions are supposed to be differentiable and not merely
    continuous. The following statements make precise the sense in
    which $\RoughMorse$ induces a gradient field and what should be
    considered flow lines of this gradient field. First we treat the
    gradient, i.e., we describe the behavior of $\RoughMorse$ on small
    scales.
    \begin{prop}\label{blueprint:gradient}
      Let $\TheEuPoint$ and $\AltEuPoint$ be two points in $\TheEuBuild$
      that lie within a common closed chamber. Then
      \[
        \RoughMorseOf{\AltEuPoint}
        -
        \RoughMorseOf{\TheEuPoint}
        \leq
        \DistanceOf{\AltEuPoint}{\TheEuPoint}
      \]
      with equality if and only if $\AltEuPoint$ lies on
      the ray $\RayFromTo{\TheEuPoint}{\RoughEndOf{\TheEuPoint}}$.
      In case of equality, moreover
      $\RoughEndOf{\TheEuPoint}=\RoughEndOf{\AltEuPoint}$.
    \end{prop}
    Thus, we define the gradient
    of $\RoughMorse$ at the point $\TheEuPoint$ to be the direction
    $\GradOf[\TheEuPoint]{\RoughMorse}\in\LinkOf{\TheEuPoint}$
    defined by the ray $\RayFromTo{\TheEuPoint}{\RoughEndOf{\TheEuPoint}}$.
    It is the unique direction of fastest ascent, and $\RoughMorse$
    grows in that direction with unit speed.
    \begin{proof}
      Choose $\TheRatChamber$ and $\TheEuApp$ so that
      $\RoughMorse$ and $\RoughMorse[\TheEuApp,\TheRatChamber]$
      agree on the segment $\LineSegment{\TheEuPoint}{\AltEuPoint}$.
      Then
      \begin{eqnarray*}
        \RoughMorseOf{\AltEuPoint}
        -
        \RoughMorseOf{\TheEuPoint}
        &=&
        \RoughMorseOf[\TheEuApp,\TheRatChamber]{\AltEuPoint}
        -
        \RoughMorseOf[\TheEuApp,\TheRatChamber]{\TheEuPoint}
        \\
        &=&
        \DistanceOf{\AltEuPoint[\TheEuApp,\TheRatChamber]}{\AltEuPoint}
        -
        \DistanceOf{\TheEuPoint[\TheEuApp,\TheRatChamber]}{\TheEuPoint}
        \\
        &\leq&
        \DistanceOf{\TheEuPoint[\TheEuApp,\TheRatChamber]}{\AltEuPoint}
        -
        \DistanceOf{\TheEuPoint[\TheEuApp,\TheRatChamber]}{\TheEuPoint}
        \\
        &\leq&
        \DistanceOf{\AltEuPoint}{\TheEuPoint}
        \Period
      \end{eqnarray*}
      We have equality in the last step  if and only if $\TheEuPoint$ lies
      on the straight line segment
      \(
        \LineSegment{\TheEuPoint[\TheEuApp,\TheRatChamber]}{\AltEuPoint}.
      \)
      In this case, however, the segment
      \(
        \LineSegment{\TheEuPoint[\TheEuApp,\TheRatChamber]}{\AltEuPoint}
      \)
      is normal to $\ConvexChunkOf[\TheEuApp,\TheRatChamber]{\TheUpperBound}$. Therefore,
      \(
        \AltEuPoint[\TheEuApp,\TheRatChamber]
        =
        \TheEuPoint[\TheEuApp,\TheRatChamber]
      \)
      whence we have equality throughout and
      $\RoughEndOf{\TheEuPoint}=\RoughEndOf{\AltEuPoint}$.
    \end{proof}
    \begin{observation}\label{blueprint:flow.line}
      Every line segment in $\TheEuBuild$ cuts through finitely many
      chambers. Hence, one can easily ``integrate'' the local information
      provided by the previous proposition: for any two points
      $\AltEuPoint$ and $\TheEuPoint$, the geodesic segment
      $\LineSegment{\TheEuPoint}{\AltEuPoint}$ can be subdivided into
      finitely many segments each of which is supported by a closed
      chamber. That is, there are points
      \[
        \TheEuPoint=\TheEuPoint[\Zero],
        \TheEuPoint[\One],\ldots,
        \TheEuPoint[\TheLastIndex]=\AltEuPoint
        \in
        \LineSegment{\TheEuPoint}{\AltEuPoint}
      \]
      such that
      $\LineSegment{\TheEuPoint}{\AltEuPoint}$
      is the concatenation of the segments
      $\LineSegment{\TheEuPoint[\TheIndex]}{\TheEuPoint[\TheIndex+\One]}$.
      Applying the previous proposition to each of those segments, we obtain
      \[
        \RoughMorseOf{\AltEuPoint}
        \leq
        \RoughMorseOf{\TheEuPoint}
        +
        \DistanceOf{\AltEuPoint}{\TheEuPoint}
      \]
      with equality if and only if $\AltEuPoint\in
      \RayFromTo{\TheEuPoint}{\RoughEndOf{\TheEuPoint}}$. In the case
      of equality, $\RoughEndOf{\AltEuPoint}=\RoughEndOf{\TheEuPoint}$.\qed
    \end{observation}
    Thus, we can regard the unit speed geodesic ray
    from $\TheEuPoint$ to $\RoughEndOf{\TheEuPoint}$ as
    a flow line of the gradient field.

    Now we are ready to give a first, albeit failing, attempt to prove
    the Rank~Theorem. The actual proof given below follows the same
    lines, and this argument will help identify exactly the shortcomings
    of $\RoughMorse$ that we have to address in the following sections.

    The $\ThePlaceSet$-arithmetic group $\TheLattice$
    acts on the euclidean building $\TheEuBuild$, which is a
    \CAT{\Zero}-space and therefore
    contractible. Cell stabilizers of this action are finite
    as $\TheEuBuild$ is a proper \CAT{\Zero}-space and $\TheLattice$ is
    discrete in its isometry group. For each positive
    real $\TheLevel$ let
    \(
      \TheEuBuildOf{\TheLevel}
    \)
    be the largest subcomplex of $\TheEuBuild$ fully contained in
    the sublevel set $\SubLevelSet{\RoughMorse}{\TheLevel}$.
    By Proposition~\ref{blueprint:cocompact}, the orbit space of
    $\TheEuBuildOf{\TheLevel}$ modulo the action of $\TheLattice$
    is compact. Should $\TheEuBuildOf{\TheLevel}$ be
    $(\TheDim-\Two)$-connected for some $\TheLevel$, then
    \cite[Propositions~1.1 and~3.1]{Brown:1987} would imply
    that $\TheLattice$ is of type \FType[\TheDim-\One].

    The aim of combinatorial Morse theory is to describe how the
    homotopy type of sublevel complexes
    \(
      \TheEuBuildOf{\TheLevel}
    \)
    change as $\TheLevel$ varies. This description should be in terms
    of purely local information about
    the function $\RoughMorse$. Crucial are \notion{descending links},
    i.e., the set of directions at a given point $\TheEuPoint$ along
    which the function $\RoughMorse$ decreases. More precisely, a
    cell $\HorEuCell$ containing the vertex $\TheEuVertex$ is considered
    descending if $\RoughMorse$ assumes its maximum on $\HorEuCell$
    at and only at $\TheEuVertex$. The descending cells at $\TheEuVertex$
    form the descending link at $\TheEuVertex$. If all vertices of
    $\TheEuVertex$ with
    \(
      \TheLevel\leq
      \RoughMorseOf{\TheEuVertex}
      \leq
      \AltLevel
    \)
    have $\LocalConnectivity$-connected descending links, the inclusion
    \(
      \TheEuBuildOf{\TheLevel}
      \subseteq
      \TheEuBuildOf{\AltLevel}
    \)
    of sublevel complexes induces isomorphisms in homotopy groups
    $\HomotopyGroup[\TheIndex]$ for all $\TheIndex\leq\LocalConnectivity$.
    As the euclidean building
    $\TheEuBuild$ is contractible, these isomorphisms in the
    $\HomotopyGroup[\TheIndex]$ imply that
    already some $\TheEuBuildOf{\TheLevel}$ is
    $\LocalConnectivity$-connected. Thus, using Brown's criterion
    from above, we are reduced to the question of whether descending
    links are $(\TheDim-\Two)$-connected.

    In smooth Morse theory, the descending link is an infinitesimal
    notion and a direction is descending if it spans an obtuse angle with
    the gradient. In the combinatorial setting, whether an edge determines
    an ascending or descending direction in the link of an adjacent
    vertex depends on the values of the Morse function at the end points.
    Hence, the descending link is a local rather than an infinitesimal
    notion. Often, however, the infinitesimal behavior of the Morse function
    is good enough an approximation: at a \emph{generic}
    vertex $\TheEuVertex$ an adjacent edge is descending if and only if
    it spans an
    obtuse angle with the gradient $\GradOf[\TheEuVertex]{\RoughMorse}$.
    At those vertices where predictions based solely on
    the gradient are correct, descending links
    are therefore hemisphere complexes, whose connectivity
    properties are given in \cite{Schulz:2010}.

    \begin{figure}[ht]
      \null\hfill
      \includegraphics{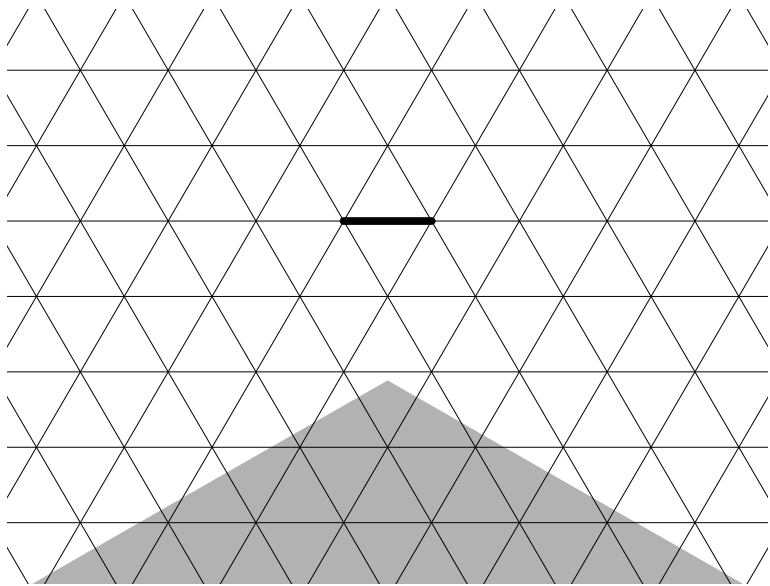}%
      \hfill\null
      \caption{a bad edges\label{blueprint:euclidean.fig}}
      {\small
      The shaded area in the bottom
      is $\ConvexChunkOf[\TheEuApp,\TheRatChamber]{\TheUpperBound}$. The ``corridor''
      above the tip yields an infinite family of ``bad edges''.}
    \end{figure}
    This strategy almost succeeds. Unfortunately, the
    \notion{gradient criterion} is sometimes wrong.
    Figure~\ref{blueprint:euclidean.fig} shows an apartment
    $\TheEuApp$ with
    the convex set
    $\ConvexChunkOf[\TheEuApp,\TheRatChamber]{\TheUpperBound}$ drawn
    in. Suppose, $\TheRatChamber$ reduces the marked edge. Then,
    the edge spans an obtuse angle with the gradient at either
    end point. Hence, the gradient criteria for both vertices are
    in direct conflict with one another. The reason is that the
    gradient criterion only makes correct predictions on an infinitesimal
    scale. The Morse function $\RoughMorse$ actually decreases along the edge
    from either end toward the center. Beyond the center point, however,
    $\RoughMorse$ increases again, spoiling the prediction based on the
    gradient. The picture~\ref{blueprint:euclidean.fig} also shows
    that this problem occurs ``arbitrarily far out'', i.e., we cannot
    avoid it by considering $\TheEuBuildOf{\TheLevel}$ for some
    high value of $\TheLevel$: since the reduction datum is geometric,
    the chamber $\TheRatChamber$ also reduces all the edges in the
    ``corridor'' above the marked edge whence $\RoughMorse$ along
    these edges can be read off in the picture.

    \begin{figure}[ht]
      \null\hfill
      \includegraphics{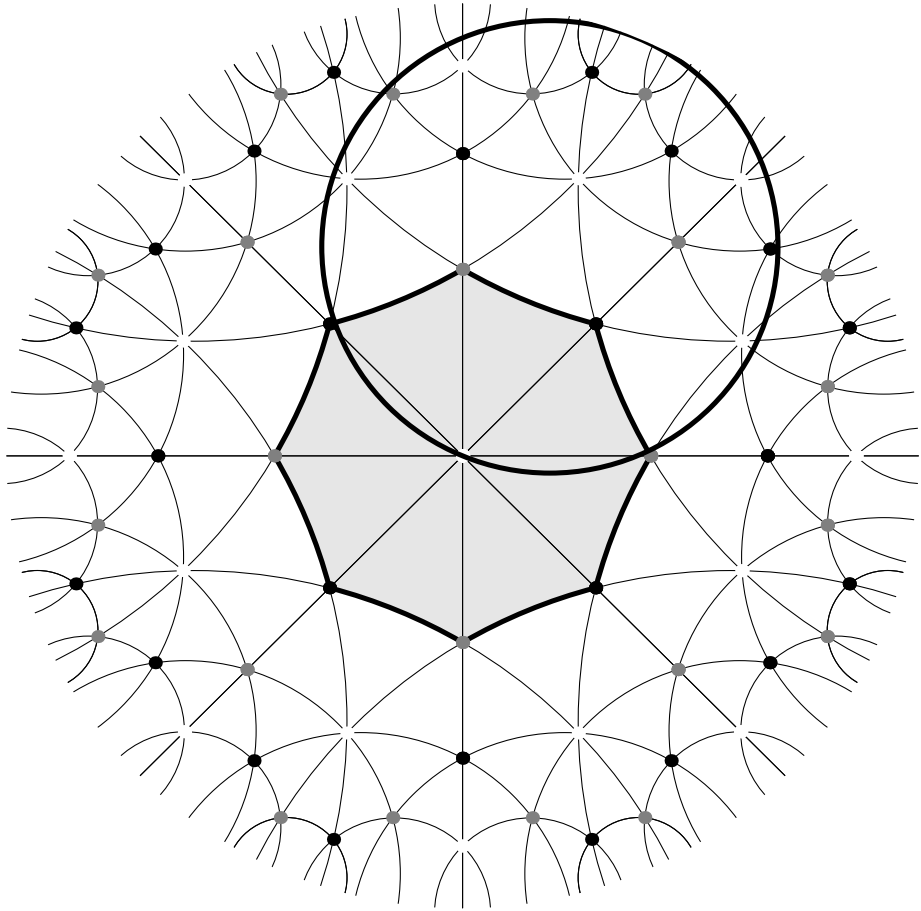}%
      \hfill\null
      \caption{a star in a hyperbolic Coxeter complex%
      \label{blueprint:hyperbolic.fig}}
    \end{figure}
    We note that the predictive powers of the gradient hinge upon
    the underlying apartments being \emph{euclidean}
    spaces. A Coxeter complex in hyperbolic space does not scale:
    its edges have lengths determined by the Coxeter diagram. If
    vertices are too far apart, the infinitesimal nature of the
    gradient renders it useless even for predicting the value of
    a Busemann function on neighbors. This phenomenon is shown
    in Figure~\ref{blueprint:hyperbolic.fig}.
    In particular, one cannot
    expect descending links to be hemisphere complexes.
    This matches examples of Abramenko in the compact hyperbolic
    case where the finiteness length of a lattice falls short of
    the dimension of the building on which the groups acts naturally.

    Our main task will be to alter the Morse function $\RoughMorse$
    to make gradients consistent, i.e., we do not want to see edges
    that span obtuse angles with gradient vectors at either end.
    The obstruction is, of course, that we want to keep high
    connectivity of descending links at vertices where they are
    already fine.

    The remainder of this paper is organized as follows. After
    some preliminaries on zonotopes in Section~\ref{sec:zonotopes}
    and spherical buildings in Section~\ref{sec:bernd},
    we define, in Section~\ref{sec:height},
    a primary Morse function $\TheMorse$ (the \notion{height}),
    which is a perturbation of $\RoughMorse$ discussed above.
    It is $\TheLattice$-invariant, $\TheLattice$-cocompact,
    continuous, and induces a gradient field with geodesic flow
    lines. It improves upon $\RoughMorse$ in that the gradient
    criterion never leads to inconsistencies. However, we cannot
    avoid \notion{$\TheMorse$-flat} cells, e.g., edges on which
    $\TheMorse$ is constant. To break ties, i.e., to determine which
    vertex of such an edge to add first in filtering the euclidean
    building, we introduce a secondary and even a tertiary Morse
    function in Section~\ref{sec:depth}. Here we rely on the notion
    of \notion{depth} introduced in \cite{Bux.Wortman:2008}
    and further developed in \cite{Witzel:2010}.
    The analysis of descending
    links for $\TheMorse$ is carried out in Section~\ref{sec:links}
    whereas Section~\ref{sec:proof} derives the Rank~Theorem.
    The final three sections are devoted to reduction theory.

  \section{A small convex geometry toolkit}\label{sec:zonotopes}
    \begin{figure}[ht]
      \null\hfill
      \includegraphics{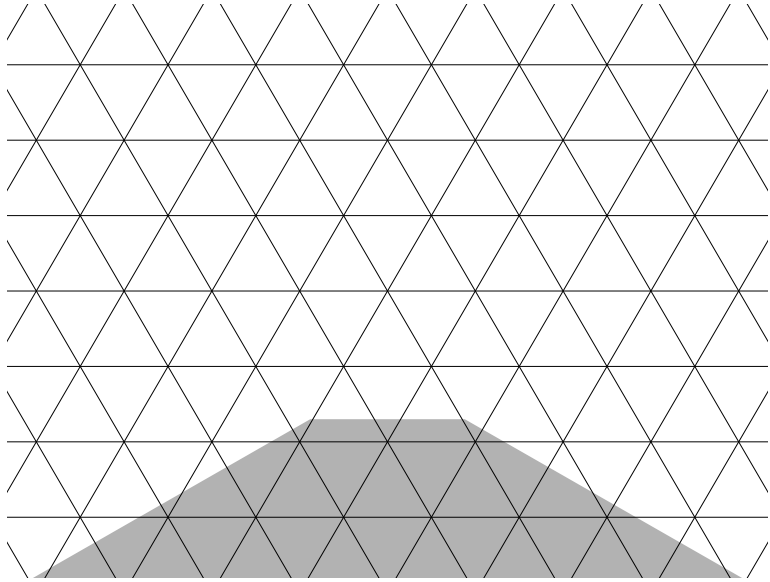}
      \hfill\null
      \caption{flattening the tip\label{zonotopes:euclidean.fig}}
      {\small
      Flattening the tip removes corridors of bad edges.}
    \end{figure}
    We address the problem illustrated in
    Figure~\ref{blueprint:euclidean.fig} by changing the shape
    of $\ConvexChunkOf[\TheEuApp,\TheRatChamber]{\TheUpperBound}$.
    We will flatten the
    tip as shown in Figure~\ref{zonotopes:euclidean.fig}. Then the
    gradient field becomes consistent along edges.

    Let $\TheEuSpace$ denote a euclidean space with inner
    product $\InnerProd{\DummyArg}{\DummyArg}$ and origin
    $\SpZero$. Let $\TheFace$ be the face of some convex
    polytope $\ThePolytope$. The \notion{normal cone}
    \[
      \NormalConeOf{\TheFace}
      :=
      \SetOf[
        \TheNormalVect\in\TheEuSpace
      ]{
        \InnerProd{\TheNormalVect}{\ThePolytopePoint}
        =
        \Max[\AltPolytopePoint\in\ThePolytope]{
          \InnerProd{\TheNormalVect}{\AltPolytopePoint}
        }
        \text{\ for all\ }
        \ThePolytopePoint\in\TheFace
      }
    \]
    is the set of all $\TheNormalVect\in\TheEuSpace$ such that
    the function $\InnerProd{\TheNormalVect}{\DummyArg}$ restricted
    to $\ThePolytope$ assumes its maximum on the points in $\TheFace$.
    It is a closed convex cone. For any point
    $\TheSpVector\in\TheEuSpace$, the closest point projection
    onto $\ThePolytope$ satisfies:
    \[
      \ClosestPointOf[\ThePolytope]{\TheSpVector}\in\TheFace
      \qquad\text{if and only if}\qquad
      \TheSpVector-\ClosestPointOf[\ThePolytope]{\TheSpVector}
      \in\NormalConeOf{\TheFace}
    \]
    For a finite subset $\TheZonoDef\subset\TheEuSpace$, the
    compact convex polytope
    \[
      \TheZonotopeOf{\TheZonoDef}
      :=
      \SetOf[{
        \Sum[\TheZonoVect\in\TheZonoDef]{
          \TheCoefficient[\TheZonoVect]
          \TheZonoVect
        }
      }]{
        \Zero\leq\TheCoefficient[\TheZonoVect]\leq\One
        \text{\ for all\ }
        \TheZonoVect\in\TheZonoDef
      }
    \]
    is called the \notion{zonotope} spanned by $\TheZonoDef$.
    This construction ensures:
    \begin{observation}\label{zonotopes:parallel.lines}
      Through every point $\TheZonoPoint\in\TheZonotopeOf{\TheZonoDef}$
      and every $\TheZonoVect\in\TheZonoDef$, there is a line segment
      parallel to $[\SpZero,\TheZonoVect]$ inside
      $\TheZonotopeOf{\TheZonoDef}$.\qed
    \end{observation}

    Let $\ThePolytope\subset\TheEuSpace$ be a compact convex polytope.
    We call $\TheZonoDef\subset\TheEuSpace$ \notion{saturated}
    with respect to $\ThePolytope$, if for any two vertices
    $\ThePolyVertex,\AltPolyVertex\in\ThePolytope$,
    at least one of
    \(
      \ThePolyVertex-\AltPolyVertex
    \)
    or
    \(
      \AltPolyVertex-\ThePolyVertex
    \)
    lies in $\TheZonoDef$.
    \begin{prop}\label{zonotopes:vertex.in}
      Let $\TheZonoDef\subset\TheEuSpace$ be saturated with respect
      to the compact convex polytope $\ThePolytope\subset\TheEuSpace$.
      Then any translate
      of $\TheZonotopeOf{\TheZonoDef}$ that intersects $\ThePolytope$
      contains a vertex of $\ThePolytope$. Equivalently: whenever
      one translates $\ThePolytope$ to intersect
      $\TheZonotopeOf{\TheZonoDef}$ one has to move a vertex
      of $\ThePolytope$ into $\TheZonotopeOf{\TheZonoDef}$.
    \end{prop}
    \begin{proof}
      Since differences of vertices
      are invariant under translation, we may also assume that
      $\ThePolytope$ intersects $\TheZonotopeOf{\TheZonoDef}$.
      Let $\ThePolytopePoint$ be a point in the intersection.
      We have to show that $\TheZonotopeOf{\TheZonoDef}$ contains
      a vertex of $\ThePolytope$.

      Let $\ThePolyVertex$ and $\AltPolyVertex$ be two vertices
      of $\ThePolytope$ on opposite parallel supporting hyperplanes.
      Choose the labels so that
      $\ThePolyVertex-\AltPolyVertex\in\TheZonoDef$.
      Any translate of the segment
      \(
        \LineSegment{\ThePolyVertex}{\AltPolyVertex}
      \)
      through the point $\ThePolytopePoint$ meets the
      boundary of $\ThePolytope$, and
      Observation~\ref{zonotopes:parallel.lines} implies that there is
      such a translate that stays inside $\TheZonotopeOf{\TheZonoDef}$.
      Hence, $\ThePolytope$ has a face $\TheFace$
      that intersects $\TheZonotopeOf{\TheZonoDef}$.
      By induction on the dimension, we may therefore conclude
      that $\TheZonotopeOf{\TheZonoDef}$ contains a vertex of $\TheFace$.
    \end{proof}
    \begin{cor}\label{zonotopes:cover}
      If $\TheZonoDef$ and $\ThePolytope$ are as in the proposition,
      i.e., $\TheZonoDef$ is saturated with respect to $\ThePolytope$,
      then
      \(
        \ThePolytope+\TheZonotopeOf{\TheZonoDef}
        =
        \ThePolytope[][(\Zero)]+\TheZonotopeOf{\TheZonoDef}
      \)
      where $\ThePolytope[][(\Zero)]$ is the set of vertices
      of $\ThePolytope$.
    \end{cor}
    \begin{proof}
      Let $\TheEuPoint\in\ThePolytope+\TheZonotopeOf{\TheZonoDef}$,
      i.e., there is a point $\ThePolytopePoint\in\ThePolytope$
      with
      $\TheEuPoint\in\ThePolytopePoint+\TheZonotopeOf{\TheZonoDef}$.
      Then
      \(
        \ThePolytopePoint\in\TheEuPoint+\TheZonotopeOf{-\TheZonoDef}
      \).
      By Proposition~\ref{zonotopes:vertex.in}, there is a vertex
      $\ThePolyVertex\in\ThePolytope[][(\Zero)]$ with
      \(
        \ThePolyVertex\in\TheEuPoint+\TheZonotopeOf{-\TheZonoDef}
      \),
      which implies
      \(
        \TheEuPoint\in\ThePolyVertex+\TheZonotopeOf{\TheZonoDef}
      \).
    \end{proof}
    Let $\TheConvexMap\mapcolon\TheEuSpace\rightarrow\TheReals$
    be a continuous function and let $\TheConvexSet$ be a non-empty
    compact convex subset of $\TheEuSpace$. Define:
    \begin{eqnarray*}
      \TheConvexMap[\TheConvexSet] \mapcolon \TheEuSpace
      &\longrightarrow&\TheReals
      \\
      \TheEuPoint & \mapsto &
      \min\SetOf[\TheConvexMapOf{\AltEuPoint}]{
        \AltEuPoint\in\TheEuPoint+\TheConvexSet
      }
    \end{eqnarray*}
    \begin{observation}\label{zonotopes:convex.function}
      If $\TheConvexMap$ is convex
      then so is $\TheConvexMap[\TheConvexSet]$.\qed
    \end{observation}
    \begin{prop}\label{zonotopes:minset.maxset}
      Let $\TheConvexMap\mapcolon\TheEuSpace\rightarrow\TheReals$ be
      continuous,
      let $\ThePolytope\subset\TheEuSpace$ be a compact convex
      polytope, and let $\TheZonoDef\subset\TheEuSpace$ be a finite
      subset containing all differences $\TheSpVector-\AltSpVector$
      for any two vertices $\TheSpVector$ and $\AltSpVector$ of
      $\ThePolytope$. Then the following hold:
      \begin{enumerate}
        \item\label{minset}
          The min-set of $\TheConvexMap[\TheZonotopeOf{\TheZonoDef}]$
          on $\ThePolytope$ contains a vertex of $\ThePolytope$.
        \item\label{maxset}
          If $\TheConvexMap$ is a convex function then the
          max-set of $\TheConvexMap[\TheZonotopeOf{\TheZonoDef}]$
          on $\ThePolytope$
          is a face of $\ThePolytope$.
      \end{enumerate}
    \end{prop}
    \begin{proof}
      Part~(\ref{minset}) follows immediately from
      Corollary~\ref{zonotopes:cover}.
      As for part~(\ref{maxset}), put:
      \[
        \TheSpVertexSet :=
        \SetOf[
          \TheSpVertex
          \in
          \MaxSetOf{\TheConvexMap[\TheZonotopeOf{\TheZonoDef}]|_{\ThePolytope}}
        ]{
          \TheSpVertex \text{\ is vertex of\ }
          \ThePolytope
        }
      \]
      Since $\TheConvexMap[\TheZonotopeOf{\TheZonoDef}]$
      is a convex function, $\TheSpVertexSet$
      is not empty. Let $\TheFace$ be the smallest face of
      $\ThePolytope$ containing $\TheSpVertexSet$, and let
      $\ConvexHullOf{\TheSpVertexSet}$ be the convex hull of
      $\TheSpVertexSet$. By part~(a),
      the function
      \(
        \TheConvexMap[\TheZonotopeOf{\TheZonoDef}]
        |_{\ConvexHullOf{\TheSpVertexSet}}
      \)
      assumes its minimum in a vertex, i.e., a point in $\TheSpVertexSet$.
      Hence $\TheConvexMap$ is constant on $\ConvexHullOf{\TheSpVertexSet}$,
      as $\TheSpVertexSet$ consists of points in $\ThePolytope$ where
      $\TheConvexMap[\TheZonotopeOf{\TheZonoDef}]$ is maximal.
      As $\TheFace$ is the smallest face of $\ThePolytope$ containing
      $\ConvexHullOf{\TheSpVertexSet}$ and
      $\TheConvexMap[\TheZonotopeOf{\TheZonoDef}]$ is a convex function,
      $\TheConvexMap[\TheZonotopeOf{\TheZonoDef}]$ is constant on
      $\TheFace$.
    \end{proof}

    \begin{example}\label{zonotopes:distance.to.body}
      Let $\ConvexChunk\subset\TheEuSpace$
      be a non-empty closed convex set. Consider the
      distance from $\ConvexChunk$ as a function
      \begin{eqnarray*}
        \TheConvexMap \mapcolon \TheEuSpace &
        \longrightarrow &
        \TheReals \\
        \TheSpPoint & \mapsto &
          \DistanceOf{\TheSpPoint}{\ConvexChunk}
          \Period
      \end{eqnarray*}
      Then, for any non-empty compact convex set
      \(
        \TheConvexSet\subset\TheEuSpace
      \)
      and for any $\TheSpPoint\in\TheEuSpace$, we have:
      \begin{eqnarray*}
        \TheConvexMapOf[\TheConvexSet]{\TheSpPoint}
        &=&
        \MinOf[\AltSpPoint\in\TheConvexSet]{
          \TheConvexMapOf{\TheSpPoint+\AltSpPoint}
        }
        \\
        &=&
        \MinOf[\AltSpPoint\in\TheConvexSet]{
          \DistanceOf{\TheSpPoint+\AltSpPoint}{\ConvexChunk}
        }
        \\
        &=&
        \MinOf[\AltSpPoint\in\TheConvexSet,\ThrSpPoint\in\ConvexChunk]{
          \DistanceOf{\TheSpPoint+\AltSpPoint}{\ThrSpPoint}
        }
        \\
        &=&
        \MinOf[\AltSpPoint\in\TheConvexSet,\ThrSpPoint\in\ConvexChunk]{
          \DistanceOf{\TheSpPoint}{\ThrSpPoint-\AltSpPoint}
        }
        \\
        &=&
        \DistanceOf{\TheSpPoint}{\ConvexChunk-\TheConvexSet}
      \end{eqnarray*}
      Hence, the functions
      \(
        \TheConvexMap[\TheConvexSet]
      \)
      and
      \(
        \DistanceOf{\DummyArg}{\ConvexChunk-\TheConvexSet}
      \)
      coincide.\qed
    \end{example}

  \section{Some subcomplexes of spherical buildings}
    \label{sec:bernd}
    To deduce finiteness properties, we use the well-established
    technique of filtering a complex upon which the group acts.
    The main task, as usual, is to control the homotopy type of
    relative links that arise in the filtration.
    In this section, we collect the results concerning connectivity
    properties of those subcomplexes of spherical buildings that we
    will encounter.

    Let $\TheSpace$ be euclidean or hyperbolic space or a round
    sphere. We call an intersection of a non-empty family of closed
    half-spaces (or hemispheres in the latter case)
    \notion{demi-convex}. We call a subset
    of $\TheSpace$ \notion{fat} if it has non-empty interior.
    Note that a proper open convex subset
    of $\TheSpace$ is contained in an open hemisphere.
    \begin{observation}\label{bernd:deformation.retraction}
      Let $\NewPart\subset\TheSpace$ be fat and demi-convex
      and let $\BadSet\subset\TheSpace$ be proper, open, and convex.
      If $\NewPart$ and $\BadSet$ intersect, then
      $\NewPart\setminus\BadSet$ strongly deformation retracts
      onto the boundary part
      $\BoundaryOf{\NewPart}\setminus\BadSet$.
    \end{observation}
    \begin{proof}
      Note that $\BadSet$ intersects the interior of $\NewPart$
      since every boundary point of the convex set $\NewPart$ is
      an accumulation point of interior points because $\NewPart$
      is fat. Choose $\TheRadialCenter$ in the intersection.
      Note that $\NewPart$ is star-like with
      regard to $\TheRadialCenter$, and the geodesic projection away from
      $\TheRadialCenter$ restricts to the deformation retraction we need.
    \end{proof}
    We call a CW-complex \notion{geometric} if its cells carry a spherical,
    euclidean, or hyperbolic structure in which they are demi-convex (i.e.,
    each cell is an intersection of half-spaces in the model geometry).
    Also, we require attaching maps to be isometric embeddings.
    Iterated application of the projection trick yields:
    \begin{prop}\label{bernd:push_in}
      Suppose that $\TheGeometricComplex$ is a geometric CW-complex.
      Let $\BadSet$ be an open subset of $\TheGeometricComplex$
      that intersects each cell in a convex set.
      Then there is a strong deformation retraction
      \[
        \TheDeformation[\TheGeometricComplex]
        \mapcolon
        \TheGeometricComplex\setminus\BadSet
        \longrightarrow
        \TheGeometricComplex[][\BadSet]
      \]
      of $\TheGeometricComplex\setminus\BadSet$ onto its maximal subcomplex.
    \end{prop}
    \begin{proof}
      First, we assume that $\TheGeometricComplex$ has finite dimension.
      Let $\GeometricCell$ be a maximal cell of $\TheGeometricComplex$.
      If $\GeometricCell\subseteq\BadSet$, the cell $\GeometricCell$
      does not intersect
      $\TheGeometricComplex\setminus\BadSet$ and we do not need
      to do anything.
      If $\GeometricCell$ avoids $\BadSet$, the map $\TheDeformation$ must
      be the identity on $\GeometricCell$.
      Otherwise, let $\TheRadialCenter$
      be a point in the intersection $\GeometricCell\intersect\BadSet$
      chosen in the relative interior of $\GeometricCell$.
      Projecting away from $\TheRadialCenter$, as in
      Observation~\ref{bernd:deformation.retraction},
      deformation retracts $\GeometricCell\setminus\BadSet$ onto
      $\BoundaryOf{\GeometricCell}\setminus\BadSet$. The maps constructed
      for two maximal cells agree on their intersection.
      Hence we can paste all these maps together to get a deformation
      retraction of $\TheGeometricComplex\setminus\BadSet$ onto
      $\AltGeometricComplex\setminus\BadSet$ where $\AltGeometricComplex$ is
      $\TheGeometricComplex$ with the interiors of all maximal
      cells intersecting $\BadSet$ removed.

      Now, $\AltGeometricComplex$ has other maximal cells, which
      might intersect $\BadSet$. Using the same construction
      for $\AltGeometricComplex$, we obtain another deformation retraction
      $\AltGeometricComplex\setminus\BadSet
      \rightarrow\ThrGeometricComplex\setminus\BadSet$. We keep going,
      removing more and more cells intersecting $\BadSet$.
      Since the dimension of $\TheGeometricComplex$ is finite, the process
      terminates
      after finitely many steps. The composition of the maps
      thus obtained is the strong deformation retraction from
      $\TheGeometricComplex\setminus\BadSet$ onto
      $\TheGeometricComplex[][\BadSet]$.
      This proves the claim for finite dimensional $\TheGeometricComplex$.

      Note that the construction is local: what it does on a cell
      is only determined by the intersection of this cell with the set
      $\BadSet$.
      Hence, the deformation retraction is compatible with
      subcomplexes. More precisely, if $\TheGeometricSubcomplex$ is a subcomplex of
      $\TheGeometricComplex$, then the deformation retractions
      $\TheDeformation[\TheGeometricComplex]$ and
      $\TheDeformation[\TheGeometricSubcomplex]$ from above
      are constructed
      such that $\TheDeformation[\TheGeometricSubcomplex]$ is the restriction
      of
      $\TheDeformation[\TheGeometricComplex]$ to $\TheGeometricSubcomplex$.
      It follows that the
      pair
      $(\TheGeometricComplex\setminus\BadSet,\TheGeometricSubcomplex\setminus\BadSet)$
      is homotopy equivalent
      to $(\TheGeometricComplex[][\BadSet],\TheGeometricSubcomplex[][\BadSet])$.
      Applying this observation to pairs of skeleta, the
      claim follows by standard arguments in the case that
      $\TheGeometricComplex$ has infinite dimension.
    \end{proof}

    Let $\TheSphBuild$ be a spherical building. We regard $\TheSphBuild$
    as a metric space with the angular metric. So each apartment is
    a round sphere of radius~$\One$. When $\TheSphBuild$ is a
    finite building,
    the topology induced by the metric agrees with the weak topology
    it carries as a simplicial complex. For locally infinite buildings, both
    topologies differ and we will use the weak topology throughout
    for the building and all its subcomplexes.
    \begin{prop}\label{bernd:complement}
      Let $\TheSphBuild$ be a spherical building and fix a chamber
      $\TheSphChamber$ in $\TheSphBuild$.
      Let $\BadSet\subset\TheSphBuild$
      be a subset such that, for any apartment $\TheSphApp$
      containing $\TheSphChamber$ the intersection
      $\BadSet\intersect\TheSphApp$ is a proper, open, and
      convex subset of the sphere $\TheSphApp$. Then the
      space $\TheSphComplement:=\TheSphBuild\setminus\BadSet$ and its
      maximal subcomplex $\TheSphBuild[][\BadSet]$ are both
      $(\DimOf{\TheSphBuild}-\One)$-connected. The complex
      $\TheSphBuild[][\BadSet]$ has dimension $\DimOf{\TheSphBuild}$
      and hence is spherical of this dimension.
    \end{prop}
    \begin{rem}
      Using $\BadSet=\emptyset$ in Proposition~\ref{bernd:complement},
      we obtain the Solomon-Tits Theorem as a special case.
      Theorem~A of \cite{Schulz:2010}, whose proof
      inspired the argument given below, is the special case where
      $\BadSet$ is open, convex, and of diameter strictly less
      than $\pi$.
    \end{rem}
    \begin{proof}[of Proposition~\ref{bernd:complement}]
      We observe first that Proposition~\ref{bernd:push_in}
      implies that
      the subset $\TheSphComplement$ and its maximal subcomplex
      $\TheSphBuild[][\BadSet]$ are homotopy equivalent.
      Therefore, it suffices to prove that $\TheSphComplement$ is
      $(\DimOf{\TheSphBuild}-\One)$-connected.

      We have to contract spheres of dimensions up to
      $\DimOf{\TheSphBuild}-\One$. Let
      $\TheSphere\subseteq\TheSphComplement$ be such a sphere.
      Since $\TheSphere$ is compact in $\TheSphBuild$, it is
      covered by a finite family of apartments and we can
      apply \cite[Lemma~3.5]{Heydebreck:2003}: there is a finite
      sequence $\TheSphApp[\One],\TheSphApp[\Two],\ldots,
      \TheSphApp[\TheLastIndex]$ such that
      (a) each $\TheSphApp[\TheIndex]$ contains $\TheSphChamber$,
      (b) the sphere $\TheSphere$ is contained in the union
      $\Union[\TheIndex]{\TheSphApp[\TheIndex]}$, and most
      importantly,
      (c) for each $\TheIndex\geq\Two$ the intersection
      $\TheSphApp[\TheIndex]\intersect
      (\TheSphApp[\One]\union\cdots
      \union\TheSphApp[\TheIndex-\One])$
      is a union of closed half-apartments, each of which contains
      $\TheSphChamber$.
      Put
      \(
        \TheGeometricComplex[\TheIndex] :=
        \TheSphApp[\One] \union\cdots\union
        \TheSphApp[\TheIndex]
      \)
      and observe that $\TheGeometricComplex[\TheIndex]$ is obtained
      from $\TheGeometricComplex[\TheIndex-\One]$ by gluing in the
      closure
      \(
        \NewPart[\TheIndex] :=
        \ClosureOf{
          \TheSphApp[\TheIndex]\setminus
          (\TheSphApp[\One]\union\cdots
          \union\TheSphApp[\TheIndex-\One])
        }
      \)
      along the boundary $\BoundaryOf{\NewPart[\TheIndex]}$
      of $\NewPart[\TheIndex]$ in $\TheSphApp[\TheIndex]$.
      Note that $\NewPart[\TheIndex]$ is fat and demi-convex.

      Now, we can build $\TheGeometricComplex[\TheLastIndex]\setminus\BadSet$
      inductively. We begin with $\TheGeometricComplex[\One]\setminus\BadSet$,
      which is contractible. The space
      $\TheGeometricComplex[\TheIndex]\setminus\BadSet$ is obtained
      from
      $\TheGeometricComplex[\TheIndex-\One]\setminus\BadSet$ by gluing
      in $\NewPart[\TheIndex]\setminus\BadSet$ along
      $\BoundaryOf{\NewPart[\TheIndex]}\setminus\BadSet$.
      If $\NewPart[\TheIndex]$ and $\BadSet$ are disjoint,
      this is a cellular
      extension of dimension $\DimOf{\TheSphBuild}$
      as $\NewPart[\TheIndex]$ is fat. Otherwise,
      Observation~\ref{bernd:deformation.retraction} implies that
      $\NewPart[\TheIndex]\setminus\BadSet$ deformation retracts
      onto $\BoundaryOf{\NewPart[\TheIndex]}\setminus\BadSet$, whence
      $\TheGeometricComplex[\TheIndex]\setminus\BadSet$ and
      $\TheGeometricComplex[\TheIndex-\One]\setminus\BadSet$ are
      homotopy equivalent in this case. In the end,
      the sphere $\TheSphere$ can be contracted inside
      $\TheGeometricComplex[\TheLastIndex]\setminus\BadSet$.
    \end{proof}
    \begin{cor}\label{bernd:complement:closed}
      Let $\TheSphBuild$ be a \emph{finite}
      spherical building and fix a chamber
      $\TheSphChamber$ in $\TheSphBuild$.
      Let $\ClosedBadSet\subset\TheSphBuild$
      be a subset such that, for any apartment $\TheSphApp$
      containing $\TheSphChamber$ the intersection
      $\ClosedBadSet\intersect\TheSphApp$ is a closed
      convex subset of diameter strictly less than $\pi$
      in the sphere $\TheSphApp$. Then the
      space $\TheSphComplement:=\TheSphBuild\setminus\ClosedBadSet$ and its
      maximal subcomplex $\TheSphBuild[][\ClosedBadSet]$ are both
      $(\DimOf{\TheSphBuild}-\One)$-connected. The complex
      $\TheSphBuild[][\ClosedBadSet]$ has dimension $\DimOf{\TheSphBuild}$
      and hence is spherical of this dimension.
    \end{cor}
    \begin{proof}
      The building $\TheSphBuild$ is finite, hence $\ClosedBadSet$
      is compact. Let $\BadSet$ be an $\TheEpsilon$-neighborhood of
      $\ClosedBadSet$. Choosing $\TheEpsilon$ sufficiently small,
      we can ensure that $\BadSet$ satisfies the hyptheses of
      Proposition~\ref{bernd:complement}, that
      \(
        \TheSphBuild\setminus\BadSet
      \)
      and
      \(
        \TheSphBuild\setminus\ClosedBadSet
      \)
      are homotopy equivalent, and that
      \(
        \TheSphBuild[][\BadSet]
        =
        \TheSphBuild[][\ClosedBadSet]
      \)
    \end{proof}

    An interesting special case of Proposition~\ref{bernd:complement},
    also already noted in
    \cite{Schulz:2010}, is obtained when $\BadSet$ is chosen
    as the open $\frac{\pi}{\Two}$-ball around a fixed point
    $\NorthPole\in\TheSphBuild$, which we think of as the
    \notion{north pole}. Then the complex
    $\TheSphBuildOf[][\geq\frac{\pi}{\Two}]{\NorthPole} :=
    \TheSphBuild[][\BadSet]$
    is a \notion{closed hemisphere complex} and
    $\DimOf{\TheSphBuild}$-spherical by Proposition~\ref{bernd:complement}.
    The argument fails badly if $\BadSet$ is chosen as the
    \emph{closed ball} of radius $\frac{\pi}{\Two}$ around
    $\NorthPole$. In fact, the \notion{open hemisphere complex}
    $\TheSphBuildOf[][>\frac{\pi}{\Two}]{\NorthPole}$ spanned by
    all vertices avoiding the closed ball $\BadSet$ generally is
    not $\DimOf{\TheSphBuild}$-spherical: the dimension of
    $\TheSphBuildOf[][>\frac{\pi}{\Two}]{\NorthPole}$ might be too
    small. The main result of Schulz is that this is the
    only obstruction.
    \begin{prop}[{see \cite[Theorems~A and~B]{Schulz:2010}}]\label{bernd:open.hemi}
      The open hemisphere complex
      $\TheSphBuildOf[][>\frac{\pi}{\Two}]{\NorthPole}$
      is spherical of dimension
      $\DimOf{\TheSphBuild[\ver]}$.
      If $\TheSphBuild$ is thick,
      then neither open nor closed hemisphere complexes in
      $\TheSphBuild$ are contractible.
    \end{prop}
    The subcomplex $\TheSphBuildOf[\ver]{\NorthPole}$ is defined as
    follows: The \notion{equator}
    $\TheSphBuildOf[][=\frac{\pi}{\Two}]{\NorthPole}$ is the
    subcomplex spanned by
    those points in $\TheSphBuild$ of distance $\frac{\pi}{\Two}$
    from $\NorthPole$. Recall that
    $\TheSphBuild$ decomposes as a join of unique irreducible
    factors. The \notion{horizontal part}
    $\TheSphBuildOf[\hor]{\NorthPole}$
    is the join of all factors fully contained in the equator.
    The complex $\TheSphBuildOf[\ver]{\NorthPole}$ is the join of the
    other irreducible factors. In particular,
    \begin{equation}\label{eq:join.dec}
      \TheSphBuild = \TheSphBuildOf[\hor]{\NorthPole}
      \join
      \TheSphBuildOf[\ver]{\NorthPole}.
    \end{equation}

  \section{Height}\label{sec:height}
    We now begin the proof of the Rank~Theorem proper.
    Let $\TheModelApp$ be a euclidean Coxeter complex upon which
    the apartments of $\TheEuBuild$ are modeled, and let
    $\TheEuSpace$ be the underlying euclidean space where the
    origin $\SpZero$ shall correspond to a special vertex in
    $\TheModelApp$. Let $\TheFullWeyl$ be the spherical Weyl
    group generated by the walls of $\TheModelApp$ through
    $\SpZero$. For constructing zonotopes, we shall choose an
    \notion{admissible} subset $\TheZonoDef\subset\TheEuSpace$,
    i.e., we require that $\TheZonoDef$ be finite,
    $\TheFullWeyl$-invariant, and symmetric with respect to
    the origin $\SpZero$. In the course of the argument, we
    will need to strengthen the requirements on $\TheZonoDef$,
    but we begin with any admissible $\TheZonoDef$.
    Since $\TheZonoDef$ is invariant with respect to the maximal
    Weyl group, the subset
    \(
      \TheEuPoint+\TheZonotopeOf{\TheZonoDef}
    \)
    is well-defined in any apartment $\TheEuApp$ of
    $\TheEuBuild$ containing $\TheEuPoint$.

    \begin{figure}[ht]
      \null\hfill
      \setlength{\unitlength}{1cm}
      \begin{picture}(9,7.3)
        \put(-3.5,0){\includegraphics{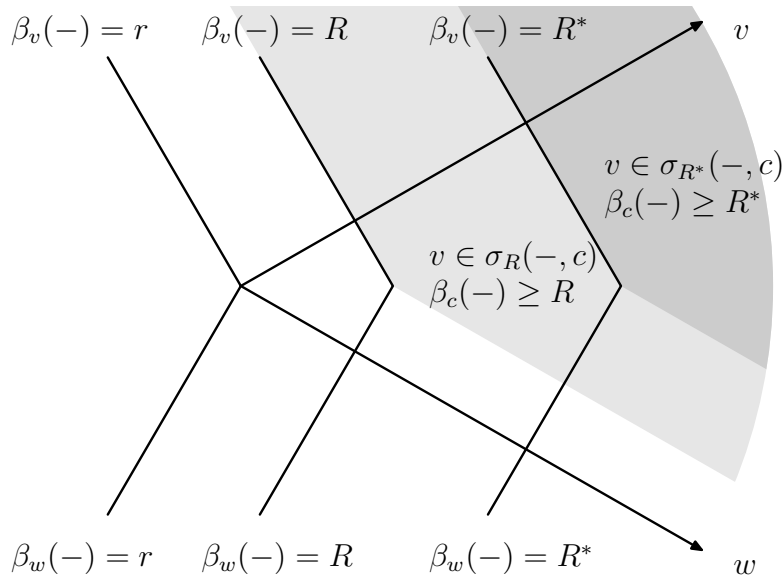}}
        \put(8,7){$\TheRatVertex$}
        \put(8,-0.1){$\AltRatVertex$}
        \put(-1.5,7){$\TheBusemannOf[\TheRatVertex]{\DummyArg}=\TheLowerBound$}
        \put(-1.5,0){$\TheBusemannOf[\AltRatVertex]{\DummyArg}=\TheLowerBound$}
        \put(1,7){$\TheBusemannOf[\TheRatVertex]{\DummyArg}=\TheUpperBound$}
        \put(1,0){$\TheBusemannOf[\AltRatVertex]{\DummyArg}=\TheUpperBound$}
        \put(4,7){$\TheBusemannOf[\TheRatVertex]{\DummyArg}=\TheHugeBound$}
        \put(4,0){$\TheBusemannOf[\AltRatVertex]{\DummyArg}=\TheHugeBound$}
        \put(4,3.5){$\TheBusemannOf[\TheRatChamber]{\DummyArg}\geq\TheUpperBound$}
        \put(4,4){$\TheRatVertex\in\CloseSimplexOf[\TheUpperBound]{\DummyArg,\TheRatChamber}$}
        \put(6.3,4.7){$\TheBusemannOf[\TheRatChamber]{\DummyArg}\geq\TheHugeBound$}
        \put(6.3,5.2){$\TheRatVertex\in\CloseSimplexOf[\TheHugeBound]{\DummyArg,\TheRatChamber}$}
      \end{picture}
      \hfill\null
      \caption{moving out the tip\label{height:05.fig}}
      {\small
      Given a fixed distance bound $\TheHausdorffDist$, we can choose
      $\TheHugeBound$ large enough so that
      \(
        \TheBusemannOf[\TheRatChamber]{\TheEuPoint}
        \geq
        \TheHugeBound
      \)
      and
      \(
        \TheRatVertex\in
        \CloseSimplexOf[\TheHugeBound]{\TheEuPoint,\TheRatChamber}
      \)
      implies
      \(
        \TheRatVertex\in
        \CloseSimplexOf[\TheUpperBound]{\AltEuPoint,\TheRatChamber}
      \)
      for any two points $\TheEuPoint$ and $\AltEuPoint$ of
      distance at most $\TheHausdorffDist$:
      every reduced point $\TheHausdorffDist$-close to the darker
      area still lies within the lightly shaded region.}
    \end{figure}
    \begin{observation}\label{height:second.level}
      Given a reduction datum
      \(
        \TupelOf{
          \FamOf{\TheBusemann[\DummyVertex]},
          \TheLowerBound,
          \TheUpperBound
        }
      \)
      for an absolutely almost simple group,
      there exists a bound $\TheHugeBound\geq\TheUpperBound$ such
      that the following implication holds:
      \[
        \TheBusemannOf[\TheRatChamber]{\TheEuPoint}
        \geq\TheHugeBound
        \text{\ and\ }
        \AltEuPoint\in\TheEuPoint+\TheZonotopeOf{\TheZonoDef}
        \qquad\Longrightarrow\qquad
        \CloseSimplexOf[\TheHugeBound]{\AltEuPoint,\TheRatChamber}
        \subseteq
        \CloseSimplexOf[\TheUpperBound]{\TheEuPoint,\TheRatChamber}
      \]
      Here, $\TheRatChamber\in\ChambersOf{\TheRatBuild}$ is any
      chamber and $\TheEuPoint,\AltEuPoint\in\TheEuBuild$ are
      any two points that lie in a common
      apartment whose visual boundary contains $\TheRatChamber$.

      See Figure~\ref{height:05.fig} for a picture that explains
      how to choose $\TheHugeBound$ for a single vertex.
      Since there are only finitely many types of vertices in the
      building $\TheRatBuild$, one can choose $\TheHugeBound$ large
      enough for every type.\qed
    \end{observation}

    We modify the construction given in Section~\ref{sec:blueprint}.
    Let
    \(
      \TupelOf{
        \FamOf{\TheBusemann[\DummyVertex]},
        \TheLowerBound,
        \TheUpperBound
      }
    \)
    be a reduction datum that reduces chambers in $\TheEuBuild$
    uniformly. Fix $\TheHugeBound$ as in
    Observation~\ref{height:second.level}.
    For a point $\TheEuPoint\in\TheEuBuild$, an apartment $\TheEuApp$
    containing $\TheEuPoint$, and a chamber
    $\TheRatChamber\in\ChambersOf{\TheRatBuild}$ in the visual
    boundary $\BoundaryOf{\TheEuApp}$ and reducing $\TheEuPoint$,
    we define
    \[
      \TheMorseOf[\TheEuApp,\TheRatChamber]{\TheEuPoint}
      :=
      \DistanceOf{
        \TheEuPoint+\TheZonotopeOf{\TheZonoDef}
       }{
        \ConvexChunkOf[\TheEuApp,\TheRatChamber]{\TheHugeBound}
      }
      =
      \DistanceOf{
        \TheEuPoint
       }{
        \ConvexChunkOf[\TheEuApp,\TheRatChamber]{\TheHugeBound}
        -
        \TheZonotopeOf{\TheZonoDef}
      }
    \]
    Let
    \(
      \TheEuPoint[\TheEuApp,\TheRatChamber][*]
    \)
    be the point in
    \(
        \ConvexChunkOf[\TheEuApp,\TheRatChamber]{\TheHugeBound}
        -
        \TheZonotopeOf{\TheZonoDef}
    \)
    closest to $\TheEuPoint$, and for
    \(
      \TheEuPoint\not\in
        \ConvexChunkOf[\TheEuApp,\TheRatChamber]{\TheHugeBound}
        -
        \TheZonotopeOf{\TheZonoDef}
    \)
    let
    \(
      \TheEndOf[\TheEuApp,\TheRatChamber]{\TheEuPoint}
    \)
    be the visual limit of the ray from
    \(
      \TheEuPoint[\TheEuApp,\TheRatChamber][*]
    \)
    through
    \(
      \TheEuPoint.
    \)
    \begin{prop}\label{height:well.put}
      Let $\TheEuPoint\in\TheEuBuild$ be a point. Let $\TheEuApp$
      and $\AltEuApp$ be two apartments containing $\TheEuPoint$,
      let $\TheRatChamber$ and $\AltRatChamber$ be two chambers
      of $\TheRatBuild$ reducing $\TheEuPoint$. Assume that
      $\TheRatChamber\subset\BoundaryOf{\TheEuApp}$ and
      $\AltRatChamber\subset\BoundaryOf{\AltEuApp}$.
      Then
      \(
        \TheMorseOf[\TheEuApp,\TheRatChamber]{\TheEuPoint}
        =
        \TheMorseOf[\AltEuApp,\AltRatChamber]{\TheEuPoint}
      \).
      Moreover,
      \(
        \TheEndOf[\TheEuApp,\TheRatChamber]{\TheEuPoint}
        =
        \TheEndOf[\AltEuApp,\AltRatChamber]{\TheEuPoint}
        \in
        \CloseSimplexOf{\TheEuPoint}
      \)
      provided $\TheMorseOf[\TheEuApp,\TheRatChamber]{\TheEuPoint}>\Zero$.
    \end{prop}
    \begin{proof}
      Assume first $\TheBusemannOf[\TheRatChamber]{\TheEuPoint}\leq
      \TheHugeBound$. Since $\TheHugeBound\geq\TheUpperBound$,
      Corollary~\ref{behr:low} implies
      \(
        \TheBusemannOf[\AltRatChamber]{\TheEuPoint}
        \leq
        \TheHugeBound
      \).
      Hence
      \[
        \TheMorseOf[\TheEuApp,\TheRatChamber]{\TheEuPoint}
        =
        \Zero
        =
        \TheMorseOf[\AltEuApp,\AltRatChamber]{\TheEuPoint}
        \Period
      \]

      It remains to argue the case
      $\TheBusemannOf[\TheRatChamber]{\TheEuPoint}>
      \TheHugeBound\geq\TheUpperBound$.
      First, we work inside $\TheEuApp$.
      Choose $\AltEuPoint\in\TheEuPoint+\TheZonotopeOf{\TheZonoDef}$
      so that it minimizes the distance to
      $\ConvexChunkOf[\TheEuApp,\TheRatChamber]{\TheHugeBound}$,
      and let $\AltEuPoint[\TheEuApp,\TheRatChamber]$ the point
      in $\ConvexChunkOf[\TheEuApp,\TheRatChamber]{\TheHugeBound}$
      closest to $\AltEuPoint$. By choice of $\TheHugeBound$,
      Observation~\ref{height:second.level} applies whence
      \(
        \CloseSimplexOf[\TheHugeBound]{\AltEuPoint,\TheRatChamber}
        \subseteq
        \CloseSimplexOf[\TheUpperBound]{\TheEuPoint,\TheRatChamber}
        =
        \CloseSimplexOf{\TheEuPoint}
      \).
      It follows from Observation~\ref{behr:restricting.simplex}
      that $\AltEuPoint[\TheEuApp,\TheRatChamber]$ is also
      the point in
      $\ConvexChunkOf[\TheEuApp,\CloseSimplexOf{\TheEuPoint}]{\TheHugeBound}$
      closest to $\AltEuPoint$.

      Now the isometry argument applies: since
      \(
        \TheBusemannOf[\TheRatChamber]{\TheEuPoint}\geq
        \TheHugeBound\geq\TheUpperBound
      \)
      we have $\CloseSimplexOf{\TheEuPoint}\subseteq
      \TheRatChamber\intersect\AltRatChamber$. Hence there is a
      Coxeter isomorphism
      $\TheIsometry\mapcolon\TheEuApp\rightarrow\AltEuApp$
      fixing $\TheEuPoint$ and $\CloseSimplexOf{\TheEuPoint}$.
      Since
      \(
        \ConvexChunkOf[\TheEuApp,\CloseSimplexOf{\TheEuPoint}]{
          \TheHugeBound
        }
      \)
      is only defined in terms of Busemann functions indexed by
      vertices in $\CloseSimplexOf{\TheEuPoint}$, the isometry
      $\TheIsometry$ identifies
      \(
        \ConvexChunkOf[\TheEuApp,\CloseSimplexOf{\TheEuPoint}]{
          \TheHugeBound
        }
      \)
      with
      \(
        \ConvexChunkOf[\AltEuApp,\CloseSimplexOf{\TheEuPoint}]{
          \TheHugeBound
        }
      \).
      As $\TheIsometry$ is a Coxeter isomorphism, it identifies
      the two sets $\TheEuPoint+\TheZonotopeOf{\TheZonoDef}$
      as drawn in $\TheEuApp$ and $\AltEuApp$. It follows that
      \begin{eqnarray*}
        \TheMorseOf[\TheEuApp,\TheRatChamber]{\TheEuPoint}
        &=&
        \DistanceOf[\TheEuApp]{
          \TheEuPoint+\TheZonotopeOf{\TheZonoDef}
        }{
          \ConvexChunkOf[\TheEuApp,\CloseSimplexOf{\TheEuPoint}]{
            \TheHugeBound
          }
        }
        \\
        &=&
        \DistanceOf[\AltEuApp]{
          \TheEuPoint+\TheZonotopeOf{\TheZonoDef}
        }{
          \ConvexChunkOf[\AltEuApp,\CloseSimplexOf{\TheEuPoint}]{
            \TheHugeBound
          }
        }
        \\
        &=&
        \TheMorseOf[\AltEuApp,\AltRatChamber]{\TheEuPoint}
        \text{.}
      \end{eqnarray*}

      If $\TheMorseOf[\TheEuApp,\TheRatChamber]{\TheEuPoint}>\Zero$,
      then
      \(
        \AltEuPoint[\TheEuApp,\TheRatChamber]
        \neq
        \AltEuPoint
      \)
      and the ray from
      $\AltEuPoint[\TheEuApp,\TheRatChamber]$ through
      $\AltEuPoint$ is parallel to the ray
      from $\TheEuPoint[\TheEuApp,\TheRatChamber][*]$
      through $\TheEuPoint$. Hence, it defines the same visual
      end, which lies in
      $\CloseSimplexOf[\TheHugeBound]{\AltEuPoint,\TheRatChamber}
      \subseteq\CloseSimplexOf[\TheUpperBound]{\TheEuPoint,\TheRatChamber}
      =\CloseSimplexOf{\TheEuPoint}$.
      The isometry $\TheIsometry$ identifies the ray
      from
      $\AltEuPoint[\TheEuApp,\TheRatChamber]$ through
      $\AltEuPoint$
      with its counter part in $\AltEuApp$. Hence,
      $\TheEndOf[\TheEuApp,\TheRatChamber]{\TheEuPoint}
      =\TheEndOf[\AltEuApp,\AltRatChamber]{\TheEuPoint}$.
    \end{proof}
    \begin{observation}\label{height:gradient.in.app}
      Note that the ray
      from
      $\AltEuPoint[\TheEuApp,\TheRatChamber]$ through
      $\AltEuPoint$
      gives the direction of fastest ascent for the
      function $\TheMorse[\TheEuApp,\TheRatChamber]$ in
      the point $\TheEuPoint$. Also, moving $\TheEuPoint$
      in that direction increases
      $\TheMorseOf[\TheEuApp,\TheRatChamber]{\TheEuPoint}$
      with unit speed and the gradient of
      $\TheMorse[\TheEuApp,\TheRatChamber]$ does not change
      along this ray.\qed
    \end{observation}
    Hence, we can define
    \(
      \TheMorseOf{\TheEuPoint}
      :=
      \TheMorseOf[\TheEuApp,\TheRatChamber]{\TheEuPoint}
    \)
    and
    \(
      \TheEndOf{\TheEuPoint}
      :=
      \TheEndOf[\TheEuApp,\TheRatChamber]{\TheEuPoint}
    \).
    Here $\TheEuApp$ is any apartment of $\TheEuBuild$ containing
    $\TheEuPoint$ and
    $\TheRatChamber$ is any chamber in $\TheRatBuild$ lying in
    $\BoundaryOf{\TheEuApp}$ and reducing $\TheEuPoint$.

    \begin{observation}\label{height:invariant}
      Since the reduction datum used in the construction is
      $\TheLattice$-invariant, so is the function $\TheMorse$.\qed
    \end{observation}
    \begin{observation}\label{height:cocompact}
      There is a constant $\TheShiftConst$, depending on
      $\TheZonoDef$ and $\TheHugeBound$, such that
      \(
        \RoughMorseOf{\TheEuPoint}\leq
        \TheMorseOf{\TheEuPoint}+\TheShiftConst
      \)
      for each $\TheEuPoint\in\TheEuBuild$.
      Hence Proposition~\ref{blueprint:cocompact} implies that
      each sublevel set
      \(
        \SubLevelSet{\TheMorse}{\TheBound}
        \subseteq
        \SubLevelSet{\RoughMorse}{\TheBound+\TheShiftConst}
      \)
      has compact quotient modulo the action of $\TheLattice$.\qed
    \end{observation}
    As for continuity and the gradient field, nothing essential
    changes.
    \begin{observation}\label{height:continuous}
      The same reasoning as in the proof of
      Proposition~\ref{blueprint:continuous}
      shows that
      the function $\TheMorse$ is continuous.\qed
    \end{observation}
    \begin{prop}\label{height:gradient}
      Let $\TheEuPoint$ and $\AltEuPoint$ be two points in
      $\TheEuBuild$ that lie in a common closed chamber of
      $\TheEuBuild$. Then
      \[
        \TheMorseOf{\AltEuPoint}
        -
        \TheMorseOf{\TheEuPoint}
        \leq
        \DistanceOf{\AltEuPoint}{\TheEuPoint}
      \]
      with equality if and only if $\AltEuPoint$ lies on the
      ray $\RayFromTo{\TheEuPoint}{\TheEndOf{\TheEuPoint}}$.
      In case of equality, moreover
      \(
        \TheEndOf{\AltEuPoint}
        =
        \TheEndOf{\TheEuPoint}
      \).
    \end{prop}
    \begin{proof}
      By uniformity of the reduction datum, choose
      $\TheRatChamber$ and $\TheEuApp$ so that
      $\TheMorse$ and $\TheMorse[\TheEuApp,\TheRatChamber]$
      agree on the segment $\LineSegment{\TheEuPoint}{\AltEuPoint}$.
      Then
      \begin{eqnarray*}
        \TheMorseOf{\AltEuPoint}
        -
        \TheMorseOf{\TheEuPoint}
        &=&
        \TheMorseOf[\TheEuApp,\TheRatChamber]{\AltEuPoint}
        -
        \TheMorseOf[\TheEuApp,\TheRatChamber]{\TheEuPoint}
        \\
        &=&
        \DistanceOf{\AltEuPoint}{
          {\ConvexChunkOf[\TheEuApp,\TheRatChamber]{\TheHugeBound}
          -\TheZonotopeOf{\TheZonoDef}}}
        -
        \DistanceOf{\TheEuPoint}{
          {\ConvexChunkOf[\TheEuApp,\TheRatChamber]{\TheHugeBound}
          -\TheZonotopeOf{\TheZonoDef}}}
        \\
        &\leq&
        \DistanceOf{\AltEuPoint}{\TheEuPoint}
        \Period
      \end{eqnarray*}
      By Observation~\ref{height:gradient.in.app}, we have
      equality if $\AltEuPoint\in
      \RayFromTo{\TheEuPoint}{\TheEndOf{\TheEuPoint}}$,
      and in this case
      $\TheEndOf{\TheEuPoint}=\TheEndOf{\AltEuPoint}$.
    \end{proof}
    For $\TheEuPoint\in\TheEuBuild$, we define the gradient
    $\GradOf[\TheEuPoint]{\TheMorse}\in\LinkOf{\TheEuPoint}$
    to be the direction defined by the geodesic ray
    $\RayFromTo{\TheEuPoint}{\TheEndOf{\TheEuPoint}}$.
    Along this ray, the function
    $\TheMorse$ increases with unit speed and all other directions
    show a slower increase. Thus, the geodesic ray
    $\RayFromTo{\TheEuPoint}{\TheEndOf{\TheEuPoint}}$
    can be regarded as the flow line of the gradient field
    $\GradOf{\TheMorse}$ starting at $\TheEuPoint$.

    Let us call a \notion{brick} any subset of $\TheEuBuild$ that
    arises as the
    convex hull of a set of vertices of a common chamber
    in $\TheEuBuild$.
    \begin{observation}\label{height:orthogonal.gradient}
      Let $\TheEuPoint\in\TheEuBuild$
      be a point in a brick $\TheEuBrick$ such that
      $\GradOf[\TheEuPoint]{\TheMorse}$ is perpendicular to
      $\TheEuBrick$. Then $\TheEuPoint$ is a point of lowest
      height in $\TheEuBrick$.
    \end{observation}
    \begin{proof}
      Choose $\TheEuApp$ and $\TheRatChamber$ so that
      $\TheMorse$ agrees
      with $\TheMorse[\TheEuApp,\TheRatChamber]$
      on $\TheEuBrick$. Hence, $\TheMorse$ is a convex function
      on $\TheEuBrick$ and the claim follows.
    \end{proof}

    To actually ensure that $\TheMorse$ is superior to
    $\RoughMorse$, we have to strengthen the requirement on
    $\TheZonoDef$. Of course, we have to adjust $\TheHugeBound$
    accordingly.

    We call $\TheZonoDef$ \notion{almost rich} if for any
    two vertices $\TheSpVertex$ and $\AltSpVertex$ of
    $\TheModelApp$ that belong to a common chamber, the
    difference $\TheSpVertex-\AltSpVertex\in\TheZonoDef$.
    Note that one can obtain a finite, admissible, almost rich
    set $\TheZonoDef$ by starting with the finite set of
    difference vectors arising from the vertices of a fixed
    chamber (note that this is automatically symmetric with
    respect to the origin) and then closing the set with respect
    to the action of $\TheFullWeyl$: since $\SpZero$ is a special
    vertex, $\TheFullWeyl$ acts transitively on parallelism classes
    of chambers in $\TheModelApp$.
    \begin{figure}[ht]
      \null\hfill
      \includegraphics{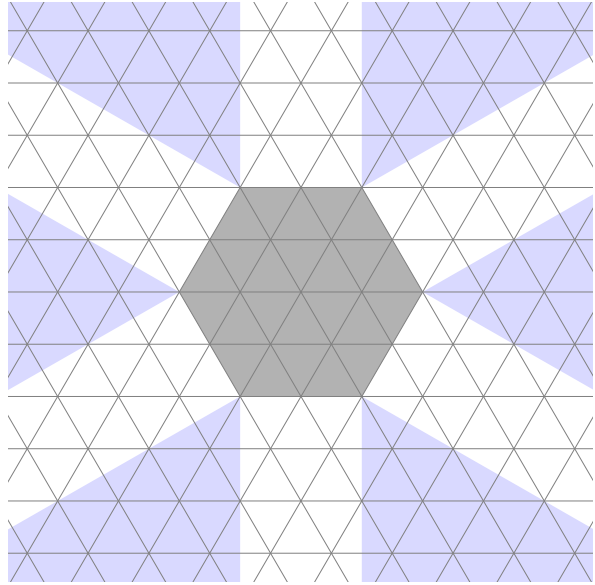}\hfill\null
      \caption{an almost rich zonotope
      for $\tilde{\text{\textsf{A}}}_2$\label{fig:almost.rich.zono}}
      {\small
      The lightly shaded areas and the white corridors are the
      normal cones for the zonotope.}
    \end{figure}

    When $\TheZonoDef$ is almost rich,
    the results of Section~\ref{sec:zonotopes} apply to bricks.
    Let us spell out the consequences of
    Proposition~\ref{zonotopes:minset.maxset} combined with
    Example~\ref{zonotopes:distance.to.body}.
    \begin{cor}[to Proposition~\ref{zonotopes:minset.maxset}]\label{height:minmax}
      If $\TheZonoDef$ is almost rich, the function $\TheMorse$
      assumes its minimum on any brick in $\TheEuBuild$ in a
      vertex, i.e., the subset of point of minimum height contains
      a vertex. The subset of points of maximum height, on the other
      hand, is a face of the brick.\qed
    \end{cor}

    \begin{prop}[Gradient Criterion]\label{height:criterion}
      We still assume that $\TheZonoDef$ is almost rich.
      Let $\TheEuPoint$ and $\AltEuPoint$ be two distinct vertices in
      $\TheEuBuild$ that lie in a common chamber. Then the following
      hold:
      \begin{enumerate}
        \item The function $\TheMorse$ is monotonic on the
          line segment
          $\TheSegment := \LineSegment{\TheEuPoint}{\AltEuPoint}$.
        \item
          The angle
          \(
            \AngleOf[\TheEuPoint]{\TheSegment}{
              \GradOf[\TheEuPoint]{\TheMorse}}
            >\PiHalf
          \)
          if and only if
          \(
            \TheMorseOf{\AltEuPoint} < \TheMorseOf{\TheEuPoint}
          \).
      \end{enumerate}
    \end{prop}
    \begin{proof}
      By Corollary~\ref{height:minmax}, $\TheMorse$ attains its
      minimum along $\TheSegment$ at a boundary point. As $\TheMorse$
      is a convex function, this proves the first claim.

      If
      \(
        \AngleOf[\TheEuVertex]{\TheSegment}{
          \GradOf[\TheEuVertex]\TheMorse
        }
        \neq\PiHalf
      \)
      the height $\TheMorse$ changes when one moves from the vertex
      $\TheEuVertex$ infinitesimally into the segment $\TheSegment$.
      If the
      angle is obtuse, the height decreases; if the angle is
      acute, the height increases. By monotonicity, $\TheEuVertex$ must
      be the highest or lowest point on $\TheSegment$, respectively.
      Observation~\ref{height:orthogonal.gradient} covers the remaining
      case that
      $\GradOf[\TheEuVertex]{\TheMorse}$ is orthogonal to $\TheSegment$.
    \end{proof}

  \section{Cells of constant height}\label{sec:flat}
    Let $\TheEuPoint\in\TheEuBuild$.
    We think of the link $\LinkOf{\TheEuPoint}$
    as the space of directions issuing from $\TheEuPoint$.
    It is a spherical building and we regard it as a metric
    space via the angular metric.

    Now suppose that $\TheEuPoint$ is carried by the cell
    $\HorEuCell$. The link $\LinkOf{\HorEuCell}$ is also a spherical
    building. Its simplicial structure corresponds to the poset
    of cofaces of $\HorEuCell$ in $\TheEuBuild$.
    We can realize $\LinkOf{\HorEuCell}$ as the space of directions
    at $\TheEuPoint$ orthogonal to $\HorEuCell$. For two different
    points carried by $\HorEuCell$, the corresponding realizations of
    $\LinkOf{\HorEuCell}$ are canonically identified and
    we may think of elements in $\LinkOf{\HorEuCell}$ as
    parallel fields of directions perpendicular to $\HorEuCell$.
    This way, $\LinkOf{\HorEuCell}$ carries an angular metric.
    We thus consider $\LinkOf{\HorEuCell}$ as a metric space.

    The point link splits as a spherical join
    \begin{equation}\label{eq:point.link.split}
      \LinkOf{\TheEuPoint}
      =
      \BoundaryOf{\HorEuCell}
      \join
      \LinkOf{\HorEuCell}
    \end{equation}
    where $\BoundaryOf{\HorEuCell}$ is the subspace of
    $\LinkOf{\TheEuPoint}$ consisting of those directions that
    do not leave $\HorEuCell$. It is a round sphere in the
    angular metric. As $\HorEuCell$ is a poly-simplicial
    cell, there is also an obvious cell structure
    on $\BoundaryOf{\HorEuCell}$.

    A cell $\HorEuCell$ in $\TheEuBuild$ is \notion{$\TheMorse$-flat} if
    $\TheMorse$ restricts to a constant function on $\HorEuCell$.
    \begin{observation}\label{flat:parallel.and.normal}
      Let $\HorEuCell$ be an $\TheMorse$-flat cell. Then all flow lines
      issuing in $\HorEuCell$ are pairwise parallel and
      orthogonal to $\HorEuCell$.
    \end{observation}
    \begin{proof}
      For flow lines issuing from points carried by
      $\HorEuCell$, the claim is clear; and it
      follows for points on the boundary by continuity.
    \end{proof}
    Hence, we can talk about the \notion{gradient}
    $\GradOf[\HorEuCell]{\TheMorse}$ of a $\TheMorse$-flat cell as a
    point in the link $\LinkOf{\HorEuCell}$. Regarding the
    gradient as the north pole in the spherical building
    $\LinkOf{\HorEuCell}$, the link decomposes as in
    (\ref{eq:join.dec})
    \begin{equation}\label{eq:cell.link.split}
      \LinkOf{\HorEuCell}
      =
      \LinkOf[\hor]{\HorEuCell}
      \join
      \LinkOf[\ver]{\HorEuCell}
    \end{equation}
    into the horizontal and vertical parts of $\LinkOf{\HorEuCell}$
    relative to the north pole $\GradOf[\HorEuCell]{\TheMorse}$.
    We call the horizontal part
    $\LinkOf[\hor]{\HorEuCell}$
    the \notion{horizontal link}
    of the $\TheMorse$-flat cell $\HorEuCell$, and we call
    the vertical part $\LinkOf[\ver]{\HorEuCell}$
    its \notion{vertical link}. Beware that the
    vertical link can contain equatorial cells; and consequently
    not every $\TheMorse$-flat coface of a $\TheMorse$-flat
    cell $\HorEuCell$ defines
    a simplex in $\LinkOf[\hor]{\HorEuCell}$. It can also
    happen that a cell in the horizontal link is
    not $\TheMorse$-flat.

  \section{Depth}\label{sec:depth}
    Horizontal cells are the main obstacle for the analysis of
    the cocompact filtration of $\TheEuBuild$ by height. We will
    use the method of \cite{Bux.Wortman:2008} as extended to
    non-irreducible buildings in \cite{Witzel:2010}
    to cope with this difficulty. Here, we mostly
    follow \cite[Section~2]{Witzel:2010}.

    Let $\HorEuCell$ be an $\TheMorse$-flat cell in
    $\TheEuBuild$. By Observation~\ref{flat:parallel.and.normal},
    the flow lines starting in $\HorEuCell$ are pairwise parallel
    geodesic
    rays in $\TheEuBuild$ and therefore, they define a point
    $\TheEndOf{\HorEuCell}$
    in the spherical building at infinity. Let $\Busemann$ be
    a Busemann function centered at that point. Since the
    flow lines are orthogonal to $\HorEuCell$, the function
    $\Busemann$ is constant on $\HorEuCell$, i.e., the simplex
    $\HorEuCell$ is \notion{$\Busemann$-flat}.
    The notion of the
    horizontal and vertical link of $\HorEuCell$ defined above
    agree with the notions in \cite[Section~2]{Witzel:2010}, whence
    we can use some results therein directly.

    The Busemann function $\Busemann$ is not constant on any
    factor $\TheEuBuild[\ThePlace]$. In the Rank~Theorem,
    the group $\TheGroupScheme$ is assumed to be absolutely almost
    simple. Hence, the factors $\TheEuBuild[\ThePlace]$
    are all irreducible. It follows that $\Busemann$ is not constant
    on any irreducible factor $\TheEuBuild$, i.e., the Busemann
    function is in \notion{general position},
    see~Proposition~\ref{us:metric}.
    \begin{lemma}\label{equivalence}
      For any $\TheMorse$-flat cell $\HorEuCell$, there
      is a unique face $\HorEuCell[][\min]$ such that for any proper
      face $\HorEuFace<\HorEuCell$, the following equivalence holds
      \[
        \HorEuCell
        \text{\ defines a simplex in the horizontal link of\ }
        \HorEuFace
        \quad\text{\ if and only if\ }\quad
        \HorEuCell[][\min]
        \faceof
        \HorEuFace
        \Period
      \]
    \end{lemma}
    \begin{proof}
      Note that $\HorEuCell$ is $\Busemann$-flat for any Busemann
      function $\Busemann$ centered at $\TheEndOf{\HorEuCell}$.
      Then the statement follows from \cite[Lemma~2.7]{Witzel:2010}.
    \end{proof}
    In the same way, the following lemma is an immediate consequence
    of \cite[Observation~2.11]{Witzel:2010}.
    \begin{lemma}\label{wedging}
      Suppose $\HorEuCell[][\min]\faceof\HorEuFace\faceof\HorEuCell$,
      i.e., $\HorEuCell$ defines a simplex in the horizontal link
      of $\HorEuFace$.
      Then $\HorEuCell[][\min]=\HorEuFace[][\min]$.\qed
    \end{lemma}
    \begin{figure}[ht]
      \null\hfill
      \includegraphics[scale=1]{min_01.ps}\hfill
      \includegraphics[scale=1]{min_02.ps}\hfill\null
      \caption{the face $\HorEuCell[][\min]$\label{fig:min.illustrated}}
      {\small
      Both figures take place inside the Coxeter
      complex~$\tilde{\text{\textsf{B}}}_3$. In the picture on
      the left hand side,
      the black vertex is the face $\HorEuCell[][\min]$ of the
      horizontal solidly colored $\Two$-simplex $\HorEuCell$.
      The two edges of $\HorEuCell$ containing $\HorEuCell[][\min]$
      illustrate
      Lemma~\ref{wedging}. In the picture on the right, the
      horizontal simplex $\HorEuCell$ is the center edge. Here, we
      have $\HorEuCell=\HorEuCell[][\min]$.}
    \end{figure}

    For any two $\Busemann$-flat cells $\HorEuCell$
    and $\HorEuFace$, we define \notion{going up} as
    \[
      \HorEuFace \goesup \HorEuCell
      \quad
      :\Longleftrightarrow
      \quad
      \HorEuFace=\HorEuCell[][\min] \neq \HorEuCell
    \]
    and \notion{going down} as
    \[
      \HorEuCell \goesdown \HorEuFace
      \quad
      :\Longleftrightarrow
      \quad
      \HorEuCell[][\min] \not\faceof \HorEuFace
      \strictfaceof \HorEuCell.
    \]
    We define a \notion{$\Busemann$-move} as either going up or
    going down.
    \begin{figure}[ht]
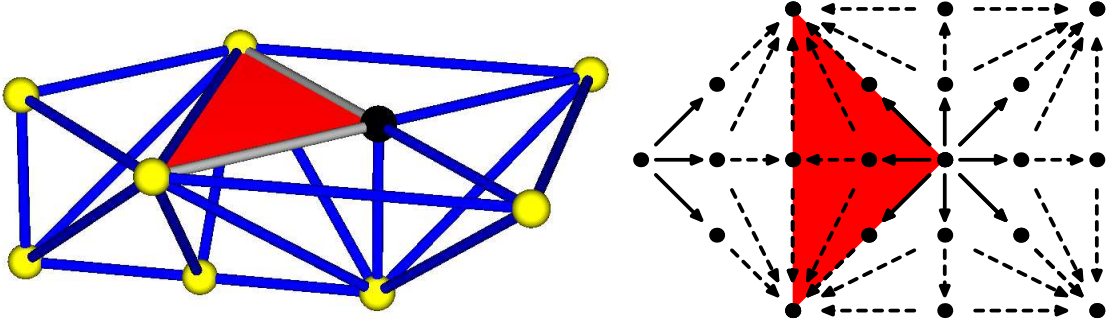

      \null\hfill
      \includegraphics{min_03.ps}\hfill
      \includegraphics{move.ps}\hfill\null
      \caption{sequences of moves\label{fig:move.sequence}}
      {\small
      This figure continues Figure~\ref{fig:min.illustrated}; it also
      takes place inside a Coxeter complex of
      type~$\tilde{\text{\textsf{B}}}_3$. It shows a possible
      patch of horizontal $\Two$-cells. Each dot in the
      picture on the right represents
      a cell; for orientation, one horizontal $\Two$-cells
      has been filled in. Arrows indicate moves: solid arrows
      represent going up, whereas dashed arrows represent going
      down. Note that there are no moves between triangles and
      short edges.}
    \end{figure}
    \begin{observation}\label{move}
      If there is a move from $\HorEuCell$ to $\HorEuCell'$,
      then either $\HorEuCell$ is a face of $\HorEuCell'$ or
      $\HorEuCell'$ is a face of $\HorEuCell$. In
      either case, we have
      $\TheEndOf{\HorEuCell}=\TheEndOf{\HorEuCell'}$.\qed
    \end{observation}
    The following is the statement
    of \cite[Proposition~2.9]{Witzel:2010}:
    \begin{prop}
      There is a uniform upper bound, depending only on the
      building $\TheEuBuild$, on the length of any sequence
      of $\Busemann$-moves.\qed
    \end{prop}
    We define the \notion{depth} $\DepthOf{\HorEuCell}$
    of an $\TheMorse$-flat cell $\HorEuCell$
    as the maximum length of a sequence of $\Busemann$-moves
    starting at $\HorEuCell$ for the corresponding Busemann
    function $\Busemann$ given by the flow lines of the gradient
    field on $\HorEuCell$.

  \section{The Morse function}\label{sec:morse}
    From now on, we assume that $\TheZonoDef$ is almost rich.
    Let $\TheEuCell$ be any cell of $\TheEuBuild$. By
    Proposition~\ref{height:minmax}, the max-set of $\TheMorse$
    on $\TheEuCell$ is a face $\RoofOf{\TheEuCell}$, which we call
    the \notion{roof} of $\TheEuCell$. The roof is $\TheMorse$-flat.
    We define the depth of $\TheEuCell$ of $\TheEuBuild$ as follows:
    \[
      \DepthOf{\TheEuCell}
      :=
      \begin{cases}
        \DepthOf{\TheEuCell} & \text{if\ } \TheEuCell\text{\ is flat} \\
        \DepthOf{\RoofOf{\TheEuCell}} - \OneHalf & \text{otherwise}
      \end{cases}
    \]

    We define the following Morse function on cells of $\TheEuBuild$:
    \begin{eqnarray*}
      \CombinedHeight \mapcolon \CellsOf{\TheEuBuild} & \longrightarrow &
      \TheReals\crossprod\TheReals\crossprod\TheReals \\
      \TheEuCell  & \mapsto &
      \FamOf{
        \MaxOf[\TheEuCell]{\TheMorse},
        \DepthOf{\TheEuCell},
        \DimOf{\TheEuCell}
      }
    \end{eqnarray*}
    \begin{observation}\label{morse:no.horizontal.edges}
      The dimension component assures that comparable but distinct cells
      (i.e., one is a strict face of the other or vice versa) are not
      assigned the same triple.\qed
    \end{observation}
    The cells of $\TheEuBuild$ are in one-to-one correspondence
    to the vertices of the barycentric subdivision
    $\TheEuSubdiv$ of
    $\TheEuBuild$. Ordering
    $\TheReals\crossprod\TheReals\crossprod\TheReals$
    lexicographically, we regard $\CombinedHeight$ as a Morse
    function on $\TheEuSubdiv$.

    For each cell $\TheEuCell$ let
    $\CenterOf{\TheEuCell}$ denote its barycenter,
    i.e., the vertex of $\TheEuSubdiv$ that corresponds
    to $\TheEuCell$. The link of the vertex
    $\CenterOf{\TheEuCell}$ in
    $\TheEuSubdiv$ decomposes as a join of the
    boundary $\BoundaryOf{\TheEuCell}$ and the
    link $\LinkOf{\TheEuCell}$. This corresponds to
    the decomposition~(\ref{eq:point.link.split}).

    The \notion{descending link} of a vertex
    $\CenterOf{\TheEuCell}\in\TheEuSubdiv$ consists
    of those simplices in the link all of whose
    vertices lie strictly below $\CenterOf{\TheEuCell}$
    with respect to the Morse function $\CombinedHeight$.
    Hence, the descending link also decomposes as
    a join:
    \begin{equation}\label{eq:desc.link.coarse.join}
      \LinkOf[][\downarrow]{\CenterOf{\TheEuCell}}
      =
      \BoundaryOf[][\downarrow]{\TheEuCell}
      \join
      \LinkOf[][\downarrow]{\TheEuCell}
    \end{equation}
    Here, $\BoundaryOf[][\downarrow]{\TheEuCell}$ can be
    regarded as the poset of strict faces of $\TheEuCell$
    with smaller $\CombinedHeight$-value and
    $\LinkOf[][\downarrow]{\TheEuCell}$ can likewise be
    viewed as the poset of strict cofaces with
    smaller $\CombinedHeight$-value.

  \section{Descending links}\label{sec:links}
    We call $\TheEuCell$ \notion{insignificant} if
    \(
      \TheEuCell
      \neq
      \RoofOf[\min]{\TheEuCell}
    \)
    and \notion{significant} otherwise.
    We shall deal with the insignificant cells first. Here, the
    descending link is always contractible. In fact, in
    the decomposition
    \(
      \LinkOf[][\downarrow]{\CenterOf{\TheEuCell}}
      =
      \BoundaryOf[][\downarrow]{\TheEuCell}
      \join
      \LinkOf[][\downarrow]{\TheEuCell}
    \),
    already the boundary part
    \(
      \BoundaryOf[][\downarrow]{\TheEuCell}
    \)
    is contractible.
    \begin{prop}\label{links:insignificant.prop}
      If
      \(
        \TheEuCell
        \neq
        \RoofOf[\min]{\TheEuCell}
      \)
      then
      \(
        \BoundaryOf[][\downarrow]{\TheEuCell}
      \)
      is contractible. More precisely,
      the complex
      \(
        \BoundaryOf[][\downarrow]{\TheEuCell}
      \)
      deformation retracts onto the subcomplex
      \(
        \BoundaryOf{\TheEuCell}
        \setminus
        \StarOf{\RoofOf[\min]{\TheEuCell}}
      \).
    \end{prop}
    \begin{figure}[h]
      \null\hfill
      \setlength{\unitlength}{1cm}
      \begin{picture}(10,4)
        \put(0,0){\includegraphics{deformation.ps}}
      \end{picture}
      \hfill\null
      \caption{the deformation of
      \(
        \BoundaryOf[][\downarrow]{\TheEuCell}
      \)
      onto
      \(
        \BoundaryOf{\TheEuCell}
        \setminus
        \StarOf{\RoofOf[\min]{\TheEuCell}}
      \) from Proposition~\ref{links:insignificant.prop}%
      \label{links:deformation.fig}}
      {\small
      On the left hand, the solid vertex at the top is
      \(
        \RoofOf[\min]{\TheEuCell}
      \).
      The shaded triangle is
      \(
        \BoundaryOf{\TheEuCell}
        \setminus
        \StarOf{\RoofOf[\min]{\TheEuCell}}
      \)
      and guaranteed to be descending. Note that the
      barycentric subdivision is drawn but
      \(
        \StarOf{\RoofOf[\min]{\TheEuCell}}
      \)
      denotes the open star of
      \(
        \RoofOf[\min]{\TheEuCell}
      \)
      with respect to the original cell structure.
      The hollow vertices are
      \emph{unknown}; they could be descending or
      ascending. On the right hand, one possibility for
      \(
        \BoundaryOf[][\downarrow]{\TheEuCell}
      \)
      is shown.}
    \end{figure}
    \begin{proof}
      First, we note that
      \(
        \RoofOf[\min]{\TheEuCell}
      \)
      cannot correspond to a vertex in the descending link.
      The height does not decide as
      \[
        \MaxOf[\TheEuCell]{\TheMorse}
        =
        \MaxOf[\RoofOf{\TheEuCell}]{\TheMorse}
        =
        \MaxOf[{\RoofOf[\min]{\TheEuCell}}]{\TheMorse}
        \Period
      \]
      As for the depth, we have
        \
      \[
        \DepthOf{\TheEuCell}
        \leq
        \DepthOf{\RoofOf{\TheEuCell}}
        \leq
        \DepthOf{\RoofOf[\min]{\TheEuCell}}
        \Period
      \]
      If $\TheEuCell\neq\RoofOf{\TheEuCell}$, then
      the first inequality is strict. If
      $\TheEuCell=\RoofOf{\TheEuCell}$, the hypothesis
      that $\TheEuCell$ is insignificant implies
      $\RoofOf{\TheEuCell}\neq\RoofOf[\min]{\TheEuCell}$
      whence there is a move
      \(
        \RoofOf[\min]{\TheEuCell}
        \goesup
        \RoofOf{\TheEuCell}
      \)
      and
      \(
        \DepthOf{\RoofOf{\TheEuCell}}
        <
        \DepthOf{\RoofOf[\min]{\TheEuCell}}
      \).
      In either case the strict inequality
       \(
        \DepthOf{\TheEuCell}
        <
        \DepthOf{\RoofOf[\min]{\TheEuCell}}
      \)
      follows, whence
      \(
        \CombinedHeightOf{\TheEuCell}
        <
        \CombinedHeightOf{\RoofOf[\min]{\TheEuCell}}
      \),
      i.e.,
      \(
        \RoofOf[\min]{\TheEuCell}
      \)
      is not descending.

      We now turn to the ``opposite part'', i.e., we
      identify a subcomplex of
      $\BoundaryOf{\TheEuCell}$ that is descending.
      Let $\AltEuCell$ be a cell with
      \(
        \RoofOf[\min]{\TheEuCell}
        \not\faceof
        \AltEuCell
        \strictfaceof
        \TheEuCell
      \).
      As
      \(
        \AltEuCell
        \strictfaceof
        \TheEuCell
      \)
      we have
      \(
        \MaxOf[\AltEuCell]{\TheMorse}
        \leq
        \MaxOf[\TheEuCell]{\TheMorse}
      \).
      Also,
      \(
        \RoofOf[\min]{\TheEuCell}
        \not\faceof
        \AltEuCell
      \)
      implies
      \(
        \RoofOf[\min]{\TheEuCell}
        \not\faceof
        \RoofOf{\AltEuCell}
      \).
      Hence, there is a move
      \(
        \RoofOf{\TheEuCell}
        \goesdown
        \RoofOf{\AltEuCell}
      \)
      whence
      \(
        \DepthOf{\RoofOf{\AltEuCell}}
        <
        \DepthOf{\RoofOf{\TheEuCell}}
      \).
      Since depths of flat cells are integer valued,
      \(
        \DepthOf{\AltEuCell}
        \leq
        \DepthOf{\RoofOf{\AltEuCell}}
        <
        \DepthOf{\RoofOf{\TheEuCell}}-\OneHalf
        \leq
        \DepthOf{\TheEuCell}
      \).

      So, let $\TheDescBoundary$ be the subcomplex
      of $\BoundaryOf{\TheEuCell}$ spanned by vertices
      $\CenterOf{\AltEuCell}$ with
      \(
        \RoofOf[\min]{\TheEuCell}
        \not\faceof
        \AltEuCell
        \strictfaceof
        \TheEuCell
      \).
      We have seen that $\TheDescBoundary$ is descending.
      Let $\TheAscVertex$ be the barycenter of
      $\RoofOf[\min]{\TheEuCell}$. We have seen that
      $\TheAscVertex$ is ascending, not descending.
      Radial projection inside $\TheEuCell$ away from
      $\TheAscVertex$ defines
      a deformation retraction of
      $\BoundaryOf[][\downarrow]{\TheEuCell}$
      onto
      $\TheDescBoundary$. Since $\TheDescBoundary$ is
      a sphere with an open star of a cell removed,
      $\BoundaryOf[][\downarrow]{\TheEuCell}$
      is contractible.
    \end{proof}
    \begin{cor}\label{links:insignificant}
      If $\TheEuCell$ is insignificant then the descending
      link $\DescLinkOf{\CenterOf{\TheEuCell}}$ of
      its barycenter is contractible.\qed
    \end{cor}

    The remainder of this section is devoted to
    the analysis of descending links
    $\DescLinkOf{\CenterOf{\TheEuCell}}$ when
    $\TheEuCell$ is a significant cell.
    \begin{observation}\label{links:significant.boundary.part}
      If $\TheEuCell$ is significant, i.e.,
      $\TheEuCell=\RoofOf[\min]{\TheEuCell}$, then
      $\TheEuCell$ is flat whence
      $\MaxOf[\AltEuCell]{\TheMorse}=
      \MaxOf[\TheEuCell]{\TheMorse}$
      for any face $\AltEuCell\faceof\TheEuCell$.
      In particular, the depth and the dimension
      determine which part of $\BoundaryOf{\TheEuCell}$
      is descending.

      It follows that
      $\BoundaryOf{\TheEuCell}$ is completely
      descending: for any proper face
      $\AltEuCell\strictfaceof\TheEuCell$, there
      is a move $\TheEuCell\goesdown\AltEuCell$ whence
      $\DepthOf{\AltEuCell}<\DepthOf{\TheEuCell}$.
      Thus,
      $\BoundaryOf[][\downarrow]{\TheEuCell}$
      is a sphere of dimension
      $\DimOf{\TheEuCell}-\One$.\qed
    \end{observation}
    \begin{observation}\label{links:key}
      If $\TheEuCoface\strictcofaceof\TheEuCell=
      \RoofOf[\min]{\TheEuCell}$ is a flat coface of
      a significant cell, then either
      $\TheEuCell\goesup\TheEuCoface$ or
      $\TheEuCoface\goesdown\TheEuCell$:
      If $\TheEuCoface[][\min]\faceof\TheEuCell$,
      then $\TheEuCoface[][\min]=\TheEuCell[][\min]
      =\TheEuCell$ by Lemma~\ref{wedging}. In this
      case, $\TheEuCell\goesup\TheEuCoface$. If
      $\TheEuCoface[][\min]\not\faceof\TheEuCell$, then
      $\TheEuCoface\goesdown\TheEuCell$.\qed
    \end{observation}
    \begin{prop}\label{links:significant.full}
      Assume that $\TheEuCell$ is significant. Fix cofaces
      $\TheEuCoface$ and $\AltEuCoface$ with
      \(
        \TheEuCell\strictfaceof\TheEuCoface
        \faceof\AltEuCoface
      \).
      If $\CombinedHeightOf{\AltEuCoface}<\CombinedHeightOf{\TheEuCell}$,
      then
      $\CombinedHeightOf{\TheEuCoface}<\CombinedHeightOf{\TheEuCell}$.
    \end{prop}
    \begin{proof}
      First note that
      \(
        \MaxOf[\TheEuCell]{\TheMorse}
        \leq
        \MaxOf[\TheEuCoface]{\TheMorse}
        \leq
        \MaxOf[\AltEuCoface]{\TheMorse}
      \)
      as
      \(
        \TheEuCell\strictfaceof\TheEuCoface
        \faceof\AltEuCoface
      \).
      By hypothesis,
      \(
        \MaxOf[\AltEuCoface]{\TheMorse}
        \leq
        \MaxOf[\TheEuCell]{\TheMorse}
      \).
      Thus, we have equality throughout.

      As $\DimOf{\AltEuCoface}>\DimOf{\TheEuCell}$, the
      hypothesis
      \(
        \CombinedHeightOf{\AltEuCoface}<\CombinedHeightOf{\TheEuCell}
      \)
      implies
      \(
        \DepthOf{\AltEuCoface}<\CombinedHeightOf{\TheEuCell}
      \).
      Passing to roofs, we have the inclusions
      \(
        \TheEuCell
        \faceof
        \RoofOf{\TheEuCoface}
        \faceof
        \RoofOf{\AltEuCoface}
      \)
      of flat cells. If
      \(
        \TheEuCell
        =
        \RoofOf{\TheEuCoface}
      \)
      then
      \(
        \TheEuCoface
        \neq
        \RoofOf{\TheEuCoface}
      \).
      Hence,
      \(
        \DepthOf{\TheEuCoface}
        <
        \DepthOf{\RoofOf{\TheEuCoface}}
        =
        \DepthOf{\TheEuCell}
      \)
      and $\TheEuCoface$ is descending.

      If, on the other hand,
      \(
        \TheEuCell
        \neq
        \RoofOf{\TheEuCoface}
      \)
      then $\RoofOf{\AltEuCoface}$ is a proper flat coface of the
      significant cell $\TheEuCell$. By Observation~\ref{links:key},
      there is a move $\TheEuCell\goesup\RoofOf{\AltEuCoface}$ or a
      move $\RoofOf{\AltEuCoface}\goesdown\TheEuCell$. In the latter
      case,
      \(
        \DepthOf{\AltEuCoface}
        \geq
        \DepthOf{\RoofOf{\AltEuCoface}}
        >
        \DepthOf{\TheEuCell}
      \)
      contradicting the hypothesis that
      $\CombinedHeightOf{\AltEuCoface}<\CombinedHeightOf{\TheEuCell}$.
      Therefore, there is a move
      \(
        \TheEuCell\goesup\RoofOf{\AltEuCoface}
      \),
      that is, $\TheEuCell=\RoofOf[\min]{\AltEuCoface}$.
      Then Lemma~\ref{wedging} implies
      \(
        \TheEuCell=\RoofOf[\min]{\TheEuCoface}
      \)
      whence
      \(
        \TheEuCell\goesup\RoofOf[\min]{\TheEuCoface}
      \)
      and
      \(
        \DepthOf{\TheEuCoface}
        \leq
        \DepthOf{\RoofOf{\TheEuCoface}}
        <\DepthOf{\TheEuCell}
      \).
    \end{proof}
    Proposition~\ref{links:significant.full} justifies a
    notational vagueness of which we are guilty. In
    Section~\ref{sec:flat}, particularly in the
    decomposition~\ref{eq:point.link.split}, we used
    $\LinkOf{\TheEuCell}$ to denote a spherical building.
    In Section~\ref{sec:morse}, we switched to its
    barycentric subdivision (the geometric realization of
    the poset of strict cofaces of $\TheEuCell$).
    Since barycenters of insignificant cells have contractible
    descending links just by their boundary part
    $\BoundaryOf[][\downarrow]{\TheEuCell}$, the
    precise structure of $\DescLinkOf{\TheEuCell}$ did
    not matter in this case. If $\TheEuCell$ is significant,
    $\DescLinkOf{\TheEuCell}$ does matter. Although it
    is defined as a subcomplex of the barycentric subdivision,
    Proposition~\ref{links:significant.full} implies
    that we can regard the descending link as a subcomplex
    of the spherical building $\LinkOf{\TheEuCell}$.
    Hence, we put:
    \begin{eqnarray*}
      \DescLinkOf[\hor]{\TheEuCell}
      &:=&
      \DescLinkOf{\TheEuCell}
      \intersect
      \LinkOf[\hor]{\TheEuCell}
      \\
      \DescLinkOf[\ver]{\TheEuCell}
      &:=&
      \DescLinkOf{\TheEuCell}
      \intersect
      \LinkOf[\ver]{\TheEuCell}
    \end{eqnarray*}
    \begin{prop}\label{links:vertical}
      If $\TheEuCell$ is significant,
      $\DescLinkOf[\ver]{\TheEuCell}$ is an
      open hemisphere complex in $\LinkOf{\TheEuCell}$
      with respect to the north pole
      $\GradOf[\TheEuCell]{\TheMorse}$.
    \end{prop}
    \begin{proof}
      Let $\TheHemisphereCx$ be the open hemisphere complex
      of $\LinkOf{\TheEuCell}$
      with respect to $\GradOf[\TheEuCell]{\TheMorse}$.
      The Gradient~Criterion~\ref{height:criterion} implies
      \(
        \TheHemisphereCx
        \subseteq
        \DescLinkOf[\ver]{\TheEuCell}
      \).
      To show equality, it suffices to argue that no
      flat coface $\TheEuCoface\strictcofaceof\TheEuCell$
      in $\LinkOf[\ver]{\TheEuCell}$ is descending, i.e.,
      $\CombinedHeightOf{\TheEuCoface}>\CombinedHeightOf{\TheEuCell}$.
      However, as $\TheEuCoface$ belongs to
      $\LinkOf[\ver]{\TheEuCell}$, it does not belong to
      $\LinkOf[\hor]{\TheEuCell}$ whence, by
      Lemma~\ref{equivalence}, there is a move
      $\TheEuCoface\goesdown\TheEuCell$. Thus,
      $\DepthOf{\TheEuCoface}>\DepthOf{\TheEuCell}$. As
      $\TheEuCoface$ is flat,
      $\CombinedHeightOf{\TheEuCoface}>\CombinedHeightOf{\TheEuCell}$ follows.
    \end{proof}
    \begin{prop}
      The decomposition~\ref{eq:join.dec} at the end of
      Section~\ref{sec:bernd} induces the decomposition
      \[
        \DescLinkOf{\TheEuCell}
        =
        \DescLinkOf[\hor]{\TheEuCell}
        \join
        \DescLinkOf[\ver]{\TheEuCell}
      \]
      provided $\TheEuCell$ is significant.
    \end{prop}
    \begin{proof}
      As each vertex of $\DescLinkOf{\TheEuCell}$ lies in
      \(
        \DescLinkOf[\hor]{\TheEuCell}
      \)
      or
      \(
        \DescLinkOf[\ver]{\TheEuCell}
      \),
      it follows that
      \(
        \DescLinkOf{\TheEuCell}
        \subseteq
        \DescLinkOf[\hor]{\TheEuCell}
        \join
        \DescLinkOf[\ver]{\TheEuCell}
      \).

      To see the converse, let $\TheEuCoface[\hh]$ and
      $\TheEuCoface[\vv]$ denote strict cofaces of
      $\TheEuCell$ where $\TheEuCoface[\hh]$ determines
      a simplex in $\DescLinkOf[\hor]{\TheEuCell}$ and
      $\TheEuCoface[\vv]$ determines a simplex in
      $\DescLinkOf[\ver]{\TheEuCell}$. We need to show
      that the join
      \(
        \TheEuCoface:=
        \TheEuCoface[\hh]\join\TheEuCoface[\vv]
      \)
      lies in $\DescLinkOf{\TheEuCell}$, i.e.,
      $\CombinedHeightOf{\TheEuCoface}<\CombinedHeightOf{\TheEuCell}$. The cell
      $\TheEuCoface$ is the smallest coface of $\TheEuCell$
      containing $\TheEuCoface[\hh]$ and $\TheEuCoface[\vv]$.

      Since $\TheEuCoface[\hh]$ is descending, it is flat.
      As $\DescLinkOf[\ver]{\TheEuCell}$ is an open hemisphere
      complex, all directions from $\TheEuCell$ that have
      a non-vanishing component into a direction of
      $\TheEuCoface[\vv]$ are descending with respect to
      the height $\TheMorse$. Hence
      $\TheEuCoface[\hh]=\RoofOf{\TheEuCoface}$ and
      \(
        \MaxOf[\TheEuCoface]{\TheMorse}
        =
        \MaxOf[{\TheEuCoface[\hh]}]{\TheMorse}
        =
        \MaxOf[\TheEuCell]{\TheMorse}
      \).
      Since $\TheEuCoface[\hh]$ is descending,
      \(
        \DepthOf{\TheEuCoface}
        \leq
        \DepthOf{\TheEuCoface[\hh]}
        <
        \DepthOf{\TheEuCell}
      \).
      Thus,
      \(
        \CombinedHeightOf{\TheEuCoface}<\CombinedHeightOf{\TheEuCell}
      \).
    \end{proof}
    \begin{cor}\label{links:fine.decomposition}
      For significant $\TheEuCell$, the descending link
      $\DescLinkOf{\CenterOf{\TheEuCell}}$
      decomposes as
      \[
        \DescLinkOf{\CenterOf{\TheEuCell}}
        =
        \BoundaryOf{\TheEuCell}
        \join
        \DescLinkOf[\ver]{\TheEuCell}
        \join
        \DescLinkOf[\hor]{\TheEuCell}
        \Period
        \qed
      \]
    \end{cor}
    In Proposition~\ref{links:vertical}, we have determined that
    $\DescLinkOf[\ver]{\TheEuCell}$ is an open hemisphere complex.
    It remains to analyze $\DescLinkOf[\hor]{\TheEuCell}$.

    \begin{lemma}\label{links:horizontal.descending.iff.flat}
      Let $\TheEuCoface\in\LinkOf[\hor]{\TheEuCell}$ for
      a significant cell $\TheEuCell$, i.e., $\TheEuCoface$ is
      a proper coface of $\TheEuCell$ with
      \(
        \TheEuCoface[][\min]\faceof
        \TheEuCell
        =
        \TheEuCell[][\min]
        \strictfaceof
        \TheEuCoface
      \).
      Then, the following are
      equivalent:
      \begin{enumerate}
        \item
          The cell $\TheEuCoface$ is descending, i.e.,
          $\TheEuCoface\in\DescLinkOf[\hor]{\TheEuCell}$.
        \item
          The cell $\TheEuCoface$ is $\TheMorse$-flat.
        \item
          We have
          \(
            \MaxOf[\TheEuCoface]{\TheMorse}
            =
            \MaxOf[\TheEuCell]{\TheMorse}
          \).
      \end{enumerate}
    \end{lemma}
    \begin{figure}[t]
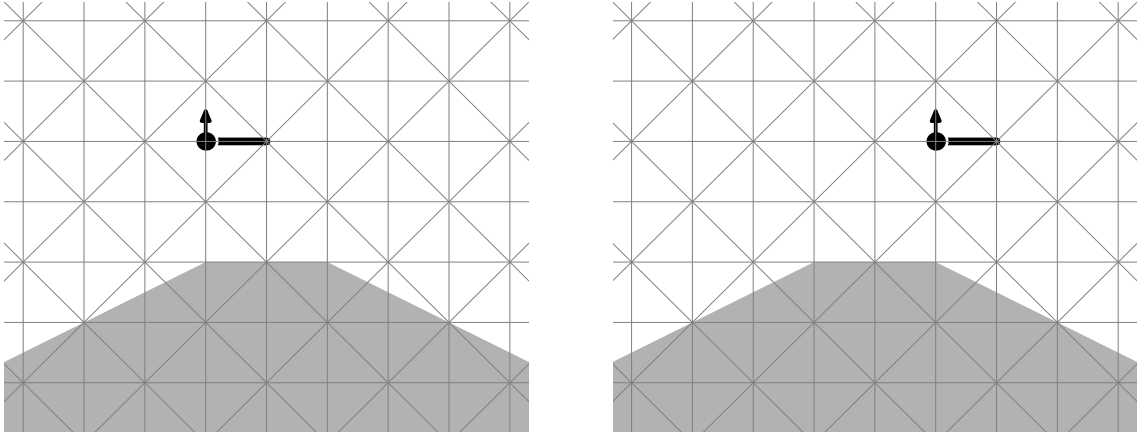

      \null\hfill
      \setlength{\unitlength}{1cm}
      \begin{picture}(16,6)
        \put(0,0){\includegraphics{horizontal_flat.ps}}
        \put(8,0){\includegraphics{horizontal_nonflat.ps}}
      \end{picture}
      \hfill\null
      \caption{a flat vs.\ a non-flat horizontal coface%
      \label{links:flat.vs.nonflat.fig}}
      {\small
      The primary height is distance from the shaded
      area, the significant cell is the fat vertex,
      the little arrow indicates the gradient of the
      height function, and
      the horizontal coface is the marked edge issuing
      to the right. To see that the edge is in the
      horizontal link of its left vertex, recall that
      the link decomposes as a join into the vertical
      and horizontal parts. The horizontal part is
      the maximal join factor that is perpendicular to
      the gradient.}
    \end{figure}
    \begin{proof}
      First assume that $\TheEuCoface$ is flat.
      Then
          \(
            \MaxOf[\TheEuCoface]{\TheMorse}
            =
            \MaxOf[\TheEuCell]{\TheMorse}
          \).
      Also, by Lemma~\ref{wedging},
      \(
        \TheEuCoface[][\min]
        =
        \TheEuCell[][\min]
        =
        \TheEuCell
      \).
      Hence, there is a move
      \(
        \TheEuCell\goesup\TheEuCoface
      \)
      whence
      \(
        \DepthOf{\TheEuCoface}<\DepthOf{\TheEuCell}.
      \)
      Hence $\CombinedHeightOf{\TheEuCoface}<\CombinedHeightOf{\TheEuCell}$.

      Now assume that $\TheEuCoface$ is not flat.
      The Gradient Criterion~\ref{height:criterion} implies
      that
      \(
        \MaxOf[\TheEuCoface]{\TheMorse}
        >
        \MaxOf[\TheEuCell]{\TheMorse}
      \).
      In particular,
      \(
        \CombinedHeightOf{\TheEuCoface}
        >
        \CombinedHeightOf{\TheEuCell}
      \).
    \end{proof}

    For the last part of the analysis, we need to strengthen
    the hypothesis on $\TheZonoDef$ one last time. We call
    $\TheZonoDef$ \notion{rich} if it contains the differences
    $\TheSpVertex-\AltSpVertex$ of any to vertices
    $\TheSpVertex,\AltSpVertex\in\TheModelApp$
    whose closed stars intersect.

    Also, at last, we have to enlarge the diameter
    $\TheDiam$ of uniformity for the reduction datum.
    Using Theorem~\ref{behr:main}, we assume that any
    closed star of any cell can be uniformly reduced.
    This affects the constants $\TheLowerBound$ and
    $\TheUpperBound$. Using a rich $\TheZonoDef$, the
    construction of Section~\ref{sec:height} then will
    yield an appropriate $\TheHugeBound$.

    Let $\TheEuCell$ be a significant cell, let
    $\TheRatChamber$ be a chamber in $\TheRatBuild$
    uniformly reducing the closed star of $\TheEuCell$,
    and let $\TheEuApp$ be an apartment containing
    $\TheEuCell$
    with $\TheRatChamber\subset\BoundaryOf{\TheEuApp}$.
    Put:
    \[
      \UpSetOf[\TheEuApp]{\TheEuCell}
      :=
      \SetOf[{
        \TheAppVertex\in\TheEuApp
      }]{
        \TheAppVertex\text{\ is a vertex,\ }
        \TheAppVertex\commoncell\TheEuCell
        \text{\ defines a cell in\ }
        \LinkOf{\TheEuCell},\,\,
        \TheMorseOf{\TheAppVertex}
        >
        \TheMorseOf{\TheEuCell}
      }
    \]
    Let
    \(
      \TheHull[\TheEuApp]
    \)
    denote the convex hull of
    $\UpSetOf[\TheEuApp]{\TheEuCell}$ in the euclidean space
    $\TheEuApp$.
    \begin{observation}\label{links:upset.hull}
      Assume that $\TheZonoDef$ is rich. Then
      $\TheHull[\TheEuApp]$
      is a convex polytope satisfying the hypotheses
      of Proposition~\ref{zonotopes:minset.maxset}.
      Hence, $\TheMorse$ assumes its minimum on
      $\TheHull[\TheEuApp]$ in a vertex, which is still
      higher than the flat cell $\TheEuCell$.\qed
    \end{observation}
    \begin{cor}\label{links:upset.disjoint}
      Provided that $\TheZonoDef$ is rich,
      $\TheHull[\TheEuApp]$
      is disjoint from the affine subspace
      of $\TheEuApp$ spanned by $\TheEuCell$.\qed
    \end{cor}
    The convex set $\TheHull[\TheEuApp]$ induces a closed
    subset $\TheHullProj[\TheEuApp]$ in
    $\LinkOf[\TheEuApp]{\TheEuCell}\subset
    \LinkOf{\TheEuCell}$ by
    projection onto an orthogonal complement of the span
    of $\TheEuCell$.
    \begin{cor}\label{links:horizontal.upset}
      Also under the hypothesis that $\TheZonoDef$ is rich,
      the subset $\TheHullProj[\TheEuApp]\subset
      \LinkOf[\TheEuApp]{\TheEuCell}$
      is closed, convex, and has diameter strictly less
      than $\pi$.\qed
    \end{cor}
    We can extract a little more information:
    \begin{observation}\label{links:upset.not.descending}
      By Lemma~\ref{links:horizontal.descending.iff.flat},
      a horizontal coface $\TheEuCoface$
      of $\TheEuCell$ is descending if
      and only if it is flat. Hence, the value of $\TheMorse$
      on $\TheEuCoface$ cannot exceed the value on $\TheEuCell$.
      Therefore,
      if $\TheZonoDef$ is rich,
      the descending horizontal link
      $\DescLinkOf[\hor]{\TheEuCell}$ and
      $\TheHullProj[\TheEuApp]$ are disjoint
      by Observation~\ref{links:upset.hull}.\qed
    \end{observation}
    Let $\AltEuApp$ be another apartment in $\TheEuBuild$
    containing $\TheEuCell$ and satisfying
    $\TheRatChamber\subset\BoundaryOf{\AltEuApp}$.
    \begin{observation}\label{links:isometry}
      Any Coxeter isomorphism
      \(
        \TheIsometry\mapcolon
        \TheEuApp \rightarrow \AltEuApp
      \)
      that is the identity on the intersection
      \(
        \TheEuApp\intersect\AltEuApp
      \)
      makes the diagram
      \[
        \xymatrix{
          {\TheEuApp}
          \ar[rr]^{\TheIsometry}
          \ar[rd]_{\TheMorse[\TheEuApp,\TheRatChamber]}
          &
          &
          {\AltEuApp}
          \ar[ld]^{\TheMorse[\AltEuApp,\TheRatChamber]}
          \\
          &
          {\TheReals}
          &
        }
      \]
      commute.
    \end{observation}
    \begin{proof}
      The height $\TheMorse$ is defined in terms of
      (a)~Busemann functions $\TheBusemann[\TheRatVertex]$ for
      points $\TheRatVertex\in\TheRatChamber$ and (b)~zonotopes
      of the form $\TheSpPoint+\TheZonotopeOf{\TheZonoDef}$.
      The Busemann functions are clearly preserved under the
      Coxeter isomorphism $\TheIsometry$ since the intersection
      $\TheEuApp\intersect\AltEuApp$ contains a sector bounding
      $\TheRatChamber$. Since $\TheIsometry$ is a Coxeter
      isomorphism, it also preserves $\TheZonoDef$ which is
      invariant under the full spherical Weyl group.
    \end{proof}
    Now, we fix a chamber $\TheBoundaryChamber$ in the visual
    boundary $\BoundaryOf{\TheEuBuild}$ that contains
    $\TheRatChamber$. If $\TheEuApp$ and $\AltEuApp$ are
    two apartments both containing the convex cone
    $\ConvexHullOf{\TheEuCell,\TheBoundaryChamber}$ then the
    retraction
    \(
      \TheRetraction \mapcolon \TheEuBuild
      \rightarrow
      \TheEuApp
    \)
    of the building $\TheEuBuild$ onto $\TheEuApp$ centered
    at the chamber $\TheBoundaryChamber$ restricts to a Coxeter
    isomorphism
    \(
      \TheIsomorphism \mapcolon \AltEuApp
      \rightarrow
      \TheEuApp
    \)
    to which Observation~\ref{links:isometry} applies.
    Now, we use the hypothesis that $\TheRatChamber$ uniformly
    reduces the closed star of $\TheEuCell$. Hence:
    \begin{eqnarray*}
      \TheMorse|_{\StarOf{\TheEuCell}\intersect\TheEuApp}
      &=&
      \TheMorse[\TheEuApp,\TheRatChamber]|_{\StarOf{\TheEuCell}\intersect\TheEuApp}
      \\
      \TheMorse|_{\StarOf{\TheEuCell}\intersect\AltEuApp}
      &=&
      \TheMorse[\AltEuApp,\TheRatChamber]|_{\StarOf{\TheEuCell}\intersect\AltEuApp}
    \end{eqnarray*}
    In particular, $\TheIsomorphism$ identifies
    $\TheHull[\AltEuApp]$ with $\TheHull[\TheEuApp]$.
    \begin{observation}\label{links:apply.bernd}
      We assume that $\TheZonoDef$ is rich so that we can use the
      previous results.
      Let $\TheLinkChamber$ be the projection of $\TheBoundaryChamber$
      in the spherical building $\LinkOf[\hor]{\TheEuCell}$. It is a
      chamber. Let $\TheLinkApp$ be the apartment
      $\TheEuApp\intersect\LinkOf[\hor]{\TheEuCell}$. Then
      $\TheLinkChamber\subset\TheLinkApp$. Let
      \(
        \TheLinkRetraction\mapcolon\LinkOf[\hor]{\TheEuCell}
        \rightarrow\TheLinkApp
      \)
      be the retraction onto $\TheLinkApp$ centered at
      $\TheLinkChamber$.
      Put
      \(
        \TheHullProj:=\TheLinkRetractionOf[][-1]{
          \TheHullProj[\TheEuApp]
        }
      \).
      Then, any apartment in $\LinkOf[\hor]{\TheEuCell}$ that contains
      $\TheLinkChamber$ is
      of the form
      \(
        \AltLinkApp :=
        \AltEuApp
        \intersect
        \LinkOf[\hor]{\TheEuCell}
      \)
      where $\AltEuApp$ is an apartment in $\TheEuBuild$
      containing the convex cone
      $\ConvexHullOf{\TheEuCell,\TheBoundaryChamber}$; moreover
      $\TheHullProj\intersect\AltEuApp
      =\TheHullProj[\AltEuApp]$ is a closed convex subset of
      $\AltLinkApp$ of diameter less than $\pi$ by
      Corollary~\ref{links:horizontal.upset}.

      Hence~Corollary~\ref{bernd:complement:closed} applies; and the
      the maximal subcomplex of the complement
      \(
        \LinkOf[\hor]{\TheEuCell}
        \setminus
        \TheHullProj
      \)
      is
      \(
        (\DimOf{\LinkOf[\hor]{\TheEuCell}} - \One)
      \)-connected and of dimension
      $\DimOf{\LinkOf[\hor]{\TheEuCell}}$.\qed
    \end{observation}
    \begin{cor}\label{links:horizontal}
      Assume that $\TheZonoDef$ is rich.
      The horizontal descending link
      \(
        \DescLinkOf[\hor]{\TheEuCell}
      \)
      of a significant cell $\TheEuCell$ is contractible or
      spherical of dimension
      \(
        \DimOf{\LinkOf[\hor]{\TheEuCell}}
      \).
    \end{cor}
    \begin{proof}
      By the preceding Observation~\ref{links:apply.bernd}, we
      have to argue that
      \(
        \DescLinkOf[\hor]{\TheEuCell}
      \)
      is the maximal subcomplex of
      \(
        \LinkOf[\hor]{\TheEuCell}
        \setminus
        \TheHullProj
      \).
      Observation~\ref{links:upset.not.descending} implies
      the inclusion
      \(
        \DescLinkOf[\hor]{\TheEuCell}
        \subseteq
        \LinkOf[\hor]{\TheEuCell}
        \setminus
        \TheHullProj
      \).
      On the other hand, any vertex of
      $\LinkOf[\hor]{\TheEuCell}\setminus\DescLinkOf[\hor]{\TheEuCell}$
      lies within $\TheHullProj$ by definition of the sets
      $\UpSetOf[\TheEuApp]{\TheEuCell}$.
    \end{proof}

    We can summarize the analysis of descending links:
    \begin{prop}\label{links:connectivity}
      Assume that $\TheZonoDef$ is rich. Then the descending link
      $\DescLinkOf{\CenterOf{\TheEuCell}}$ of any barycenter
      is contractible or spherical of dimension
      $\DimOf{\TheEuBuild}-\One$.
    \end{prop}
    \begin{proof}
      If $\TheEuCell$ is insignificant, then
      $\DescLinkOf{\TheEuCell}$ is contractible by
      Corollary~\ref{links:insignificant}.

      If $\TheEuCell$ is significant, then the descending link
      decomposes as
      \[
        \DescLinkOf{\CenterOf{\TheEuCell}}
        =
        \BoundaryOf{\TheEuCell}
        \join
        \DescLinkOf[\ver]{\TheEuCell}
        \join
        \DescLinkOf[\hor]{\TheEuCell}
      \]
      by Corollary~\ref{links:fine.decomposition}.
      The part $\BoundaryOf{\TheEuCell}$ is a sphere of dimension
      $\DimOf{\TheEuCell}-\One$ (or empty if $\TheEuCell$ is a vertex).
      The other parts are treated in Proposition~\ref{links:vertical}
      and Corollary~\ref{links:horizontal}. Their join is contractible
      or spherical of dimension
      $\DimOf{\TheEuBuild}-\DimOf{\TheEuCell}-\One$. Hence,
      $\DescLinkOf{\CenterOf{\TheEuCell}}$ is contractible or
      spherical of dimension $\DimOf{\TheEuBuild}-\One$.
    \end{proof}

    \begin{observation}\label{links:generically.not.contractible}
      A vertex $\TheEuVertex$ is always a significant cell as
      $\TheEuVertex=\TheEuVertex[][\min]$ since $\TheEuVertex[][\min]$
      is a non-empty face of $\TheEuVertex$.
      Also, a vertex has empty boundary. Hence
      \[
        \DescLinkOf{\CenterOf{\TheEuVertex}}
        =
        \DescLinkOf[\ver]{\TheEuVertex}
        \join
        \DescLinkOf[\hor]{\TheEuVertex}
        \Period
      \]
      Generically, $\LinkOf{\TheEuVertex}$ will not have a horizontal
      component: the gradient $\GradOf[\TheEuVertex]{\TheMorse}$ will
      be in general position. In those cases,
      \(
        \DescLinkOf{\CenterOf{\TheEuVertex}}
        =
        \DescLinkOf[\ver]{\TheEuVertex}
      \)
      is an (open) hemisphere complex, which is \emph{not} contractible
      by Proposition~\ref{bernd:open.hemi}
      since the building $\TheEuBuild$ is thick (in the Rank~Theorem,
      the group $\TheGroupScheme$ is assumed to be non-commutative).
      In particular, there exist arbitrary high vertices with
      non-contractible descending links.\qed
    \end{observation}

  \section{Proof of the Rank Theorem}\label{sec:proof}
    We assume that $\TheZonoDef$ is rich and invariant under the
    full spherical Weyl group. E.g., one could choose $\TheZonoDef$
    to consist of difference vectors of any pair of vertices
    in $\TheModelApp=\TheEuSpace$ whose closed stars intersect.
    \begin{observation}\label{proof:stabilizers}
      The $\ThePlaceSet$-arithmetic group $\TheLattice$ acts
      on the product $\TheEuBuild$ by cell-permuting homeomorphisms.
      Cell stabilizers are finite.\qed
    \end{observation}
    \begin{observation}\label{proof:cocompactness}
      The function $\CombinedHeight$ is $\TheLattice$-invariant by
      Observation~\ref{height:invariant}, and
      its sublevel complexes are $\TheLattice$-cocompact
      by Observation~\ref{height:cocompact}.\qed
    \end{observation}
    \begin{proof}[of the Rank Theorem]
      Given the topological properties of descending links, the
      deduction of finiteness properties is routine.

      Since $\TheLattice$ acts cocompactly, there are only finitely
      many $\TheLattice$-orbits of cells in $\TheEuBuild$ below any
      given $\CombinedHeight$-bound in
      $\TheReals\crossprod\TheReals\crossprod\TheReals$. In particular,
      only finitely many elements in
      $\TheReals\crossprod\TheReals\crossprod\TheReals$
      arise as values of $\CombinedHeight$ below any given bound.
      Define $\TheFiltrationOf{\TheIndex}$ to be the subcomplex
      of $\CenterOf{\TheEuBuild}$ spanned by the barycenters
      $\CenterOf{\TheEuCell}$ of cells $\TheEuCell$ for which
      there are at most $\TheIndex$ values in the image
      $\ImageOf{\CombinedHeight}$ that are strictly below
      $\CombinedHeightOf{\TheEuCell}$.

      By Observation~\ref{morse:no.horizontal.edges}, there are no
      $\CombinedHeight$-flat edges in $\CenterOf{\TheEuBuild}$.
      Thus,
      \(
        \TheFiltrationOf{\TheIndex+\One}
        \setminus\TheFiltrationOf{\TheIndex}
      \)
      does not contain adjacent vertices.
      For any vertex
      \(
        \CenterOf{\TheEuCell}
        \in
        \TheFiltrationOf{\TheIndex+\One}
        \setminus\TheFiltrationOf{\TheIndex}
      \),
      the descending link
      $\DescLinkOf{\CenterOf{\TheEuCell}}$ is precisely the
      relative link
      \(
        \LinkOf{\CenterOf{\TheEuCell}}
        \intersect
        \TheFiltrationOf{\TheIndex}
      \).
      This relative link is contractible or spherical of
      dimension $\DimOf{\TheEuBuild}-\One$ by
      Proposition~~\ref{links:connectivity}.
      Thus, the complex
      $\TheFiltrationOf{\TheIndex+\One}$ is obtained from
      $\TheFiltrationOf{\TheIndex}$ up to homotopy equivalence
      by attaching $\TheDim$-cells -- recall that $\TheDim$ is
      the dimension of $\TheEuBuild$.
      Observation~\ref{links:generically.not.contractible}
      ensures that the extension is nontrivial at infinitely
      many stages.

      The group $\TheLattice$ acts on $\CenterOf{\TheEuBuild}$
      by cell-permuting homeomorphisms and with finite cell
      stabilizers. Thus, all hypotheses of Brown's criterion
      \cite[Corollary~3.3]{Brown:1987} are satisfied and
      $\TheLattice$ is of type $\FType[\TheDim-\One]$ but not
      of type $\FType[\TheDim]$.
    \end{proof}

  \section{Reduction theory: the adelic version}\label{sec:harder}
    In this section, we describe Harder's version of reduction
    theory for reductive groups over global function fields. Thus, we
    relax the hypotheses of the Rank~Theorem: the group scheme
    $\TheGroupScheme$ is assumed to be connected and reductive.
    After Theorem~\ref{harder:cocompact.case}, we shall add the
    requirement that $\TheGroupScheme$ be $\TheGlobalField$-isotropic.

    Let $\TheConstField$ be the finite field of constants of the
    global function field $\TheGlobalField$.
    For any place $\ThePlace$ on $\TheGlobalField$, let
    $\TheLocalField[\ThePlace]$ be the completion of $\TheGlobalField$
    at $\ThePlace$. The field $\TheLocalField[\ThePlace]$ is a local
    field on which we can regard $\ThePlace$ as a normalized discrete
    valuation. Let $\TheOka[\ThePlace]$ be the corresponding valuation
    ring and $\TheMaxIdeal[\ThePlace]$ its unique maximal ideal.
    The residue field
    \(
      \TheResField[\ThePlace]
      :=
      \TheOka[\ThePlace]\rmod\TheMaxIdeal[\ThePlace]
    \)
    is a finite extension of the field $\TheConstField$ of constants.
    Let
    \(
      \TheLocalDeg[\ThePlace]
      :=\ExtDegreeOf{\TheResField[\ThePlace]}{\TheConstField}
    \)
    denote its degree. The \notion{modulus} map
    \begin{eqnarray*}
      \AbsValueOf[\ThePlace]{\DummyArg}
      \mapcolon \TheLocalField[\ThePlace] & \longrightarrow &\TheReals\\
      \TheFieldElement & \mapsto &
      \CardOf{\TheConstField}^{-\TheLocalDeg[\ThePlace]\ThePlaceOf{\TheFieldElement}}
    \end{eqnarray*}
    describes how multiplication by $\TheFieldElement$ changes the
    Haar measure on $\TheLocalField[\ThePlace]$.

    For any finite set of places $\ThePlaceSet$,
    the product
    \[
      \TheAdeles[\ThePlaceSet]
      :=
      \Product[\ThePlace\in\ThePlaceSet]{
        \TheLocalField[\ThePlace]
      }
      \,\,\,\crossprod\,\,\,
      \Product[\ThePlace\not\in\ThePlaceSet]{
        \TheOka[\ThePlace]
      }
    \]
    is the ring of \notion{$\ThePlaceSet$-adeles}. Note that
    the functor $\ThePlaceSet\mapsto\TheAdeles[\ThePlaceSet]$ is
    a directed system indexed by the family of finite sets of
    places. The ring $\TheAdeles$ of \notion{adeles} is by definition
    the direct limit of this system. As each $\TheAdeles[\ThePlaceSet]$
    is a topological ring, so is $\TheAdeles$, and
    $\TheOka:=
    \TheAdeles[\emptyset]=\Product[\ThePlace]{\TheOka[\ThePlace]}$
    is a compact subring.

    For any adele
    \(
      \TheAdel = \FamOf[\ThePlace]{\TheFieldElement[\ThePlace]}
      \in
      \TheAdeles
    \)
    we define the \notion{idele norm} as
    \[
      \AbsValueOf{\TheAdel}
      :=
      \Product[\ThePlace]{
        \AbsValueOf[\ThePlace]{\TheFieldElement[\ThePlace]}
      }
      \Period
    \]
    Taking logarithms, we obtain:
    \begin{equation}\label{eq:idele.log}
      \LogOf{\AbsValueOf{\TheAdel}}
      =
      \Sum[\ThePlace]{
        -
        \TheLocalDeg[\ThePlace]
        \ThePlaceOf{\TheFieldElement[\ThePlace]}
      }
    \end{equation}
    For any $\TheFieldElement\in\TheGlobalField$, there are only
    finitely many places $\ThePlace$ for which
    $\TheFieldElement\not\in\TheOka[\ThePlace]$. Hence,
    $\TheGlobalField$ diagonally embeds into $\TheAdeles$,
    and with respect to this inclusion,
    \(
      \TheOkaRing = \TheAdeles[\ThePlaceSet]\intersect
      \TheGlobalField
    \).
    Also, the idele norm is trivial on $\UnitsOf{\TheGlobalField}$,
    i.e., we have the \notion{product formula}
    \begin{equation}\label{eq:product.formula}
      \AbsValueOf{\TheFieldElement}
      =
      \Product[\ThePlace]{
        \AbsValueOf[\ThePlace]{\TheFieldElement}
      }
      =
      \One
      \qquad\text{for any\ }
      \TheFieldElement\in\UnitsOf{\TheGlobalField}
      \Period
    \end{equation}

    Let $\ExtGlobalField$ be a finite Galois extension of
    $\TheGlobalField$. In particular, $\ExtGlobalField$ is
    a global function field in its own right. Let
    $\ExtAdeles$ denote the ring of adeles associated to
    $\ExtGlobalField$. Since every discrete valuation on
    $\TheGlobalField$ extends to at least one valuation
    on $\ExtGlobalField$, there is a diagonal embedding
    $\TheAdeles\subseteq\ExtAdeles$. Let
    \(
      \ExtNorm \mapcolon \ExtGlobalField\rightarrow\TheGlobalField
    \)
    denote the norm map. For any
    $\ExtFieldElement\in\ExtGlobalField$ one has
    \begin{equation}\label{eq:extending.valuations}
      \AbsValueOf[\ThePlace]{ \ExtNormOf{\ExtFieldElement} }
      =
      \Product[\ExtPlace\text{\ extends\ }\ThePlace]{
        \AbsValueOf[\ExtPlace]{\ExtFieldElement}
      }
      \Period
    \end{equation}
    As the idele norm is defined in terms of the modulus maps,
    we infer:
    \[
      \AbsValueOf[\TheGlobalField]{ \ExtNormOf{\ExtFieldElement} }
      =
      \AbsValueOf[\ExtGlobalField]{\ExtFieldElement}
      \Period
    \]

    Let $\TheGroupScheme$ be a reductive group defined over
    $\TheGlobalField$. Then
    \(
      \TheGroupSchemeOf{\TheAdeles[\ThePlaceSet]}
      =
      \Product[\ThePlace\in\ThePlaceSet]{
        \TheGroupSchemeOf{\TheLocalField[\ThePlace]}
      }
      \,\,\,\crossprod\,\,\,
      \Product[\ThePlace\in\ThePlaceSet]{
        \TheGroupSchemeOf{\TheOka[\ThePlace]}
      }
    \)
    and
    \(
      \TheGroupSchemeOf{\TheAdeles}
      =
      \DirLimit[\ThePlaceSet]{
        \TheGroupSchemeOf{\TheAdeles[\ThePlaceSet]}
      }
    \).
    Let $\Mult$
    denote the multiplicative group regarded as a group scheme
    over $\TheGlobalField$. A \notion{character} on a
    $\TheGlobalField$-group is a homomorphism into $\Mult$.
    One defines:
    \begin{eqnarray*}
      \CharKernel{
        \TheGroupSchemeOf{\TheAdeles}
      }
      &:=&
      \SetOf[
        \TheAdeleGroupElement\in\TheGroupSchemeOf{\TheAdeles}
      ]{
        \AbsValueOf{
          \TheCharacterOf{
            \TheAdeleGroupElement
          }
        }
        =
        \One
        \text{\ for any character\ }
        \TheCharacter\mapcolon\TheGroupScheme
        \rightarrow\Mult
        \text{\ defined over\ }
        \TheGlobalField
      }
      \\
      \CharKernel{
        \TheGroupSchemeOf{\ExtAdeles}
      }
      &:=&
      \SetOf[
        \TheAdeleGroupElement\in\TheGroupSchemeOf{\ExtAdeles}
      ]{
        \AbsValueOf{
          \TheCharacterOf{
            \TheAdeleGroupElement
          }
        }
        =
        \One
        \text{\ for any character\ }
        \TheCharacter\mapcolon\TheGroupScheme
        \rightarrow\Mult
        \text{\ defined over\ }
        \ExtGlobalField
      }
    \end{eqnarray*}
    Of course, there may be more characters defined over $\ExtGlobalField$
    than there are defined over $\TheGlobalField$. Hence the latter group
    appears smaller in this regard. However
    using the norm map $\ExtNorm$ to average
    over $\ExtGlobalField$-characters on $\TheGroupScheme$,
    one can deduce from~(\ref{eq:extending.valuations}) that
    the inclusion
    \(
      \TheGroupSchemeOf{\TheAdeles}
      \subgroup
      \TheGroupSchemeOf{\ExtAdeles}
    \)
    induced by $\TheAdeles\subseteq\ExtAdeles$ restricts to
    an inclusion of
    \(
      \CharKernel{
        \TheGroupSchemeOf{\TheAdeles}
      }
    \)
    in
    \(
      \CharKernel{
        \TheGroupSchemeOf{\ExtAdeles}
      }
    \)
    as a closed topological subgroup.
    \begin{lemma}\label{harder:2.2.3}
      The inclusion
      \(
        \TheGroupSchemeOf{\TheAdeles}\subseteq
        \TheGroupSchemeOf{\ExtAdeles}
      \)
      induces proper maps
      \(
        \TheGroupSchemeOf{\TheAdeles}
        \rmod
        \TheGroupSchemeOf{\TheGlobalField}
        \rightarrow
        \TheGroupSchemeOf{\ExtAdeles}
        \rmod
        \TheGroupSchemeOf{\ExtGlobalField}
      \)
      and
      \(
        \CharKernel{\TheGroupSchemeOf{\TheAdeles}}
        \rmod
        \TheGroupSchemeOf{\TheGlobalField}
        \rightarrow
        \CharKernel{\TheGroupSchemeOf{\ExtAdeles}}
        \rmod
        \TheGroupSchemeOf{\ExtGlobalField}
      \).
    \end{lemma}
    \begin{proof}
      This follows from \cite[Lemma~2.2.3]{Harder:1969}.
    \end{proof}

    The following statement says everything there is to say
    (from the reduction theory point of view)
    about $\TheGlobalField$-anisotropic groups:
    \begin{theorem}[{\cite[Korollar~2.2.7]{Harder:1969}}]\label{harder:cocompact.case}
      $\TheGroupScheme$ is $\TheGlobalField$-anisotropic
      if and only if
      $\CharKernel{\TheGroupSchemeOf{\TheAdeles}}\rmod
      \TheGroupSchemeOf{\TheGlobalField}$ is compact.\qed
    \end{theorem}
    From now on, we assume that $\TheGroupScheme$ is
    $\TheGlobalField$-isotropic.

    Note that $\TheGroupSchemeOf{\TheOka[\ThePlace]}$ is
    an open compact subgroup
    of $\TheGroupSchemeOf{\TheLocalField[\ThePlace]}$.
    Following Harder, we call a subgroup $\TheCompactSubgroup$
    of $\TheGroupSchemeOf{\TheAdeles}$ \notion{standard} if
    $\TheCompactSubgroup$ is of the form
    $\Product[\ThePlace]{\TheCompactSubgroup[\ThePlace]}$
    where each $\TheCompactSubgroup[\ThePlace]$ is an
    open compact subgroup of
    $\TheGroupSchemeOf{\TheLocalField[\ThePlace]}$.
    In particular, the canonical subgroup
    \(
      \TheGroupSchemeOf{\TheOka}
      =
      \Product[\ThePlace]{
        \TheGroupSchemeOf{\TheOka[\ThePlace]}
      }
    \)
    is standard. Let $\TheParabolic$ be a
    $\TheGlobalField$-parabolic subgroup with unipotent
    radical $\UniRadical$. Starting with a non-vanishing volume form
    $\TheHarderForm$ on $\UniRadical$ (in the sense of algebraic
    geometry and defined over $\TheGlobalField$), the associated
    measure
    $\TheHarderMeasure$ on $\UniRadicalOf{\TheAdeles}$ is independent
    of $\TheHarderForm$ because of the product
    formula \cite[Theorem~2.3.1]{Weil:1982}; in fact,
    $\TheHarderMeasure$ is proportional to the Tamagawa measure.
    Harder defines for any parabolic $\TheParabolic$ and any
    standard subgroup $\TheCompactSubgroup$:
    \begin{equation}\label{eq:harder.pi}
      \HarderPiOf{\TheParabolic,\TheCompactSubgroup}
      :=
      \VolumeOf[\TheHarderMeasure]{
        \UniRadicalOf{\TheAdeles}\intersect\TheCompactSubgroup
      }
    \end{equation}
    As the measure $\TheHarderMeasure$ is canonical, this definition
    is invariant under the conjugacy
    action of $\TheGroupSchemeOf{\TheGlobalField}$ on
    $\TheGroupSchemeOf{\TheAdeles}$, i.e., for each
    element $\TheGroupElement\in\TheGroupSchemeOf{\TheGlobalField}$
    we have:
    \begin{equation}\label{eq:invariance}
      \HarderPiOf{\TheParabolic,\TheCompactSubgroup}
      =
      \HarderPiOf{\Conj{\TheGroupElement}{
                    \TheParabolic
                  }
                  ,
                  \Conj{\TheGroupElement}{
                    \TheCompactSubgroup
                  }}
    \end{equation}

    The unipotent radical $\UniRadical$ is a weight space for the
    adjoint representation of the parabolic group $\TheParabolic$.
    We call the associated character
    \(
      \TheHarderChar[\TheParabolic]
      \mapcolon\TheParabolic\rightarrow\Mult
    \)
    the \notion{canonical character} of $\TheParabolic$.
    Its idele norm
    is the functional determinant of the conjugacy action of
    $\TheParabolic$ on $\UniRadical$. Hence:
    \begin{prop}[Transformation Formula {\cite[Satz~1.3.2]{Harder:1969}}]\label{harder:transformation.formula}
      For any standard subgroup
      \(
        \TheCompactSubgroup\subgroup
        \TheGroupSchemeOf{\TheAdeles}
      \)
      and any
      $\TheAdeleGroupElement\in\TheParabolicOf{\TheAdeles}$,
      we have
      \begin{equation}\label{eq:transformation.formula}
        \HarderPiOf{\TheParabolic,\TheCompactSubgroup}
        =
        \HarderPiOf{\TheParabolic,
                    \Conj{\TheAdeleGroupElement}{\TheCompactSubgroup}}
        \AbsValueOf{
          \TheHarderCharOf[\TheParabolic]{\TheAdeleGroupElement}
        }
        \Period
      \end{equation}
    \end{prop}

    \begin{construction}[{\cite[page~47]{Harder:1969}}]\label{harder:construction}
    Assume that $\TheParabolic$ is a
    minimal $\TheGlobalField$-parabolic,
    let $\TheRadical$
    be its radical, $\UniRadical$ be its unipotent radical, and put
    \(
      \TheTorus:=\TheRadical\rmod\UniRadical
    \).
    Let $\TheSplitTorus\subgroup\TheTorus$ be the maximal
    $\TheGlobalField$-split torus.
    We think of $\TheTorus$ and $\TheSplitTorus$ not just as
    abstract tori but as tori inside of a Levi subgroup $\TheLevy$
    of $\TheParabolic$. In particular, $\TheTorus$ is a maximal
    torus in $\TheLevy$ and $\TheSplitTorus$ is a maximal
    $\TheGlobalField$-split torus inside $\TheLevy$.
    Let
    \(
      \SetOf{
        \TheRoot[\One],\ldots,
        \TheRoot[\TheLastRoot]
      }
      \subset
      \CharactersOf{\TheSplitTorus}
    \)
    be the set of the simple roots on $\TheSplitTorus$. With
    $\CharactersOf{\TheParabolic}:=
    \HomOf[\TheGlobalField]{\TheParabolic;\Mult}$, we have
    \(
      \CharactersOf{\TheParabolic}\tensor\TheReals
      =
      \CharactersOf{\TheSplitTorus}\tensor\TheReals
    \)
    hence, we can regard each $\TheRoot[\TheRootIndex]$ as an
    element of $\CharactersOf{\TheParabolic}\tensor\TheReals$.
    The minimal parabolic $\TheParabolic$ corresponds to a chamber
    of the spherical
    building $\TheRatBuild=\TheRatBuild[\TheGlobalField]$.
    The roots $\TheRoot[\TheRootIndex]$ correspond to faces.
    Hence $\TheRootIndex$ can be regarded as a cotype. Let
    $\TheParabolic[\TheRootIndex]$ be the maximal parabolic
    above $\TheParabolic$ of type~$\TheRootIndex$ (i.e.,
    the face corresponding to $\TheRoot[\TheRootIndex]$ and
    the vertex corresponding to $\TheParabolic[\TheRootIndex]$
    span the chamber for $\TheParabolic$). Let
    \(
      \TheHarderChar[\TheRootIndex]
      \mapcolon
      \TheParabolic\rightarrow\Mult
    \)
    be the restriction of the canonical character
    \(
      \TheHarderChar[{\TheParabolic[\TheRootIndex]}]
      \mapcolon
      \TheParabolic[\TheRootIndex]
      \rightarrow
      \Mult
    \).

    The set of roots
    \(
      \SetOf{
        \TheRoot[\One],\ldots,
        \TheRoot[\TheLastRoot]
      }
    \)
    is a basis for
    \(
      \CharactersOf{\TheParabolic}\tensor\TheReals
    \)
    and so is the set
    \(
      \SetOf{
        \TheHarderChar[\One],\ldots,
        \TheHarderChar[\TheLastRoot]
      }
    \).
    This determines real (in fact rational) numbers
    $\TheBaseChange[\TheCotype\TheType]$ and
    $\InvBaseChange[\TheType\TheCotype]$ such that:
    \begin{eqnarray*}
      \TheRoot[\TheCotype]
      &=&
      \Sum[\TheType]{
        \TheBaseChange[\TheCotype\TheType]
        \TheHarderChar[\TheType]
      }
      \\
      \TheHarderChar[\TheType]
      &=&
      \Sum[\TheCotype]{
        \InvBaseChange[\TheType\TheCotype]
        \TheRoot[\TheCotype]
      }
    \end{eqnarray*}
    These two bases are almost dual. Let
    $\InnerProd{\DummyArg}{\DummyArg}$ be an inner product
    on
    \(
      \CharactersOf{\TheParabolic}\tensor\TheReals
    \)
    invariant under the action of the Weyl group. The
    $\TheRoot[\TheCotype]$ are simple roots and the
    $\TheHarderChar[\TheType]$ point in the direction of
    the fundamental weights. Thus:
    \begin{eqnarray*}
      \Zero
      &\leq&
      \InvBaseChange[\TheType\TheCotype]
      \\
      \Zero
      &<&
      \InvBaseChange[\TheType\TheType]
      \\
      \Zero
      &<&
      \TheBaseChange[\TheType\TheType]
      \\
      \InnerProd{\TheHarderChar[\TheType]}{\TheRoot[\TheCotype]}
      &=&
      \Zero
      \qquad
      \text{if\ }\TheType\neq\TheCotype
      \\
      \InnerProd{\TheHarderChar[\TheType]}{\TheRoot[\TheCotype]}
      &>&
      \Zero
      \qquad
      \text{if\ }\TheType=\TheCotype
      \\
      \InnerProd{\TheRoot[\TheCotype]}{\TheRoot[\AltCotype]}
      &\leq&
      \Zero
      \qquad\text{for all\ }\TheCotype,\AltCotype
      \\
      \InnerProd{\TheHarderChar[\TheType]}{\TheHarderChar[\AltType]}
      &\geq&
      \Zero
      \qquad\text{for all\ }\TheType,\AltType
    \end{eqnarray*}
    \end{construction}
    \begin{observation}\label{harder:geometry}
      Let $\TheAdmissibleSet$ be a subset of
      \[
        \SetOf[{
          \TupelOf{
            \TheBase[\One],\ldots,\TheBase[\TheLastRoot],
            \DualBase[\One],\ldots,\DualBase[\TheLastRoot]
          }
          \in
          \TheReals[][\Two\TheLastRoot]
        }]{
          \DualBase[\TheType]
          =
          \Sum[\TheCotype]{
            \InvBaseChange[\TheType\TheCotype]
            \TheBase[\TheCotype]
          }
        }
      \]
      with the coefficients $\InvBaseChange[\TheType\TheCotype]$
      as above. The $\DualBase[\TheType]$ depend on the
      $\TheBase[\TheCotype]$. Hence, $\TheAdmissibleSet$ is bounded
      if and only if its projection onto the first $\TheLastRoot$
      coordinates is. Moreover, the coefficients
      $\InvBaseChange[\TheType\TheCotype]$ are non-negative and
      strictly positive for $\TheType=\TheCotype$. Hence,
      $\DualBase[\TheCotype]$ tends to $\infty$ if
      $\TheBase[\TheCotype]$ tends to $\infty$ while all other
      $\TheBase[\AltCotype]$ stay bounded from below. Thus, the
      following are equivalent:
      \begin{enumerate}
        \item
          $\TheAdmissibleSet$ is bounded.
        \item
          There exists constants
          \(
            \Const[\TheBase][-]
          \)
          and
          \(
            \Const[\TheBase][+]
          \)
          with
          \[
            \Const[\TheBase][-] \leq
            \TheBase[\TheCotype] \leq
            \Const[\TheBase][+]
            \qquad\text{for all\ }
            \TheCotype
          \]
          for all
          \(
            \TupelOf{
              \TheBase[\One],\ldots,\TheBase[\TheLastRoot],
              \DualBase[\One],\ldots,\DualBase[\TheLastRoot]
            }
            \in
            \TheAdmissibleSet
          \).
        \item
          There exists constants
          \(
            \Const[\TheBase][-]
          \)
          and
          \(
            \Const[\DualBase][+]
          \)
          with
          \[
            \Const[\TheBase][-] \leq
            \TheBase[\TheCotype],\,\,\,
            \DualBase[\TheType] \leq
            \Const[\DualBase][+]
            \qquad\text{for all\ }
            \TheCotype,\TheType
          \]
          for all
          \(
            \TupelOf{
              \TheBase[\One],\ldots,\TheBase[\TheLastRoot],
              \DualBase[\One],\ldots,\DualBase[\TheLastRoot]
            }
            \in
            \TheAdmissibleSet
          \).\qed
      \end{enumerate}
    \end{observation}

    Recall, that $\TheParabolic[\TheType]$ denotes the maximal
    parabolic of type $\TheType$ containing the minimal parabolic
    $\TheParabolic$.
    Harder defines the invariants
    \[
      \HarderNuOf[\TheCotype]{
        \TheParabolic,\TheCompactSubgroup
      }
      :=
      \Product[\TheType]{
        \HarderPiOf{
          \TheParabolic[\TheType],
          \TheCompactSubgroup
        }^{\TheBaseChange[\TheCotype\TheType]}
      },
    \]
    but we find it more convenient to express his results using
    logarithms:
    \begin{eqnarray*}
      \LogPiOf{\TheParabolic[\TheType],\TheCompactSubgroup}
      &:=&
      \LogOf{
        \HarderPiOf{
          \TheParabolic[\TheType],
          \TheCompactSubgroup
        }
      }
      \\
      \LogNuOf[\TheCotype]{
        \TheParabolic,\TheCompactSubgroup
      }
      &:=&
      \LogOf{
        \HarderNuOf[\TheCotype]{
          \TheParabolic,\TheCompactSubgroup
        }
      }
    \end{eqnarray*}
    Note that
    \(
        \HarderPiOf{
          \TheParabolic[\TheType],
          \TheCompactSubgroup
        }
        >
        \Zero
    \).
    Now, we have
    \[
      \LogNuOf[\TheCotype]{
        \TheParabolic,\TheCompactSubgroup
      }
      =
      \Sum[\TheType]{
        \TheBaseChange[\TheCotype\TheType]
        \LogPiOf{\TheParabolic[\TheType],\TheCompactSubgroup}
      };
    \]
    and we say that a constant $\HarderConstOne$ is a
    \notion{lower reduction bound} if
      for any
      $\TheAdeleGroupElement\in\TheGroupSchemeOf{\TheAdeles}$, there
      exists a minimal $\TheGlobalField$-parabolic
      subgroup $\TheParabolic$ satisfying
      \(
        \LogNuOf[\TheCotype]{
          \TheParabolic,
          \Conj{\TheAdeleGroupElement}{
            \TheGroupSchemeOf{\TheOka}
          }
        }
        \geq
        \HarderConstOne
      \)
      for all $\TheCotype$.
    In this language, the main theorems of reduction theory
    read as follows:
    \begin{theorem}[{\cite[Satz~2.3.2]{Harder:1969}}]\label{harder:one}
      If $\TheGroupScheme$ is $\TheGlobalField$-isotropic it admits
      a lower reduction bound.\qed
    \end{theorem}
    For a minimal parabolic $\TheParabolic$ and an element
    $\TheAdeleGroupElement\in\TheGroupSchemeOf{\TheAdeles}$, we
    say that the parabolic \notion{reduces} $\TheAdeleGroupElement$
    \notion{with bound $\HarderConstOne$}, if
    \(
      \LogNuOf[\TheCotype]{
        \TheParabolic,
        \Conj{\TheAdeleGroupElement}{
          \TheGroupSchemeOf{\TheOka}
        }
      }
      \geq
      \HarderConstOne
    \)
    for all $\TheCotype$. We may not specify the bound if it is clear
    from the context.
    \begin{theorem}[{\cite[Satz~2.3.3]{Harder:1969}}]\label{harder:two}
      Assume that $\TheGroupScheme$ is $\TheGlobalField$-isotropic.
      For any lower reduction bound $\HarderConstOne$
      there is another constant $\HarderConstTwo$
      (which we call the \notion{upper reduction bound})
      such that:
      whenever $\TheParabolic$ is a minimal $\TheGlobalField$-parabolic
      reducing $\TheAdeleGroupElement\in
      \TheGroupSchemeOf{\TheAdeles}$ with bound $\HarderConstOne$
      and
      \(
        \LogNuOf[\TheCotype]{
          \TheParabolic,
          \Conj{\TheAdeleGroupElement}{
            \TheGroupSchemeOf{\TheOka}
          }
        }
        \geq
        \HarderConstTwo
      \)
      then \emph{any} minimal $\TheGlobalField$-parabolic that
      reduces $\TheAdeleGroupElement$
      is contained
      in the maximal $\TheGlobalField$-parabolic of type
      $\TheCotype$ above $\TheParabolic$.\qed
    \end{theorem}
    \begin{theorem}[Mahler's Compactness Criterion]\label{harder:mahler}
      Assume that $\TheGroupScheme$ is $\TheGlobalField$-isotropic.
      A subset $\TheMahlerSet\subseteq
      \CharKernel{\TheGroupSchemeOf{\TheAdeles}}$ is relatively
      compact modulo $\TheGroupSchemeOf{\TheGlobalField}$ if any
      only if there are two constants
      $\TheLowerCompactnessBound$ and
      $\TheUpperCompactnessBound$ such that for every
      $\TheMahlerElement\in\TheMahlerSet$ there exists a
      minimal $\TheGlobalField$-parabolic subgroup
      $\TheParabolic$ with
      \[
        \TheLowerCompactnessBound
        \leq
        \LogNuOf[\TheCotype]{
          \TheParabolic,
          \Conj{\TheMahlerElement}{
            \TheGroupSchemeOf{\TheOka}
          }
        }
        \leq
        \TheUpperCompactnessBound
      \]
      for each $\TheCotype$. Without loss of generality, the lower bound
      $\TheLowerCompactnessBound$ can be taken to be any lower reduction
      bound.
    \end{theorem}
    It is a little unfortunate that Harder states
    Theorem~\ref{harder:mahler} only in the case that
    $\TheGroupScheme$ is $\TheGlobalField$-split. Harder also
    provides the means of deducing the non-split case, but he
    does not carry out the argument. We provide an outline.
    \begin{proof}[of Theorem~\ref{harder:mahler}]
      For a fixed
      $\TheAdeleGroupElement\in\TheGroupSchemeOf{\TheAdeles}$,
      consider the set
      \[
        \TheLogNuSet=
        \TheLogNuSetOf{\TheAdeleGroupElement}
        :=
        \SetOf[{
          \LogNuOf[\TheCotype]{
            \TheParabolic,
            \Conj{\TheAdeleGroupElement}{
              \TheGroupSchemeOf{\TheOka}
            }
          }
        }]{
          \TheParabolic\text{\ minimal\ }
          \TheGlobalField\text{-parabolic reducing\ }
          \TheAdeleGroupElement\text{,\ }
          \TheCotype\text{\ arbitrary}
        }
        \subseteq
        \TheReals
      \]
      We claim that this set is bounded.
      It is bounded from below
      since the parabolics $\TheParabolic$ are assumed to reduce with
      respect to a fixed lower reduction bound $\HarderConstOne$.

      There are only finitely many types. So, assuming that
      $\TheLogNuSet$ is not bounded from above, there is
      a $\TheCotype$ such that
      \[
        \TheLogNuSet[\TheCotype]
        :=
        \SetOf[{
          \LogNuOf[\TheCotype]{
            \TheParabolic,
            \Conj{\TheAdeleGroupElement}{
              \TheGroupSchemeOf{\TheOka}
            }
          }
        }]{
          \TheParabolic\text{\ minimal\ }
          \TheGlobalField\text{-parabolic reducing\ }
          \TheAdeleGroupElement
        }
      \]
      is unbounded. Observation~\ref{harder:geometry} then
      implies that
      \[
        \TheLogPiSet[\TheCotype]
        :=
        \SetOf[{
          \LogPiOf{
            \TheParabolic[\TheCotype],
            \Conj{\TheAdeleGroupElement}{
              \TheGroupSchemeOf{\TheOka}
            }
          }
        }]{
          \TheParabolic[\TheCotype]\text{\ maximal\ }
          \TheGlobalField\text{-parabolic of type\ }
          \TheCotype\text{\ containing a\ }
          \TheParabolic\text{\ reducing\ }
          \TheAdeleGroupElement
        }
      \]
      is unbounded. This, however, contradicts
      Theorem~\ref{harder:two}. Hence,
      $\TheLogNuSetOf{\TheAdeleGroupElement}$ is
      bounded.

      Now assume that $\TheMahlerSet\subseteq
      \CharKernel{\TheGroupSchemeOf{\TheAdeles}}$ is relatively
      compact modulo $\TheGroupSchemeOf{\TheGlobalField}$.
      Then
      \(
        \Union[\TheAdeleGroupElement\in\TheMahlerSet]{
          \TheLogNuSetOf{\TheAdeleGroupElement}
        }
      \)
      is still bounded. The constant
      $\HarderConstOne$
      can be taken as
      $\TheLowerCompactnessBound$ and the upper bound can be
      taken as
      $\TheUpperCompactnessBound$.

      To argue the converse, we let $\ExtGlobalField$ be a finite
      Galois extension of
      $\TheGlobalField$ such that $\TheGroupScheme$ is
      $\ExtGlobalField$-split.
      Let $\TheMahlerSet$ be a subset of
      $\CharKernel{\TheGroupSchemeOf{\TheAdeles}}$.
      By Lemma~\ref{harder:2.2.3}, $\TheMahlerSet$ is relatively
      compact modulo $\TheGroupSchemeOf{\TheGlobalField}$
      in
      \(
        \CharKernel{\TheGroupSchemeOf{\TheAdeles}}
      \)
      if it is relatively compact
      modulo $\TheGroupSchemeOf{\ExtGlobalField}$ in
      \(
        \CharKernel{\TheGroupSchemeOf{\ExtAdeles}}
      \).
      As Harder argues in \cite[Lemma~2.2.2]{Harder:1969},
      this happens if there are two constants
      \(
        \TheUpperCompactnessBound'
        \geq
        \TheLowerCompactnessBound'
        >
        \Zero
      \)
      such that
      for each $\TheMahlerElement\in\TheMahlerSet$ and each
      $\ExtGlobalField$-Borel subgroup $\TheBorel$ the inequality
      \[
        \TheLowerCompactnessBound'
        \leq
        \HarderNuOf[\TheCotype,\AltIndex]{
          \TheBorel,
          \Conj{\TheMahlerElement}{
            \TheStandardSubgroup'
          }
        }
        \leq
        \TheUpperCompactnessBound'
      \]
      holds where $\TheStandardSubgroup'$ is a suitable standard
      subgroup in $\TheGroupSchemeOf{\ExtAdeles}$.

      Given bounds
      \(
        \TheLowerCompactnessBound
      \)
      and
      \(
        \TheUpperCompactnessBound
      \)
      as in the statement of Theorem~\ref{harder:mahler},
      one can find such
      \(
        \TheLowerCompactnessBound'
      \)
      and
      \(
        \TheUpperCompactnessBound'
      \)
      using \cite[Lemma~2.3.5]{Harder:1969}. We remark that
      the exponent $n$ in that statement is
      the degree of the extension $\ExtGlobalField/\TheGlobalField$,
      see \cite[Lemma~2.2.6]{Harder:1969}.
    \end{proof}
    The following alternate form of Mahler's compactness criterion
    is a consequence of Observation~\ref{harder:geometry}:
    \begin{cor}[Mahler's Compactness Criterion, alternate form]\label{harder:alternate}
      Let $\HarderConstOne$ be a lower reduction bound.
      A subset $\TheMahlerSet\subseteq
      \CharKernel{\TheGroupSchemeOf{\TheAdeles}}$ is relatively
      compact modulo $\TheGroupSchemeOf{\TheGlobalField}$ if any
      only if there is a constant
      $\TheUpperCompactnessBound$ such that for every
      $\TheMahlerElement\in\TheMahlerSet$ there exists a
      minimal $\TheGlobalField$-parabolic subgroup
      $\TheParabolic$ that reduces $\TheMahlerElement$ with
      bound $\HarderConstOne$ and satisfies
      \(
        \LogPiOf{
          \TheParabolic[\TheType],
          \Conj{\TheMahlerElement}{
            \TheGroupSchemeOf{\TheOka}
          }
        }
        \leq
        \TheUpperCompactnessBound
      \)
      for each $\TheType$.\qed
    \end{cor}

  \section{Geometric reduction theory}%
  \label{sec:us}
    In this section,
    $\TheGroupScheme$ is assumed to be connected and reductive.
    The group $\TheGroupSchemeOf{\TheLocalField[\ThePlace]}$ acts
    on the associated euclidean Bruhat-Tits
    building $\TheEuBuild[\ThePlace]$. The action is not necessarily
    type-preserving, but it is transitive on chambers;
    in particular, it has only finitely many orbits of vertices.
    The subgroup $\TheGroupSchemeOf{\TheOka[\ThePlace]}$ is
    the stabilizer of some vertex. The group
    $\TheGroupSchemeOf{\TheAdeles[\ThePlaceSet]}$ acts componentwise
    on the product
    \(
      \TheEuBuild:=\Product[\ThePlace\in\ThePlaceSet]{
        \TheEuBuild[\ThePlace]
      }
    \)
    (components corresponding to places not in $\ThePlaceSet$ act
    trivially).
    The subgroup $\TheGroupSchemeOf{\TheOka}$ is the stabilizer
    of some \ignore{poly-}vertex $\TheBasePolyVertex$ in $\TheEuBuild$.
    There are only finitely many
    $\TheGroupSchemeOf{\TheAdeles[\ThePlaceSet]}$-orbits of
    \ignore{poly-}vertices in $\TheEuBuild$; hence,
    there is a uniform upper bound
    for the distance of any point in $\TheEuBuild$ to the orbit
    $\TheGroupSchemeOf{\TheAdeles[\ThePlaceSet]}\acts\TheBasePolyVertex$.
    Heuristically, the translation of reduction theory into the language
    of buildings proceeds via \emph{pretending} that the euclidean
    building $\TheEuBuild$ can be identified with the orbit space
    \(
      \TheGroupSchemeOf{\TheAdeles[\ThePlaceSet]}
      \rmod
      \TheGroupSchemeOf{\TheOka}
    \).

    To make this more precise, let $\TheBaseVertex[\ThePlace]$ the
    vertex in $\TheEuBuild[\ThePlace]$ stabilized by the group
    $\TheGroupSchemeOf{\TheOka[\ThePlace]}$, and let
    \(
      \TheSpecVertexSet[\ThePlace]
    \)
    denote the $\TheGroupSchemeOf{\TheLocalField[\ThePlace]}$-orbit
    of $\TheBaseVertex[\ThePlace]$. Hence
    \[
      \TheSpecVertexSet[\ThePlace]
      =
      \TheGroupSchemeOf{\TheLocalField[\ThePlace]}
      \rmod
      \TheGroupSchemeOf{\TheOka[\ThePlace]}
      \text{.}
    \]
    Putting
    \(
      \TheSpecVertexSet :=
      \Product[\ThePlace\in\ThePlaceSet]{
        \TheSpecVertexSet[\ThePlace]
      }
    \),
    we have:
    \begin{eqnarray*}
      \TheSpecVertexSet
      &=&
      \Product[\ThePlace\in\ThePlaceSet]{
        \TheGroupSchemeOf{\TheLocalField[\ThePlace]}
        \rmod
        \TheGroupSchemeOf{\TheOka[\ThePlace]}
      }
      \\
      &=&
      \Product[\ThePlace\in\ThePlaceSet]{
        \TheGroupSchemeOf{\TheLocalField[\ThePlace]}
        \rmod
        \TheGroupSchemeOf{\TheOka[\ThePlace]}
      }
      \,\,\,\crossprod\,\,\,
      \Product[\ThePlace\not\in\ThePlaceSet]{
        \TheGroupSchemeOf{\TheOka[\ThePlace]}
        \rmod
        \TheGroupSchemeOf{\TheOka[\ThePlace]}
      }
      \\
      &=&
      \TheGroupSchemeOf{\TheAdeles[\ThePlaceSet]}
      \rmod
      \TheGroupSchemeOf{\TheOka}
    \end{eqnarray*}
    Conversely, for any \ignore{poly-}vertex
    $\TheSpecVertex\in\TheSpecVertexSet$, the stabilizer
    $\StabOf{\TheSpecVertex}$
    in $\TheGroupSchemeOf{\TheAdeles[\ThePlaceSet]}$ is
    a standard subgroup of $\TheGroupSchemeOf{\TheAdeles}$.

    We can now start to interpret reduction theory in terms
    of Busemann functions. Let $\TheRatBuild$ be the spherical
    building of $\TheGroupSchemeOf{\TheGlobalField}$ over
    the global field, i.e., the simplicial complex that is the
    realization of the poset of proper $\TheGlobalField$-parabolic
    subgroups of $\TheGroupScheme$.
    Any vertex $\TheRatVertex\in\VerticesOf{\TheRatBuild}$
    corresponds to a maximal $\TheGlobalField$-parabolic
    $\TheParabolic[\TheRatVertex]$ of $\TheGroupScheme$. In
    particular, the building $\TheRatBuild$ is empty if and only
    if $\TheGroupScheme$ is anisotropic over $\TheGlobalField$.
    The anisotropic case is implicitly excluded in all considerations
    that require $\TheRatBuild$ to be non-empty. Note, however, that
    any statement of the form ``for any vertex $\TheRatVertex$ in
    $\TheRatBuild$, \ldots'' is vacuously true.

    For $\TheRatVertex\in\VerticesOf{\TheRatBuild}$, we define:
    \begin{eqnarray*}
      \AlmostBusemann[\TheRatVertex] \mapcolon
      \TheSpecVertexSet & \longrightarrow & \TheReals \\
      \TheSpecVertex
      &\mapsto&
      \LogPiOf{\TheParabolic[\TheRatVertex],\StabOf{\TheSpecVertex}}
    \end{eqnarray*}
    We would like to show that $\AlmostBusemann$ can be extended to
    a Busemann function on $\TheEuBuild$.

      Let $\TheParabolic$ be a minimal $\TheGlobalField$-parabolic
      subgroup of $\TheGroupScheme$ corresponding to a chamber
      $\TheRatChamber$ of $\TheRatBuild$.
      By \cite[Theorem~13.3.6]{Springer:1998}, the group $\TheParabolic$
      contains a maximal torus $\TheTorus$
      that is defined over $\TheGlobalField$.
      Of course, $\TheTorus$ is not necessarily split
      over $\TheGlobalField$.  Let $\TheSplitTorus$ be the maximal
      $\TheGlobalField$-split subtorus of $\TheTorus$. For each
      place $\ThePlace\in\ThePlaceSet$, let $\TheSplitTorus[\ThePlace]$
      be the maximal $\TheLocalField[\ThePlace]$-subtorus of
      $\TheTorus$. Note that
      \(
        \TheSplitTorus\subgroup\TheSplitTorus[\ThePlace]
      \)
      for each $\ThePlace\in\ThePlaceSet$. Let $\TheEuApp[\ThePlace]$
      be the apartment corresponding to $\TheSplitTorus[\ThePlace]$ in
      the euclidean building $\TheEuBuild[\ThePlace]$. We put
      \(
        \TheEuApp:=\Product[\ThePlace\in\ThePlaceSet]{
          \TheEuApp[\ThePlace]
        }
      \)
      and
      \(
        \TheSpecVertexApp
        :=
        \TheSpecVertexSet\intersect\TheEuApp
      \).
    \begin{lemma}\label{us:apartment}
      For any vertex $\TheRatVertex\in\TheRatChamber$, there exists
      an affine function on $\TheEuApp$ that agrees with
      \(
        \AlmostBusemann[\TheRatVertex]
      \)
      on the set $\TheSpecVertexApp$.
    \end{lemma}
    \begin{proof}
      From the Transformation Formula in
      Proposition~\ref{harder:transformation.formula}
      we obtain
      \begin{equation}\label{eq:busemann.on.apartment}
        \LogPiOf{\TheParabolic[\TheRatVertex],
                    \Conj{\TheAdeleGroupElement}{\TheCompactSubgroup}}
        -
        \LogPiOf{\TheParabolic[\TheRatVertex],\TheCompactSubgroup}
        =
        -
        \LogOf{
          \AbsValueOf{
            \TheHarderCharOf[{\TheParabolic[\TheRatVertex]}]{
              \TheAdeleGroupElement
            }
          }
        }
        =
        \Sum[\ThePlace]{
          \TheLocalDeg[\ThePlace]
          \ThePlaceOf{
            \TheHarderCharOf[{\TheParabolic[\TheRatVertex]}]{
              \TheAdeleGroupElement
            }
          }
        }
      \end{equation}
      for each $\TheAdeleGroupElement\in\TheTorusOf{\TheAdeles}$.
      Considering this statement just for
      \(
        \TheAdeleGroupElement
        \in
        \Product[\ThePlace\in\ThePlaceSet]{
          \TheSplitTorusOf[\ThePlace]{
            \TheLocalField[\ThePlace]
          }
        }
      \)
      with $\TheCompactSubgroup$ taken to be the stabilizer of
      a \ignore{poly-}vertex in $\TheSpecVertexApp$,
      the claim follows.
    \end{proof}
    So far, we cannot speak of Busemann functions on $\TheEuBuild$
    as we did not yet fix a euclidean metric on $\TheEuBuild$. There
    is some freedom in making this choice: on the one hand, we can
    freely rescale metrics on the factors $\TheEuBuild[\ThePlace]$;
    on the other hand, the factors $\TheEuBuild[\ThePlace]$ need
    not be irreducible and if $\TheEuBuild[\ThePlace]$ decomposes
    as a product, the metrics on the irreducible factors can be
    independently scaled. However, that is the only source of
    non-uniqueness: up to scaling, there is a unique Weyl-group
    invariant metric on any irreducible euclidean building.
    In particular, we only have to choose the relative scales of the
    factors $\TheEuBuild[\ThePlace]$ if $\TheGroupScheme$ is
    absolutely almost simple.

    However, even in the case of an absolutely almost simple group
    and a single place $\ThePlaceSet=\SetOf{\ThePlace}$, we would have
    something to prove: we do not just want some metric on $\TheEuBuild$.
    Rather, we would like a metric so that the geometry of the root
    system constructed in~\ref{harder:construction} is reflected in
    the angular metric at infinity induced from the euclidean metric
    on $\TheEuBuild$. The reason for this restriction stems from the
    following: Harder's reduction theory is phrased in terms of the
    roots $\TheRoot[\TheCotype]$ and the dual characters
    $\TheCharacter[\TheType]$. In our translation, we want to dispose
    of the roots $\TheRoot[\TheCotype]$ and only work with the
    characters (to those, our Busemann function will correspond).
    The euclidean metric is supposed to supply the
    necessary duality by means of its associated inner product.
    Thus, we need to demonstrate how this can be achieved.

    As a first step, we compare the root system for $\TheGroupScheme$
    over the global field $\TheGlobalField$ to the root system over
    the local field $\TheLocalField[\ThePlace]$. If the field extension
    $\TheLocalField[\ThePlace]\fieldextends\TheGlobalField$ was normal,
    we could directly quote \cite[{\S\,6}]{Borel.Tits:1965}.
    Let $\TheLocalField[\ThePlace][\sep]$ be the abstract separable
    closure of $\TheLocalField[\ThePlace]$, and
    let $\TheGlobalField[][\sep]$ be the separable closure of
    $\TheGlobalField$ inside $\TheLocalField[\ThePlace][\sep]$. Now
    $\TheLocalField[\ThePlace][\sep]\fieldextends\TheLocalField[\ThePlace]$
    and
    $\TheGlobalField[][\sep]\fieldextends\TheGlobalField$
    are both normal extensions. The group $\TheGroupScheme$ splits
    over $\TheGlobalField[][\sep]$ and $\TheLocalField[\ThePlace][\sep]$.
    Moreover, the root systems for $\TheGroupScheme$ over
    $\TheGlobalField[][\sep]$ and $\TheLocalField[\ThePlace][\sep]$
    are canonically isomorphic: a maximal $\TheGlobalField[][\sep]$-split
    torus $\ExtSplitTorus$ in $\TheGroupScheme$ is also maximal
    $\TheLocalField[\ThePlace][\sep]$-split, and all its characters
    defined over $\TheLocalField[\ThePlace][\sep]$ are already defined
    over $\TheGlobalField[][\sep]$; hence
    \begin{equation}\label{eq:identification}
      \CharactersOf[{\TheGlobalField[][\sep]}]{
        \ExtSplitTorus
      }
      \tensor
      \TheReals
      =
      \CharactersOf[{\TheLocalField[\ThePlace][\sep]}]{
        \ExtSplitTorus
      }
      \tensor
      \TheReals
      =:
      \ExtVectorSpace
      \text{.}
    \end{equation}
    Let $\TheSplitTorus$ be a maximal $\TheGlobalField$-split torus
    in $\ExtSplitTorus$ and put
    \(
      \TheVectorSpace :=
      \CharactersOf[\TheGlobalField]{\TheSplitTorus}
      \tensor\TheReals
    \).
    Restriction of characters on $\ExtSplitTorus$ to
    $\TheSplitTorus$ induces a projection
    \(
      \ExtVectorSpace\rightarrow\TheVectorSpace
    \).
    We endow $\ExtVectorSpace$ with an inner product that is
    invariant under the full spherical Weyl group.
    By \cite[\S\,6.10]{Borel.Tits:1965}, there is a canonical way
    of realizing the abstract vector space $\TheVectorSpace$ as
    a subspace of $\ExtVectorSpace$ such that the orthogonal projection
    is the restriction homomorphism. The induced inner product on
    $\TheVectorSpace$ is invariant under the Weyl group of
    $\TheGroupScheme$ over $\TheGlobalField$.

    The same construction can be carried out for
    $\TheLocalField[\ThePlace]$, yielding a subspace
    $\TheVectorSpace[\ThePlace]$ of $\ExtVectorSpace$. We want
    to argue the inclusion
    \(
      \TheVectorSpace\subspace\TheVectorSpace[\ThePlace]
    \).
    To do so, we impose the assumption that
    the chosen $\TheGlobalField[][\sep]$-split torus
    $\ExtSplitTorus$ is defined over $\TheGlobalField$.
    It is then also defined over $\TheLocalField[\ThePlace]$.
    Moreover, \cite[\S\,6.11]{Borel.Tits:1965} applies:
    the vector space $\TheVectorSpace$ is the fixed point set
    of the Galois action of
    \(
      \GalOf{
        \TheGlobalField[][\sep]
        \fieldextends
        \TheGlobalField
      }
    \)
    on $\ExtVectorSpace$. Similarly,
    $\TheVectorSpace[\ThePlace]$ is the fixed point set
    of
    \(
      \GalOf{
        \TheLocalField[\ThePlace][\sep]
        \fieldextends
        \TheLocalField[\ThePlace]
      }
    \)
    on $\ExtVectorSpace$. Since any $\TheGlobalField$-automorphism
    of $\TheLocalField[\ThePlace][\sep]$ leaves the separable closure
    $\TheGlobalField[][\sep]\subseteq\TheLocalField[\ThePlace][\sep]$
    invariant, we have a homomorphism
    \[
      \AutOf{
        \TheLocalField[\ThePlace][\sep]
        \fieldextends
        \TheGlobalField
      }
      \rightarrow
      \GalOf{
        \TheGlobalField[][\sep]
        \fieldextends
        \TheGlobalField
      }
    \]
    and the group
    \(
      \GalOf{
        \TheLocalField[\ThePlace][\sep]
        \fieldextends
        \TheLocalField[\ThePlace]
      }
    \)
    acts on $\ExtVectorSpace$ via this projection (here, we consider
    the identification made in~\ref{eq:identification}).
    Hence $\TheVectorSpace\subspace\TheVectorSpace[\ThePlace]$.

    In particular, these considerations apply to the situation
    discussed in Lemma~\ref{us:apartment}. In that case,
    $\ExtSplitTorus$ is just the maximal $\TheGlobalField$-torus
    $\TheTorus$ within the minimal parabolic $\TheParabolic$.
    $\TheTorus$ automatically splits over the separable closures
    $\TheGlobalField[][\sep]$ and $\TheLocalField[\ThePlace][\sep]$.
    \begin{prop}\label{us:metric}
      Assume that $\TheGroupScheme$ is a connected reductive
      group. There exists a euclidean metric on
      \(
        \TheEuBuild
        =
        \Product[\ThePlace\in\ThePlaceSet]{
          \TheEuBuild[\ThePlace]
        }
      \)
      and for each vertex $\TheRatVertex\in\VerticesOf{\TheRatBuild}$,
      there is
      a positive coefficient $\TheSpeed[\TheType]\in\TheReals$,
      depending only on the type $\TheType= \TypeOf{\TheRatVertex}$,
      such that the following hold:
      \begin{enumerate}
        \item\label{metric:one}
          For each vertex $\TheRatVertex\in\VerticesOf{\TheRatBuild}$,
          the rescaled function
          \(
            \TheSpeed[\TheType]\AlmostBusemann[\TheRatVertex]\mapcolon
            \TheSpecVertexSet\rightarrow\TheReals
          \)
          is the restriction of a Busemann function
          \(
            \TheBusemann[\TheRatVertex]\mapcolon
            \TheEuBuild\rightarrow\TheReals
          \)
          to $\TheSpecVertexSet$. Let
          \(
            \TheHoroCenter[\TheRatVertex]\in
            \BoundaryOf{\TheEuBuild}
          \)
          be the \ignore{horo}center of the Busemann function, i.e., the
          visual end point of the gradient
          of $\TheBusemann[\TheRatVertex]$.
        \item\label{metric:two}
          For each vertex $\TheRatVertex\in\VerticesOf{\TheRatBuild}$, the
          Busemann function $\TheBusemann[\TheRatVertex]$ is
          non-constant on each factor $\TheEuBuild[\ThePlace]$
          of $\TheEuBuild$. In particular if all factors
          $\TheEuBuild[\ThePlace]$ are irreducible (e.g., if
          $\TheGroupScheme$ is absolutely almost simple), the Busemann
          functions $\TheBusemann[\TheRatVertex]$ are in general
          position.
        \item\label{metric:three}
          The map $\TheRatVertex\mapsto\TheHoroCenter[\TheRatVertex]$
          induces an isometric embedding of $\TheRatBuild$
          into $\BoundaryOf{\TheEuBuild}$.
        \item\label{metric:four}
          For each $\TheLatticeElement\in\TheLattice$, each vertex
          $\TheRatVertex\in\VerticesOf{\TheRatBuild}$, and each point
          $\TheEuPoint\in\TheEuBuild$, we have
          \(
            \TheBusemannOf[\TheLatticeElement\TheRatVertex]{
              \TheLatticeElement\TheEuPoint
            }
            =
            \TheBusemannOf[\TheRatVertex]{
              \TheEuPoint
            }
          \).
          In particular, the map
          \(
            \TheRatVertex\mapsto\TheHoroCenter[\TheRatVertex]
          \)
          and the induced embedding
          $\TheRatBuild\monorightarrow\BoundaryOf{\TheEuBuild}$
          are $\TheLattice$-invariant.
      \end{enumerate}
    \end{prop}
    \begin{proof}
      If $\TheGroupScheme$ is anisotropic over $\TheGlobalField$,
      the building $\TheRatBuild$ is empty, and there is nothing to
      prove: the proposition is vacuously true. So, we assume
      that $\TheGroupScheme$ is $\TheGlobalField$-isotropic.

      We choose $\TheTorus$ as our standard apartment $\StdEuApp$ for
      $\TheEuBuild$, i.e., within each factor $\TheEuBuild[\ThePlace]$
      the standard apartment $\StdEuApp[\ThePlace]$ corresponds to
      the maximal $\TheLocalField[\ThePlace]$-split torus
      $\TheSplitTorus[\ThePlace]$ within $\TheTorus$. We fixed
      an inner product on $\ExtVectorSpace$ and this inner product
      induces inner products on each $\TheVectorSpace[\ThePlace]$,
      which is the metric model for the standard
      apartment $\StdEuApp[\ThePlace]$. This way, we define a metric
      on $\StdEuApp$ and thus on the euclidean building $\TheEuBuild$.

      By Lemma~\ref{us:apartment}, the
      function $\AlmostBusemann[\TheRatVertex]$ agrees with an
      affine function on $\StdEuApp$. Affine functions on euclidean spaces
      are Busemann functions up to rescaling. Hence, we can choose a
      factor $\TheSpeed[\TheType]$ so that the rescaled function
      $\TheSpeed[\TheType]\AlmostBusemann[\TheRatVertex]$ agrees
      with a Busemann function $\TheBusemann[\TheRatVertex]\mapcolon
      \TheEuBuild\rightarrow\TheReals$ on $\StdEuApp$.
      By Equation~\ref{eq:busemann.on.apartment} from the proof of
      Lemma~\ref{us:apartment}, the rescale factor
      $\TheSpeed[\TheType]$ depends only on the length of the
      canonical character $\TheHarderChar[{\TheParabolic[\TheRatVertex]}]$
      in the product
      \(
        \Product[\ThePlace\in\ThePlaceSet]{\TheVectorSpace[\ThePlace]}
      \).
      Hence, $\TheSpeed[\TheType]$ depends only on the type of
      $\TheRatVertex$.

      Again by Equation~\ref{eq:busemann.on.apartment}, the
      \ignore{horo}center
      of $\TheBusemann[\TheRatVertex]$ is given be the visual end of
      the fundamental weight. In particular, it is stabilized by
      $\TheParabolicOf{\TheAdeles[\ThePlaceSet]}$.
      As $\TheGroupSchemeOf{\TheLocalField[\ThePlace]}$ acts strongly
      transitively on the factor $\TheEuBuild[\ThePlace]$, translates
      of $\StdEuApp[\ThePlace]$ under the action of the parabolic
      $\TheParabolicOf{\TheLocalField[\ThePlace]}$
      cover $\TheEuBuild[\ThePlace]$. Hence, the
      \(
        \TheParabolicOf{\TheAdeles[\ThePlaceSet]}
      \)-translates
      of $\StdEuApp$ cover $\TheEuBuild$. The
      Transformation~Formula~\ref{eq:transformation.formula} implies
      that $\TheSpeed[\TheType]\AlmostBusemann[\TheRatVertex]$
      agrees with $\TheBusemann[\TheRatVertex]$ on each of these
      translates. This proves the first claim.

      The second claim follows directly from
      Equation~\ref{eq:busemann.on.apartment} since
      $\TheLocalDeg[\ThePlace]\neq\Zero$ for
      each $\ThePlace\in\ThePlaceSet$.

      Claim~\ref{metric:three} follows from the discussion preceding
      this proposition. The geometry of the chamber $\TheRatChamber$
      of $\TheRatBuild$ corresponding to the minimal parabolic
      $\TheParabolic$ (i.e., the angular distances between its vertices)
      is given by the angles between the fundamental weights in
      $\TheVectorSpace$. As
      $\TheVectorSpace\subspace\TheVectorSpace[\ThePlace]$, we
      have an induced diagonal embedding of $\TheVectorSpace$ into
      the orthogonal product
      $\Product[\ThePlace\in\ThePlaceSet]{\TheVectorSpace[\ThePlace]}$.
      This is an isometric embedding modeling the map
      $\TheRatVertex\mapsto\TheHoroCenter[\TheRatVertex]$ on the
      standard apartment. Considering other minimal
      $\TheGlobalField$-parabolics
      in the same apartment, we see that
      $\TheRatBuild\monorightarrow\BoundaryOf{\TheEuBuild}$ is
      an isometric embedding on the standard apartment. The
      choice of the standard apartment was arbitrary and does
      not influence the embedding. Hence,
      $\TheRatBuild\monorightarrow\BoundaryOf{\TheEuBuild}$
      is an isometric embedding on each apartment. Since any
      two points of $\TheRatBuild$ are contained in a common
      apartment,
      $\TheRatBuild\monorightarrow\BoundaryOf{\TheEuBuild}$
      preserves distances.

      Finally, claim~\ref{metric:four} follows from
      Equation~\ref{eq:invariance}. As $\TheLattice$ might not
      act type-preserving, we have to consider the rescaling factors
      $\TheSpeed[\TheType]$
      and $\TheSpeed[\TheLatticeElement\TheRatVertex]$. These factors
      depend only on the length of the associated canonical characters.
      As $\TheLattice$ acts by isometries,
      \(
        \TheSpeed[\TheType]=
        \TheSpeed[\TheLatticeElement\TheRatVertex]
      \).
    \end{proof}

    Let $\TheRatChamber$ be a chamber of $\TheRatBuild$ and let
    $\TheParabolic$ denote the corresponding minimal
    $\TheGlobalField$-parabolic. Let $\TheEuApp$ be an apartment
    of $\TheEuBuild$ with
    $\TheRatChamber\subseteq\BoundaryOf{\TheEuApp}$.
    The Busemann functions $\TheBusemann[\TheRatVertex]$ associated
    with vertices $\TheRatVertex\in\TheRatChamber$ restrict to
    affine functions on $\TheEuApp$. Put:
    \[
      \TheEuApp[\TheOrigin]
      :=
      \SetOf[
        \TheEuPoint\in\TheEuApp
      ]{
        \TheBusemannOf[\TheRatVertex]{\TheEuPoint}
        =
        \Zero
        \text{\ for each\ }
        \TheRatVertex\in\TheRatChamber
      }
    \]
    The metric on $\TheEuApp$ is constructed from a metric on
    $\CharactersOf[\TheGlobalField]{\TheParabolic}\tensor\TheReals$
    invariant under the action of the Weyl group; and the quotient
    $\TheEuApp\rmod\TheEuApp[\TheOrigin]$ is isometric to
    $\CharactersOf[\TheGlobalField]{\TheParabolic}\tensor\TheReals$.
    The Busemann functions $\TheBusemann[\TheRatVertex]$ descend
    to the quotient $\TheEuApp\rmod\TheEuApp[\TheOrigin]$. They
    form a system of coordinates, which under the isometry
    \(
      \TheEuApp\rmod\TheEuApp[\TheOrigin]
      \isometric
      \CharactersOf[\TheGlobalField]{\TheParabolic}\tensor\TheReals
    \)
    corresponds to the set of fundamental weights up to rescaling.
    Recall that the simple roots are related to the rescaled fundamental
    weights by the matrix $(\TheBaseChange[\TheCotype\TheType])$.
    Thus, we define:
    \[
      \DualBusemann[\TheCotype][\TheRatChamber] :=
      \Sum[\TheRatVertex\in\TheRatChamber]{
        \frac{\TheBaseChange[\TheCotype\TypeOf{\TheRatVertex}]}
        {\TheSpeed[\TypeOf{\TheRatVertex}]}
        \TheBusemann[\TheRatVertex]
      }
      \mapcolon
      \TheEuBuild\longrightarrow\TheReals
    \]
    \begin{observation}\label{us:mu.analogue}
      For any $\TheSpecVertex\in\TheSpecVertexSet\subseteq\TheEuBuild$,
      we have
      \(
        \DualBusemannOf[\TheCotype][\TheRatChamber]{\TheSpecVertex}
        =
        \LogNuOf[\TheCotype]{\TheParabolic,\StabOf{\TheSpecVertex}}
      \).\qed
    \end{observation}
    Hence Theorem~\ref{harder:one}, the first main theorem
    of Harder's reduction theory, implies:
    \begin{cor}\label{us:cor.to.thm.one}
      If $\TheGroupScheme$ is $\TheGlobalField$-isotropic, there
      is a constant $\HarderConstOne\in\TheReals$ such that for any
      point $\TheEuPoint\in\TheEuBuild$, there exists a chamber
      $\TheRatChamber$ in $\TheRatBuild$ with
      \(
        \DualBusemannOf[\TheCotype][\TheRatChamber]{\TheEuPoint}
        \geq
        \HarderConstOne
      \).\qed
    \end{cor}
    Restricted to $\TheEuApp$, the
    functions $\DualBusemann[\TheCotype][\TheRatChamber]$ are
    affine and the duality between fundamental weights and simple
    roots translates into the following relationship:
    \begin{observation}\label{us:duality}
      For any real number $\TheBound\in\TheReals$ and any face
      $\TheRatSimplex\subseteq\TheRatChamber$, we consider the
      convex cone (with tip parallel to $\TheEuApp[\TheOrigin]$):
      \[
        \ConvexChunkOf[\TheEuApp,\TheRatSimplex]{\TheBound}
        :=
        \SetOf[
          \TheEuPoint\in\TheEuApp
        ]{
          \TheBusemannOf[\TheRatVertex]{\TheEuPoint}
          \leq\TheBound
          \text{\ for each\ }
          \TheRatVertex\in\TheRatSimplex
        }
      \]
      We also define:
      \[
        \DualChunkOf[\TheEuApp,\TheRatSimplex]{\TheBound}
        :=
        \SetOf[
          \TheEuPoint\in\TheEuApp
        ]{
          \DualBusemannOf[\TypeOf{\TheRatVertex}]{\TheEuPoint}
          \geq\TheBound
          \text{\ for each\ }
          \TheRatVertex\in\TheRatSimplex
        }
      \]
      Then, $\DualChunkOf[\TheEuApp,\TheRatSimplex]{\Zero}$
      is the \notion{normal cone} for
      $\ConvexChunkOf[\TheEuApp,\TheRatSimplex]{\Zero}$; i.e.,
      $\DualChunkOf[\TheEuApp,\TheRatSimplex]{\Zero}$ consists
      precisely of
      those points in $\TheEuApp$ whose closest point projection onto
      $\ConvexChunkOf[\TheEuApp,\TheRatSimplex]{\Zero}$ lies in the
      tip $\TheEuApp[\TheOrigin]$.\qed
    \end{observation}
    Let
    \[
      \ClosestPoint[\TheEuApp,\TheRatSimplex][\TheBound]
      \mapcolon\TheEuApp\longrightarrow
      \ConvexChunkOf[\TheEuApp,\TheRatSimplex]{\TheBound}
    \]
    denote the closest point projection. As seen in
    Observation~\ref{behr:eliminate.apartment}, for
    $\TheEuPoint\in\TheEuApp$, the value
    \(
      \TheBaseValueOf[\TheRatSimplex,\TheRatVertex][\TheBound]{
        \TheEuPoint
      }
      :=
      \TheBusemannOf[\TheRatVertex]{
        \ClosestPointOf[\TheEuApp,\TheRatSimplex][\TheBound]{
          \TheEuPoint
        }
      }
    \)
    is independent of the apartment $\TheEuApp$. Recall the
    definition
    \[
      \CloseSimplexOf[\TheBound]{\TheEuPoint,\TheRatSimplex}
      :=
      \SetOf[
        \TheRatVertex\in\TheRatSimplex
      ]{
        \TheBaseValueOf[\TheRatSimplex,\TheRatVertex][\TheBound]{
          \TheEuPoint
        }=\TheBound
      }
      \text{.}
    \]
    Also recall that $\TheEuPoint$ is $\TheBound$-reduced by
    $\TheRatChamber$ if
    \(
      \CloseSimplexOf[\TheBound]{\TheEuPoint,\TheRatChamber}
      =
      \TheRatChamber
    \). The closest point projection is embedded into this
    terminology so that it allows us to characterize normal
    cones to
    \(
      \ConvexChunkOf[\TheEuApp,\TheRatChamber]{\TheBound}
    \).
    \begin{observation}\label{us:normal.cone}
      The set
      \[
        \NormalChunkOf[\TheEuApp,\TheRatChamber]{\TheBound}
        :=
        \SetOf[\TheEuPoint\in\TheEuApp]{
          \TheEuPoint\text{\ is\ }
          \TheBound\text{-reduced by\ }
          \TheRatChamber
        }
      \]
      is the normal cone to
      \(
        \ConvexChunkOf[\TheEuApp,\TheRatChamber]{\TheBound}
      \).
      In particular, it is a translate of
      $\DualChunkOf[\TheEuApp,\TheRatChamber]{\Zero}$.
      Thus, there exist real constants $\AltBound[\One],\ldots,
      \AltBound[\TheLastRoot]$ only depending on $\TheBound$ such
      that
      \[
        \SetOf[\TheEuPoint\in\TheEuApp]{
          \TheEuPoint\text{\ is\ }
          \TheBound\text{-reduced by\ }
          \TheRatChamber
        }
        =
        \SetOf[
          \TheEuPoint\in\TheEuApp
        ]{
          \DualBusemannOf[\TheCotype][\TheRatChamber]{\TheEuPoint}
          \geq\AltBound[\TheCotype]
          \text{\ for each\ }\TheCotype=\One,\ldots,\TheLastRoot
        }
        \text{.}\qed
      \]
    \end{observation}
    \begin{figure}[t]
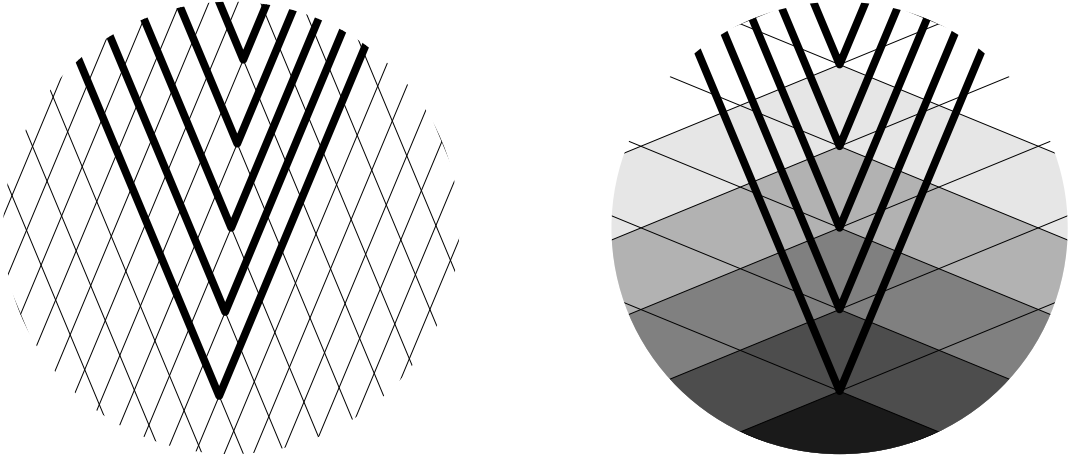

      \null\hfill
      \setlength{\unitlength}{1cm}
      \begin{picture}(14,6)
        \put(0,0){\includegraphics{mu_cone.ps}}
        \put(8,0){\includegraphics{normal_cone.ps}}
      \end{picture}
      \hfill\null
      \caption{%
        normal cones%
      \label{us:cones.fig}}
      {\small
      The family
      \(
        \DualChunkOf[\TheEuApp,\TheRatSimplex]{
          \DummyArg
        }
      \),
      shown on the left,
      is defined via the functions
      \(
        \DualBusemann[\TheType]
      \).
      Some level sets of the $\DualBusemann[\TheType]$
      are drawn. In contrast,
      the family
      \(
        \NormalChunkOf[\TheEuApp,\TheRatChamber]{
          \DummyArg
        }
      \),
      shown on the right, is defined in terms of
      the Busemann functions
      \(
        \TheBusemann[\TheRatVertex]
      \).
      The latter are normalized to have unit length
      gradient with respect to the metric. Their
      level sets are shown. The shaded areas
      indicate the family
      \(
        \ConvexChunkOf[\TheEuApp,\TheRatSimplex]{
          \DummyArg
        }
      \).}
    \end{figure}
    If $\AltEuApp$ is another apartment of $\TheEuBuild$ whose
    visual boundary contains the chamber $\TheRatChamber$, then
    the isomorphism of Coxeter complexes
    \(
      \TheIsometry\mapcolon\TheEuApp\rightarrow\AltEuApp
    \)
    from Observation~\ref{behr:eliminate.apartment} does not
    only commute with the Busemann functions
    $\TheBusemann[\TheRatVertex]$ but also with the functions
    $\DualBusemann[\TheCotype][\TheRatChamber]$. In particular,
    it identifies
    \(
      \DualChunkOf[\TheEuApp,\TheRatChamber]{\TheBound}
    \)
    with
    \(
      \DualChunkOf[\AltEuApp,\TheRatChamber]{\TheBound}
    \)
    and
    \(
      \NormalChunkOf[\TheEuApp,\TheRatChamber]{\TheBound}
    \)
    with
    \(
      \NormalChunkOf[\AltEuApp,\TheRatChamber]{\TheBound}
    \).
    We define:
    \begin{eqnarray*}
      \NormalChunkOf[\TheRatChamber]{\TheBound}
      &:=&
      \Union[\TheEuApp \,\,:\,\, \TheRatChamber\subseteq\BoundaryOf{\TheEuApp}]{
        \NormalChunkOf[\TheEuApp,\TheRatChamber]{\TheBound}
      }\\
      &=&
      \SetOf[\TheEuPoint\in\TheEuBuild]{
        \TheEuPoint \text{\ is\ }
        \TheBound\text{-reduced by\ }
        \TheRatChamber
      }\\
      \DualChunkOf[\TheRatChamber]{\AltBound}
      &:=&
      \Union[\TheEuApp \,\,:\,\, \TheRatChamber\subseteq\BoundaryOf{\TheEuApp}]{
        \DualChunkOf[\TheEuApp,\TheRatChamber]{\AltBound}
      }
    \end{eqnarray*}

    The systems $\NormalChunkOf[\TheEuApp,\TheRatChamber]{\TheBound}$
    and $\DualChunkOf[\TheEuApp,\TheRatChamber]{\AltBound}$ are
    strongly related:
    \begin{observation}\label{us:cofinal}
      For any $\TheBound$ there exist $\AltBound[+]$ and
      $\AltBound[-]$ such that
      \[
        \DualChunkOf[\TheEuApp,\TheRatChamber]{\AltBound[-]}
        \subseteq
        \NormalChunkOf[\TheEuApp,\TheRatChamber]{\TheBound}
        \subseteq
        \DualChunkOf[\TheEuApp,\TheRatChamber]{\AltBound[+]}
      \]
      for any apartment $\TheEuApp$ and any $\TheGlobalField$-rational
      chamber $\TheRatChamber$ in the visual boundary of $\TheEuApp$.

      Analogously, for any $\AltBound$ there exist $\TheBound[+]$
      and $\TheBound[-]$ such that
      \[
        \NormalChunkOf[\TheEuApp,\TheRatChamber]{\TheBound[-]}
        \subseteq
        \DualChunkOf[\TheEuApp,\TheRatChamber]{\AltBound}
        \subseteq
        \NormalChunkOf[\TheEuApp,\TheRatChamber]{\TheBound[+]}
      \]
      for any $\TheRatChamber\in\TheRatBuild$ and any euclidean
      apartment $\TheEuApp$ whose visual boundary
      contains $\TheRatChamber$.\qed
    \end{observation}
    We can also relate the
    subsets $\NormalChunkOf[\TheRatChamber]{\TheBound}$
    to their Hausdorff neighborhoods. For any subset
    $\TheEuSubset\subseteq\TheEuBuild$, let
    \(
      \MetricNbhdOf[\TheHausdorffDist]{\TheEuSubset}
    \)
    denote the Hausdorff neighborhood of radius $\TheHausdorffDist$
    around $\TheEuSubset$, i.e., the set of points in $\TheEuBuild$ of
    distance at most $\TheHausdorffDist$ to $\TheEuSubset$.
    \begin{observation}\label{us:hausdorff.distance}
      For any fixed distance $\TheHausdorffDist\geq\Zero$, there is a
      constant $\TheShiftConst$ such that
      \[
        \MetricNbhdOf[\TheHausdorffDist]{
          \NormalChunkOf[\TheRatChamber]{\TheBound}
        }
        \subseteq
        \NormalChunkOf[\TheRatChamber]{\TheBound-\TheShiftConst}
      \]
      for any $\TheBound\in\TheReals$ and any chamber
      $\TheRatChamber$ in $\TheRatBuild$.

      Conversely, for any given $\TheShiftConst$ there exists a
      constant $\TheHausdorffDist$ such that
      \[
        \NormalChunkOf[\TheRatChamber]{\TheBound-\TheShiftConst}
        \subseteq
        \MetricNbhdOf[\TheHausdorffDist]{
          \NormalChunkOf[\TheRatChamber]{\TheBound}
        }
      \]
      for any $\TheBound\in\TheReals$ and any chamber
      $\TheRatChamber$ in $\TheRatBuild$.\qed
    \end{observation}
    We are ready for the geometric version of Theorem~\ref{harder:one}.
    \begin{prop}\label{us:one}
      Assume that $\TheGroupScheme$ is $\TheGlobalField$-isotropic.
      For any fixed diameter $\TheDiam\in\TheReals$ there exists
      a constant $\TheLowerBound\in\TheReals$ such that for any
      $\TheEuPoint\in\TheEuBuild$ there is a chamber
      $\TheRatChamber$ that $\TheLowerBound$-reduces each
      point $\AltEuPoint$ of distance at most $\TheDiam$
      to $\TheEuPoint$.
    \end{prop}
    \begin{proof}
      By Corollary~\ref{us:cor.to.thm.one}, there is a constant
      $\HarderConstOne$ such that
      \[
        \TheEuBuild = \Union[\TheRatChamber\in\ChambersOf{\TheRatBuild}]{
          \DualChunkOf[\TheRatChamber]{\HarderConstOne}
        }
        \text{.}
      \]
      By Observation~\ref{us:cofinal}, there is a bound
      $\AltLowerBound$ such that
      \[
        \TheEuBuild = \Union[\TheRatChamber\in\ChambersOf{\TheRatBuild}]{
          \NormalChunkOf[\TheRatChamber]{\AltLowerBound}
        }
        \text{.}
      \]
      Now, one chooses $\TheLowerBound$ so that
      \(
        \NormalChunkOf[\TheRatChamber]{\TheLowerBound}
      \)
      contains the $\TheDiam$-Hausdorff neighborhood of
      \(
        \NormalChunkOf[\TheRatChamber]{\AltLowerBound}
      \) for any chamber $\TheRatChamber$ in $\TheRatBuild$.
    \end{proof}
    The main theorem of geometric reduction theory reads as
    follows:
    \begin{theorem}\label{us:main}
      Assume that $\TheGroupScheme$ is connected, reductive,
      and defined and isotropic over the global function
      field $\TheGlobalField$.
      For any diameter $\TheDiam$ there exist constants
      $\TheLowerBound$ and $\TheUpperBound$
      such that
      \(
        \FamOf[\TheRatVertex\in\VerticesOf{\TheRatBuild}]{
          \TheBusemann[\TheRatVertex]\mapcolon
          \TheEuBuild\rightarrow\TheReals
        }
      \)
      together with the constants $\TheLowerBound$ and $\TheUpperBound$
      is a $\TheDiam$-uniform and $\TheLattice$-invariant reduction
      datum.
      Moreover, for any $\TheEuPoint\in\TheEuBuild$ and any chamber
      $\TheRatChamber$ in $\TheRatBuild$ that
      $\TheLowerBound$-reduces $\TheEuPoint$, the simplex
      $\CloseSimplexOf[\TheUpperBound]{\TheEuPoint,\TheRatChamber}$
      is contained in any chamber $\AltRatChamber$ that
      $\TheLowerBound$-reduces $\TheEuPoint$.
    \end{theorem}
    \begin{proof}
      Using Observation~\ref{us:cofinal}, choose $\HarderConstOne$ so that
      $\NormalChunkOf[\TheRatChamber]{\TheLowerBound}\subseteq
      \DualChunkOf[\TheRatChamber]{\HarderConstOne}$ for
      any $\TheRatChamber$ in $\TheRatBuild$.
      Then, $\HarderConstOne$ is a lower reduction bound.

      We can be a little more specific: Let
      $\TheParabolic$ be the minimal $\TheGlobalField$-parabolic
      corresponding to a chamber $\TheRatChamber$ that
      $\TheLowerBound$-reduces the point
      $\TheGroupElement\TheBasePolyVertex$ for some
      $\TheAdeleGroupElement\in\TheGroupSchemeOf{\TheAdeles[\ThePlaceSet]}$.
      Then
      \(
        \LogNuOf[\TheCotype]{
          \TheParabolic,
          \Conj{\TheAdeleGroupElement}{\TheGroupSchemeOf{\TheOka}}
        }
        \geq
        \HarderConstOne
      \)
      for all $\TheCotype$.

      By Theorem~\ref{harder:two}, there
      is a corresponding $\HarderConstTwo$. Using again
      Observation~\ref{us:cofinal}, we find a constant
      $\TheUpperBound$ such that
      \(
        \NormalChunkOf[\TheRatChamber]{\TheUpperBound}\subseteq
        \DualChunkOf[\TheRatChamber]{\HarderConstTwo}
      \)
      for all chamber $\TheRatChamber$ in $\TheRatBuild$.
      Then, for any vertex $\TheRatVertex\in\TheRatChamber$ and any
      $\TheEuPoint\in\TheEuBuild$ that is $\TheLowerBound$-reduced
      by $\TheRatChamber$, we have:
      \[
        \TheRatVertex\in\CloseSimplexOf[\TheUpperBound]{
          \TheEuPoint,\TheRatChamber
        }
        \quad\Longrightarrow\quad
        \DualBusemannOf[\TypeOf{\TheRatVertex}][\TheRatChamber]{
          \TheEuPoint
        }
        \geq
        \HarderConstTwo
      \]
      For $\TheEuPoint=\TheAdeleGroupElement\TheBasePolyVertex$,
      it follows that the maximal $\TheGlobalField$-parabolic
      $\TheParabolic[\TheRatVertex]$ corresponding to $\TheRatVertex$
      contains any minimal $\TheGlobalField$-parabolic
      $\AltParabolic$ whose chamber $\AltRatChamber$
      $\TheLowerBound$-reduces $\TheAdeleGroupElement\TheBasePolyVertex$.
      Hence,
      \(
        \CloseSimplexOf[\TheUpperBound]{
          \TheAdeleGroupElement\TheBasePolyVertex,
          \TheRatChamber
        }
        \subseteq
        \AltRatChamber
      \).

      Extending coverage from $\TheSpecVertexSet$ to all of
      $\TheEuBuild$ requires changing the constants $\TheLowerBound$
      and $\TheUpperBound$ only by a little.
      Hence, we have established a reduction datum.

      We have already argued in Proposition~\ref{us:one} that this
      reduction datum is $\TheDiam$-uniform. That it is
      $\TheLattice$-invariant follows from Proposition~\ref{us:metric},
      part~\ref{metric:four}.
    \end{proof}

    It remains to discuss $\TheLattice$-cocompactness of the
    reduction datum. We consider the filtration of $\TheEuBuild$
    by subspaces
    \[
      \TheRedFiltr[\TheBound]
      :=
      \SetOf[
        \TheEuPoint\in\TheEuBuild
      ]{
        \TheBusemannOf[\TheRatVertex]{\TheEuPoint}\leq\TheBound
        \text{\ for all\ }
        \TheRatChamber
        \text{\ reducing\ }
        \TheEuPoint
        \text{\ and all\ }
        \TheRatVertex\in\TheRatChamber
      }
      \text{.}
    \]
    \begin{theorem}\label{us:mahler}
      If $\TheGroupScheme$ does not admit any non-trivial
      $\TheGlobalField$-characters,
      then $\TheRedFiltr[\TheBound]$ has compact
      quotient modulo $\TheLattice$. If there is a non-trivial
      $\TheGlobalField$-character
      $\TheGroupScheme\rightarrow\Mult$, then
      $\TheRedFiltr[\TheBound]$ does not have a compact
      quotient modulo $\TheLattice$ unless
      $\TheRedFiltr[\TheBound]$ is empty.
    \end{theorem}
    \begin{proof}
      If there are no non-trivial characters, we have
      $\TheGroupSchemeOf{\TheAdeles}=
      \CharKernel{\TheGroupSchemeOf{\TheAdeles}}$.
      Hence, cocompactness of the $\TheRedFiltr[\TheBound]$ is
      immediate from Mahler's Compactness Criterion in its
      alternate form~\ref{harder:alternate}.

      If there is a non-trivial character on $\TheGroupScheme$,
      then $\TheGroupScheme$ has a central $\TheGlobalField$-torus.
      Corresponding to this torus, the euclidean building $\TheEuBuild$
      has a euclidean space as a factor. Dirichlet's unit theorem implies
      that the $\ThePlaceSet$-arithmetic subgroup does not act cocompactly
      in the direction of this factor.
    \end{proof}
    In particular, an absolutely almost simple non-commutative
    group $\TheGroupScheme$ does not admit non-trivial
    $\TheGlobalField$-characters. Hence, Theorem~\ref{behr:main}
    follows from Theorem~\ref{us:main} and Theorem~\ref{us:mahler}.
    \begin{rem}
      Formally, the filtration is meaningful even for
      anisotropic $\TheGroupScheme$. In that case,
      $\TheRedFiltr[\TheBound]=\TheEuBuild$ independent of
      $\TheBound$. Theorem~\ref{harder:cocompact.case}
      implies that $\TheEuBuild\rmod\TheLattice$ is cocompact
      in this case. This way, one recovers Serre's proof
      that $\TheGroupSchemeOf{\TheOkaRing}$ is of type
      \FType[\infty] for $\TheGlobalField$-anisotropic
      $\TheGroupScheme$ \cite[Cas~(b), p.~126--127]{Serre:1971}.
    \end{rem}
    \begin{rem}\label{us:alternative.filtration}
      If $\TheGroupScheme$ is $\TheGlobalField$-isotropic and
      non-commutative, then $\TheRatBuild$ is non-empty and
      Corollary~\ref{behr:low}
      yields an alternative description of the filtration:
      \[
        \TheRedFiltr[\TheBound]
        =
        \SetOf[
          \TheEuPoint\in\TheEuBuild
        ]{
          \TheBusemannOf[\TheRatVertex]{\TheEuPoint}\leq\TheBound
          \text{\ for some\ }
          \TheRatChamber
          \text{\ reducing\ }
          \TheEuPoint
          \text{\ and all\ }
          \TheRatVertex\in\TheRatChamber
        }
        \qquad\text{for\ }
        \TheBound\geq\TheUpperBound\qed
      \]
    \end{rem}

  \section{The structure of the quotient \boldmath$\TheEuBuild\rmod\TheLattice$}%
  \label{sec:quotient}
    Let $\TheGroupScheme$ be a connected reductive
    $\TheGlobalField$-group.
    To simplify the exposition, we assume in this section
    that $\TheGroupScheme$
    does not admit a non-trivial $\TheGlobalField$-character so that
    the filtration $\TheRedFiltr[\TheBound]$ is $\TheLattice$-cocompact
    by Theorem~\ref{us:mahler}. Our goal is to show that in this
    case the image $\TheImage$ of $\TheLattice$ in the
    automorphism group
    $\AutOf{\TheEuBuild}$ of the euclidean building is a lattice
    in $\AutOf{\TheEuBuild}$.

      It is well-known that $\TheImage$ is a lattice.
      For a Chevalley group scheme $\TheGroupScheme$, Harder
      \cite[page~41]{Harder:1969}
      constructs a fundamental set for $\TheGroupSchemeOf{\TheGlobalField}$
      in $\TheGroupSchemeOf{\TheAdeles}$. The restriction to Chevalley
      group schemes is unnecessary. In the same paper, Harder generalizes
      the main statements of reduction theory (also first proved for
      Chevalley groups) to arbitrary reductive groups. His construction
      of a fundamental domain can then be carried out in the same vein --
      in fact, Harder points out this
      possibility \cite[page~51]{Harder:1969}.
      He also remarks that the same argument as in
      \cite{Godement:1963} following Th\'{e}or\`{e}m~7 then shows that
      $\TheGroupSchemeOf{\TheGlobalField}$ is a lattice in
      $\TheGroupSchemeOf{\TheAdeles}$. Harder also indicates
      \cite[pages~51ff]{Harder:1969} how this argument can be adapted
      to deal with the lattice $\TheGroupSchemeOf{\TheOkaRing}$ in
      $\TheGroupSchemeOf{\TheAdeles[\ThePlaceSet]}$. Since
      $\TheGroupSchemeOf{\TheAdeles[\ThePlaceSet]}$ acts cocompactly
      on $\TheEuBuild$, the result shows that
      $\TheImage$ is a lattice in
      $\AutOf{\TheEuBuild}$.

      Hence, the point of this section is to demonstrate that
      Theorem~\ref{us:main} preserves the necessary information.
      In fact, even the rough strategy of the proof is the same:
      in Proposition~\ref{quotient:fundamental.domain} a fundamental
      set for $\TheImage$ in $\TheEuBuild$ is constructed;
      an application of Serre's criterion in the proof of
      Proposition~\ref{quotient:lattice} replaces the
      covolume estimate following
      \cite[Th\'{e}or\`{e}m~7]{Godement:1963}.

    \begin{rem}\label{quotient:covolume}
      The exact covolume of $\TheGroupSchemeOf{\TheGlobalField}$
      in $\TheGroupSchemeOf{\TheAdeles}$ is not known in general.
      See \cite{Behrend.Dhillon:2009} for conjectural values and
      partial results.
    \end{rem}

    \begin{observation}\label{quotient:finite.kernel}
      The kernel of the projection
      \(
        \TheLattice\rightarrow\TheImage
      \)
      is the kernel of the action of $\TheLattice$ on
      $\TheEuBuild$. Hence, it is finite as it is clearly contained
      in any vertex stabilizer of $\TheLattice$, which is finite.\qed
    \end{observation}

    \begin{observation}\label{quotient:anisotropic}
      If $\TheGroupScheme$ is $\TheGlobalField$-anisotropic,
      the group $\TheImage$ is a lattice. In fact, the
      quotient $\TheEuBuild\rmod\TheImage$ is compact by
      Theorem~\ref{us:mahler} or by
      Theorem~\ref{harder:cocompact.case}.\qed
    \end{observation}

    Let $\TheLowerBound$ and $\TheUpperBound$ be constants
    satisfying Theorem~\ref{us:main}.

    \begin{lemma}\label{quotient:close.to.boundary}
      For any point $\TheEuPoint\in\TheEuBuild$ there is at most
      one chamber $\TheRatChamber$ in $\TheRatBuild$ with
      \(
        \TheEuPoint\in\NormalChunkOf[\TheRatChamber]{\TheUpperBound}
      \).
    \end{lemma}
    \begin{proof}
      Assume $\TheEuPoint\in\NormalChunkOf[\TheRatChamber]{\TheUpperBound}$.
      Then $\CloseSimplexOf[\TheUpperBound]{\TheEuPoint,\TheRatChamber}=
      \TheRatChamber$. Hence, $\TheRatChamber$ is contained in any
      chamber $\AltRatChamber$ that $\TheLowerBound$-reduces
      $\TheEuPoint$. In particular,
      $\TheRatChamber$ is the only chamber that $\TheLowerBound$-reduces
      $\TheEuPoint$, whence it is the only chamber that
      $\TheUpperBound$-reduces $\TheEuPoint$.
    \end{proof}
    \begin{cor}\label{quotient:stabilizer}
      For any point
      \(
        \TheEuPoint\in
        \NormalChunkOf[\TheRatChamber]{\TheUpperBound}
      \),
      the set
      $\NormalChunkOf[\TheRatChamber]{\TheUpperBound}$
      is invariant under the induced action of
      the stabilizer $\StabOf[\TheImage]{\TheEuPoint}$.\qed
    \end{cor}

    The following theorem provides a fundamental set for the
    action of $\TheLattice$ on $\TheEuBuild$. One should
    compare it to \cite[Theorem~9, page~106]{Serre:1980}.
    \begin{prop}\label{quotient:fundamental.domain}
      Assume that $\TheGroupScheme$ is $\TheGlobalField$-isotropic and
      non-commutative.
      There exists a constant $\TheHausdorffDist$, finitely many
      points $\TheEuPoint[\One],\ldots,\TheEuPoint[\TheLastPoint]$,
      and as many chambers $\TheRatChamber[\One],\ldots,
      \TheRatChamber[\TheLastPoint]$ in $\TheRatBuild$ such that
      the following hold:
      \begin{enumerate}
        \item
          The point $\TheEuPoint[\TheIndex]$ is $\TheUpperBound$-reduced
          by $\TheRatChamber[\TheIndex]$ for each $\TheIndex\in
          \SetOf{\One,\ldots,\TheLastPoint}$. In particular, the
          union $\TheSector[\TheIndex]$ of rays from $\TheEuPoint[\TheIndex]$
          with visual endpoint in $\TheRatChamber[\TheIndex]$ is isometric
          to a flat sector.
        \item
          Every point in $\TheEuBuild$ is within distance
          $\TheHausdorffDist$ to the orbit of some sector
          $\TheSector[\TheIndex]$. Equivalently, the $\TheImage$-translates
          of
          \[
            \TheFundamentalDomain
            :=
            \MetricNbhdOf[\TheHausdorffDist]{
              \Union[\TheIndex=\One][\TheLastPoint]{
                \TheSector[\TheIndex]
              }
            }
          \]
          cover $\TheEuBuild$.
        \item
          For $\TheIndex\neq\AltIndex$, the $\TheImage$-orbits of
          $\TheSector[\TheIndex]$ and $\TheSector[\AltIndex]$
          are disjoint.
      \end{enumerate}
    \end{prop}
    \begin{figure}[t]
      \null\hfill
      \setlength{\unitlength}{1cm}
      \begin{picture}(4,6)
        \put(1,0){\includegraphics{sector.ps}}
        \put(0.6,2.6){$\TheEuPoint[\TheIndex]$}
        \put(2.6,3.4){$\TheSector[\TheIndex]$}
        \put(1.6,1){$\NormalChunkOf[\TheEuApp,\TheRatChamber]{\TheUpperBound}$}
        \put(-1.5,3){$\TheEuApp$}
        \put(3.3,4.2){$\TheRatChamber[\TheIndex]$}
      \end{picture}
      \hfill\null
      \caption{a wedge
      $\NormalChunkOf[\TheEuApp,\TheRatChamber]{\TheUpperBound}$
      and a sector $\TheSector[\TheIndex]$%
      \label{quotient:wedge.fig}}
      {\small
      Here, the rational building $\TheRatBuild$ has dimension $\One$
      and the euclidean building~$\TheEuBuild$ has dimension~$\Three$.}
    \end{figure}
    \begin{proof}
      The set
      \(
        \Union[\TheRatChamber\in\ChambersOf{\TheRatBuild}]{
          \NormalChunkOf[\TheRatChamber]{\TheUpperBound}
        }
      \)
      is $\TheImage$-invariant. Hence
      \(
        \TheTips
        :=
        \TheRedFiltr[\TheUpperBound]
        \intersect
        \Union[\TheRatChamber\in\ChambersOf{\TheRatBuild}]{
          \NormalChunkOf[\TheRatChamber]{\TheUpperBound}
        }
      \)
      is $\TheImage$-invariant with compact quotient.
      By Lemma~\ref{quotient:close.to.boundary}, each point
      $\TheEuPoint\in\TheTips$ has a unique $\TheUpperBound$-reducing
      chamber in $\TheRatBuild$. Moreover, this chamber varies
      $\TheImage$-equivariantly with the point. Since
      $\TheTips\rmod\TheImage$ is compact, there exist a constant
      $\TheHausdorffDist[\Zero]$, finitely
      many pairwise $\TheImage$-inequivalent
      chambers $\TheRatChamber[\One],\ldots,
      \TheRatChamber[\TheLastPoint]$ in $\TheRatBuild$, and
      points $\TheEuPoint[\One],\ldots,\TheEuPoint[\TheLastPoint]$
      such that the following holds:
      \begin{quote}
        For
        each point $\TheEuPoint\in\TheTips$ there exists a
        $\TheImageElement\in\TheImage$ and a unique index
        $\TheIndex\in\SetOf{\One,\ldots,\TheLastPoint}$ such that
        $\TheEuPoint$ is within distance at
        most $\TheHausdorffDist[\Zero]$ of
        $\TheImageElement\TheEuPoint[\TheIndex]$ and so that
        $\TheEuPoint$ is $\TheUpperBound$-reduced by the chamber
        $\TheImageElement\TheRatChamber[\TheIndex]$.
      \end{quote}
      Let $\TheSector[\TheIndex]$ be the union of geodesic rays
      from $\TheEuPoint[\TheIndex]$ with visual endpoint in
      $\TheRatChamber[\TheIndex]$. We put
      \[
        \TheFundamentalDomain[\Zero]
        :=
        \MetricNbhdOf[{\TheHausdorffDist[\Zero]}]{
          \Union[\TheIndex=\One][\TheLastPoint]{
            \TheSector[\TheIndex]
          }
        }
      \]
      and claim that the $\TheImage$-translates
      of $\TheFundamentalDomain[\Zero]$ cover the union
      \(
        \Union[\TheRatChamber\in\ChambersOf{\TheRatBuild}]{
          \NormalChunkOf[\TheRatChamber]{\TheUpperBound}
        }
      \).
      Recall that $\NormalChunkOf[\TheRatChamber]{\TheUpperBound}$ is
      the union of \emph{wedges}
      $\NormalChunkOf[\TheEuApp,\TheRatChamber]{\TheUpperBound}$. The
      tip of such a wedge consists precisely of the subspace
      \(
        \SetOf[\TheEuPoint\in\TheEuApp]{
          \TheBusemannOf[\TheRatVertex]{\TheEuPoint}=\TheUpperBound
          \text{\ for all\ }
          \TheRatVertex\in\TheRatChamber
        }
      \), which is a subset of $\TheTips$.
      Now, let
      \(
        \AltEuPoint\in
        \NormalChunkOf[\TheEuApp,\TheRatChamber]{\TheUpperBound}
      \).
      Then, there is a unique point $\TheEuPoint\in
      \TheTips\intersect\NormalChunkOf[\TheEuApp,\TheRatChamber]{\TheUpperBound}$
      such that $\AltEuPoint$ lies on a geodesic ray from $\TheEuPoint$
      with visual endpoint in $\TheRatChamber$. There is
      $\TheImageElement\in\TheImage$ and a subscript
      $\TheIndex\in\SetOf{\One,\ldots,\TheLastPoint}$ such that
      $\TheRatChamber=\TheImageElement\TheRatChamber[\TheIndex]$
      and such that
      $\TheEuPoint$ is within distance at most $\TheHausdorffDist$ of
      $\TheImageElement\TheEuPoint[\TheIndex]$. It follows that
      $\AltEuPoint$ is within distance $\TheHausdorffDist$ of
      $\TheImageElement\TheSector[\TheIndex]$.

      By Theorem~\ref{us:one},
      \(
        \TheEuBuild=
        \Union[\TheRatChamber]{
            \NormalChunkOf[\TheRatChamber]{\TheLowerBound}
        }
      \), and by Observation~\ref{us:hausdorff.distance} there is
      a constant $\TheHausdorffDist[\One]$ such that
      \(
        \NormalChunkOf[\TheRatChamber]{\TheLowerBound}
        \subseteq
        \MetricNbhd[{\TheHausdorffDist[\One]}]{
          \NormalChunkOf[\TheRatChamber]{\TheUpperBound}
        }
      \)
      for all chambers $\TheRatChamber$ of $\TheRatBuild$. It
      follows that
      $\TheHausdorffDist:=\TheHausdorffDist[\Zero]+
      \TheHausdorffDist[\One]$ yields a
      domain
      \[
        \TheFundamentalDomain
        :=
        \MetricNbhdOf[{\TheHausdorffDist}]{
          \Union[\TheIndex=\One][\TheLastPoint]{
            \TheSector[\TheIndex]
          }
        }
      \]
      whose $\TheImage$-translates cover $\TheEuBuild$.
    \end{proof}

    \begin{lemma}\label{quotient:stabilizer.index}
      For every distance $\TheHausdorffDist$ there is a
      constant $\TheIndexBound$ such that for any two
      points $\TheEuPoint,\AltEuPoint\in\TheEuBuild$
      of distance at most $\TheHausdorffDist$
      the inequality
      \[
        \frac{\One}{\TheIndexBound}
        \CardOf{\StabOf[\TheImage]{\TheEuPoint}}
        \leq
        \CardOf{\StabOf[\TheImage]{\AltEuPoint}}
        \leq
        \TheIndexBound
        \CardOf{\StabOf[\TheImage]{\TheEuPoint}}
      \]
      holds.
    \end{lemma}
    \begin{proof}
      The ball $\TheBall := \TheBallOf[\TheHausdorffDist]{\AltEuPoint}$
      of radius $\TheHausdorffDist$ with center $\AltEuPoint$ is
      invariant under the induced action of
      \(
        \StabOf[\TheImage]{\AltEuPoint}
      \).
      The group $\TheImage$ acts by cell-permuting homeomorphisms
      on $\TheEuBuild$, hence the induced action on the barycentric
      subdivison is rigid: a cell that is stabilized is fixed pointwise.
      By restriction, the action of
      \(
        \StabOf[\TheImage]{\AltEuPoint}
      \)
      on the induced cell decomposition of $\TheBall$ inherits this
      property. In particular, the size of any orbit is bounded by
      the number $\TheIndexBound[\AltEuPoint]$ of cells in $\TheBall$.

      The index of the group
      \(
        \StabOf[\TheImage]{\TheEuPoint}\intersect
        \StabOf[\TheImage]{\AltEuPoint}
      \)
      in
      \(
        \StabOf[\TheImage]{\AltEuPoint}
      \)
      is given by the size of the orbit
      \(
        \StabOf[\TheImage]{\AltEuPoint}\acts\TheEuPoint
      \).
      Hence, it is bounded by $\TheIndexBound[\AltEuPoint]$.

      The claim follows since there is a maximum number of cells
      that a ball of radius $\TheHausdorffDist$ centered anywhere
      in $\TheEuBuild$ can meet.
    \end{proof}

    \begin{prop}\label{quotient:lattice}
      The group $\TheImage$ is a lattice in $\AutOf{\TheEuBuild}$.
    \end{prop}
    \begin{proof}
      We consider the action of $\AutOf{\TheEuBuild}$ on
      the set $\ChambersOf{\TheEuBuild}$ of chambers of
      $\TheEuBuild$. By \cite[1.6~Corollary]{Bass.Lubotzky:2001},
      we have to argue that the infinite sum
      \begin{equation}\label{serre}
        \Sum[\TheEuChamber\in\TheChRepSet]{
          \frac{\One}{\CardOf{\StabOf[\TheImage]{\TheEuChamber}}}
        }
      \end{equation}
      converges, where $\TheChRepSet$ is a set of representatives
      of $\ChambersOf{\TheEuBuild}\rmod\TheImage$.

      Let $\TheHausdorffDist$,
      $\TheEuPoint[\One],\ldots,\TheEuPoint[\TheLastPoint]$,
      $\TheRatChamber[\One],\ldots,\TheRatChamber[\TheLastPoint]$,
      $\TheSector[\One],\ldots,\TheSector[\TheLastIndex]$, and
      $\TheFundamentalDomain$ be as
      in Proposition~\ref{quotient:fundamental.domain}.
      We now choose $\TheChRepSet$ to be the collection of
      all chambers in $\TheEuBuild$ that
      intersect $\TheFundamentalDomain$.

      The set
      \(
        \TheTips[\TheIndex]
        :=
        \TheRedFiltr[\TheUpperBound]
        \intersect
        \NormalChunkOf[{\TheRatChamber[\TheIndex]}]{\TheUpperBound}
      \)
      consists of the tips of those wedges that form
      $\NormalChunkOf[{\TheRatChamber[\TheIndex]}]{\TheUpperBound}$.
      Let $\TheRay$ be some geodesic ray from $\TheEuPoint[\TheIndex]$
      in $\TheSector[\TheIndex]$, i.e., the visual endpoint $\TheEnd$
      of $\TheRay$ lies in $\TheRatChamber[\TheIndex]$. Let
      $\TheTree[\TheRay]$ be the union of all geodesic rays in
      $\TheEuBuild$ that share an infinite segment with $\TheRay$.
      Then $\TheTree[\TheRay]$ is a locally finite tree and
      intersects $\TheTips[\TheIndex]$ in a discrete set.
      Let $\TheRepSet[\TheRay]$ be the set of points
      in $\TheTips[\TheIndex]\intersect\TheTree[\TheRay]$
      that lie within distance $\TheHausdorffDist$ of
      $\TheEuPoint[\TheIndex]$.

      Consider a point $\AltEuPoint$ on the ray $\TheRay$.
      The union of all geodesic rays in $\TheEuBuild$ with
      endpoint $\TheEnd$ that pass through $\AltEuPoint$ is a
      subtree $\TheTree[\AltEuPoint]$ of $\TheTree[\TheRay]$.
      The intersection $\TheTree[\AltEuPoint]\intersect
      \TheTips[\TheIndex]$ is a finite set, on which the
      finite group $\StabOf[\TheImage]{\AltEuPoint}$ acts.
      By Proposition~\ref{quotient:fundamental.domain}, each
      $\StabOf[\TheImage]{\AltEuPoint}$-orbit has a representative
      in $\TheRepSet[\TheRay]$. As the cardinality
      \(
        \CardOf{
          \TheTips[\TheIndex]\intersect\TheTree[\AltEuPoint]
        }
      \)
      grows exponentially with the distance
      $\DistanceOf{\AltEuPoint}{\TheEuPoint[\TheIndex]}$,
      so does the size
      \(
        \CardOf{
          \StabOf[\TheImage]{\AltEuPoint}
        }
      \).
      As $\TheRatChamber[\TheIndex]$ is compact, the growth rate
      is uniformly bounded away from $\One$ for all rays $\TheRay$
      from the tip $\TheEuPoint[\TheIndex]$ into $\TheSector[\TheIndex]$.

      On the other hand, the number of chambers in $\TheEuBuild$
      intersecting $\TheSector[\TheIndex]$ at a point of distance
      $\TheDistance$ to $\TheEuPoint[\TheIndex]$ grows only
      polynomially with the distance $\TheDistance$.

      Finally, Lemma~\ref{quotient:stabilizer.index} shows that
      stabilizers of chambers of bounded distance have
      comparable sizes.

      It follows that the sum~(\ref{serre}) converges.
    \end{proof}
    The case that $\TheEuBuild$ is a tree is treated in
    \cite[Exercise~2\,a, page~110]{Serre:1980}.


  \noindent
  \parbox[t]{\textwidth}{%
  Kai-Uwe~Bux\\
  Fakult\"at f\"ur Mathematik\\
  Universit\"at Bielefeld\\
  Universit\"atsstra\ss{}e 25\\
  33501 Bielefeld\\
  Germany\\
  website: \texttt{www.kubux.net}
  }
  \bigskip

  \noindent
  \parbox[t]{\textwidth}{%
  Ralf~Gramlich\\
  Fachbereich Mathematik\\
  {\small TU} Darmstadt\\
  Schlo\ss{}gartenstra\ss{}e 7\\
  64289 Darmstadt\\
  Germany\\
  website: \texttt{www.mathematik.tu-darmstadt.de/$\sim$gramlich}\\[2mm]
  alternative address:\\
  Justus-Liebig-Universit\"at Gie\ss{}en\\
  Mathematisches Institut\\
  Arndtstra\ss e 2\\
  35392 Gie\ss{}en\\
  Germany}
  \bigskip

  \noindent
  \parbox[t]{\textwidth}{%
  Stefan~Witzel\\
  Fachbereich Mathematik\\
  {\small TU} Darmstadt\\
  Schlo\ss{}gartenstra\ss{}e 7\\
  64289 Darmstadt\\
  Germany\\
  website: \texttt{www.mathematik.tu-darmstadt.de/$\sim$switzel}
  }

\end{document}